\documentclass[11pt,a4paper,usenames,dvipsnames]{article}

\usepackage{url,amsmath,amssymb,latexsym,mathrsfs,comment,amsthm,enumerate,geometry,calc,ifthen,multicol,footnote}
\usepackage[noadjust]{cite}
\usepackage[colorlinks]{hyperref}

\usepackage[margin=10pt,font=small,labelfont=bf, labelsep=period]{caption}

\usepackage[x11names, svgnames, rgb,table]{xcolor}
\usepackage{hhline}
\usepackage{multirow}

\usepackage{enumitem}

\usepackage[utf8]{inputenc}
\usepackage[T1]{fontenc}

\usepackage{tikz}
\usetikzlibrary{decorations.markings}
\usetikzlibrary{matrix}
\usepgflibrary{snakes,arrows,shapes}

\makeatletter
\def\blfootnote{\gdef\@thefnmark{}\@footnotetext}
\makeatother

\newcounter{ncols}
\newcounter{incols}
\newenvironment{partn}[1]{
  \setcounter{ncols}{#1} \setcounter{incols}{\thencols - 1}\setlength{\arraycolsep}{1pt}
  \Bigl( \hspace{-1.5truemm}\scriptsize \renewcommand*{\arraystretch}{1}
    \begin{array}{@{\hskip 3pt} c *{\theincols}{|c} @{\hskip 3pt}  }
}{
     \end{array}
     \normalsize \hspace{-1.5truemm}\Bigr)\setlength{\arraycolsep}{6pt}
}

\usetikzlibrary{calc}

\allowdisplaybreaks

\numberwithin{equation}{section}

\newtheorem{thm}[equation]{Theorem}
\newtheorem{lemma}[equation]{Lemma}

\newtheorem{prop}[equation]{Proposition}

\theoremstyle{definition}

\newtheorem{rem}[equation]{Remark}

\newcommand{\nc}{\newcommand}
\nc{\rnc}{\renewcommand}

\let\oldproofname=\proofname
\rnc{\proofname}{\rm\bf{\oldproofname}}

\nc{\pf}{\begin{proof}}
\nc{\epf}{\end{proof}}

\nc\degrc{\deg_{\operatorname{rc}}}
\nc\opp{^{\operatorname{op}}}

\nc\JEnew[1]{\textcolor{blue}{#1}}
\nc\ka\kappa
\nc\wh\widehat
\nc\bn{{\bf n}}
\nc\bm{{\bf m}}
\nc\bk{{\bf k}}
\nc\br{{\bf r}}
\nc\bs{{\bf s}}
\nc\mt\mapsto
\nc\op\oplus
\nc\si\sigma
\nc\Si\Sigma
\rnc\implies{\ \Rightarrow\ }
\rnc\iff{\ \Leftrightarrow\ }
\nc\Iff{\ \ \Leftrightarrow\ \ }
\nc\IFF{\qquad \Leftrightarrow\qquad }
\nc\K{\mathscr K}
\nc\R{\mathscr R}
\rnc\L{\mathscr L}
\nc\J{\mathscr J}
\rnc\H{\mathscr H}
\nc\D{\mathscr D}
\nc\sm\setminus
\nc\nab\nabla
\nc\restr{{\restriction}}
\nc\mr\mathrel
\nc\ze\zeta
\nc\coker{\operatorname{coker}}
\nc\De\Delta
\rnc\P{\mathcal P}
\nc\C{\mathcal C}
\nc\RR{\mathcal R}
\nc\leqJ{\leq_{\J}}
\nc\leqL{\leq_{\L}}
\nc\leqR{\leq_{\R}}
\nc\geqJ{\geq_{\J}}
\nc\geqL{\geq_{\L}}
\nc\geqR{\geq_{\R}}
\nc\leqJU{\leq_{\J^U}}
\nc\leqLU{\leq_{\L^U}}
\nc\leqRU{\leq_{\R^U}}
\nc\geqJU{\geq_{\J^U}}
\nc\geqLU{\geq_{\L^U}}
\nc\geqRU{\geq_{\R^U}}
\nc\lam\lambda
\nc\codom{\operatorname{codom}}
\nc\es\varnothing
\nc\bP{{\bf P}}
\nc\B{\mathcal B}
\nc\PB{\P\B}
\nc\PP{\mathscr P\P}
\nc\M{\mathcal M}
\nc\TL{\T\!\mathcal L}
\rnc\th\theta
\nc\vt\vartheta
\nc\ds\displaystyle
\nc\bbD{\mathbb D}
\nc\bbM{\mathbb M}
\nc\G{\mathcal G}
\nc\Cong{\operatorname{Cong}}

\nc\pre\preceq
\nc\ol\bar

\rnc{\arraystretch}{1.2}

\nc\ord{\operatorname{ord}}
\nc\dom{\operatorname{dom}}

\nc\id{\operatorname{id}}
\nc\ba{{\bf a}}
\nc\bb{{\bf b}}
\nc\A{\mathcal A}
\rnc\S{\mathcal S}
\nc\Z{\mathbb Z}
\nc\bone{{\bf 1}}

\nc\ben{\begin{enumerate}[label=\textup{(\roman*)},leftmargin=7mm]}
\nc\bena{\begin{enumerate}[label=\textup{(\alph*)},leftmargin=12mm]}
\nc\een{\end{enumerate}}
\nc{\bit}{\begin{itemize}}
\nc{\eit}{\end{itemize}}
\nc\bmc{\begin{multicols}}
\nc\emc{\end{multicols}}
\nc{\firstpfitem}[1]{#1.}
\nc{\pfitem}[1]{\medskip \noindent #1.}
\nc{\pfcase}[1]{\medskip\noindent {\bf Case #1.}}
\nc\aftercases{\medskip\noindent}

\nc\T{\mathcal T}
\nc\PT{\P\T}

\nc{\al}{\alpha}
\nc{\be}{\beta}
\nc{\ga}{\gamma}
\nc{\de}{\delta}
\nc\ve\varepsilon
\nc\Om\Omega
\nc\Ga\Gamma

\nc{\rank}{\operatorname{rank}}

\nc{\set}[2]{\{ {#1} : {#2} \}} 
\nc{\bigset}[2]{\big\{ {#1} : {#2} \big\}} 
\nc{\sub}{\subseteq}
\nc{\la}{\langle}
\nc{\ra}{\rangle}
\nc{\pres}[2]{\la {#1} : {#2} \ra} 

\nc{\COMMA}{,\qquad}
\nc{\AND}{\qquad\text{and}\qquad}
\nc{\WHERE}{\qquad\text{where}\qquad}
\nc{\WHERe}{\quad\text{where}\quad}
\nc{\BY}{\qquad\text{by}\qquad}
\nc{\GIVENBY}{\qquad\text{given by}\qquad}
\nc{\COMMa}{,\quad}
\nc{\ANd}{\quad\text{and}\quad}
\nc{\ANDSIM}{\qquad\text{and similarly}\qquad}
\nc{\OR}{\qquad\text{or}\qquad}

\nc{\custpartn}[3]{{\lower1.4 ex\hbox{
\begin{tikzpicture}[scale=.3]
\foreach \x in {#1}
{ \uvert{\x}  }
\foreach \x in {#2}
{ \lvert{\x}  }
#3 \end{tikzpicture}
}}}
\nc{\uvert}[1]{\fill (#1,2)circle(.2);}
\rnc{\lvert}[1]{\fill (#1,0)circle(.2);}
\newcommand{\darcx}[3]{\draw(#1,0)arc(180:90:#3) (#1+#3,#3)--(#2-#3,#3) (#2-#3,#3) arc(90:0:#3);}
\newcommand{\uarcx}[3]{\draw(#1,2)arc(180:270:#3) (#1+#3,2-#3)--(#2-#3,2-#3) (#2-#3,2-#3) arc(270:360:#3);}
\nc{\darc}[2]{\darcx{#1}{#2}{.4}}
\nc{\uarc}[2]{\uarcx{#1}{#2}{.4}}
\nc\udotted[2]{\draw[dotted] (#1+.5,2)--(#2-.5,2);}
\nc\ddotted[2]{\draw[dotted] (#1+.5,0)--(#2-.5,0);}

\nc{\lv}[1]{\fill (#1,0)circle(.15);}
\nc{\lvs}[1]{{\foreach \x in {#1}{\lv{\x}}}}

\nc{\uv}[1]{\fill (#1,2)circle(.15);}
\nc{\uvs}[1]{{\foreach \x in {#1}{\uv{\x}}}}

\nc{\uvc}[2]{\fill[#2] (#1,2)circle(.15);}
\nc{\uvcs}[2]{{\foreach \x in {#1}{\uvc{\x}{#2}}}}

\nc{\stline}[2]{\draw(#1,2)--(#2,0);}
\nc{\hstline}[1]{\draw(#1,2)--(#1,.5);}
\nc{\stlines}[1]{{\foreach \x/\y in {#1}{\stline{\x}{\y} }}}
\nc{\hstlines}[1]{{\foreach \x/\y in {#1}{\hstline{\x}{\y} }}}

\nc\ur[3]{\fill[blue!20] (#1,2)--(#2,2)--(#2,2-#3)--(#1,2-#3); \draw[dotted] (#1,2)--(#2,2);}
\nc\lr[3]{\fill[blue!20] (#1,0)--(#2,0)--(#2,#3)--(#1,#3); \draw[dotted] (#1,0)--(#2,0);}

\begin{document}

\title{\vspace{-1.2cm}Minimum transformation representations of diagram monoids}

\date{}
\author{}

\maketitle

\vspace{-15mm}

\begin{center}
{\large 
Reinis Cirpons,%
\hspace{-.25em}\footnote{UFR Sciences et Techniques 2, rue de la Houssinière
BP 92208,
44322 Nantes Cedex 3,
France. {\it Email:} {\tt reinis.cirpons@inria.fr}.  This research was conducted while the first author was a doctoral student at the University of St Andrews.}
James East,%
\hspace{-.25em}\footnote{Centre for Research in Mathematics and Data Science, Western Sydney University, Locked Bag 1797, Penrith NSW 2751, Australia. {\it Email:} {\tt j.east@westernsydney.edu.au}.  Supported by ARC Future Fellowship FT190100632.}
James D.~Mitchell\footnote{Mathematical Institute, School of Mathematics and Statistics, University of St Andrews, St Andrews, Fife KY16~9SS, UK. {\it Email:} {\tt jdm3@st-andrews.ac.uk}.}\blfootnote{We thank Nik Ru\v skuc for some useful discussions.  We also thank the referee for their helpful suggestions.}
}
\end{center}

\vspace{1mm}

\begin{abstract}
We obtain formulae for the minimum transformation degrees of the most well-studied families of finite diagram monoids, including the partition, Brauer, Temperley--Lieb and Motzkin monoids.  For example, the partition monoid $\mathcal P_n$ has degree $1 + \frac{B(n+2)-B(n+1)+B(n)}2$ for $n\geq2$, where these are Bell numbers.  The proofs involve constructing explicit faithful representations of the minimum degree, many of which can be realised as (partial) actions on projections.


\end{abstract}

\tableofcontents

\section{Introduction}\label{sect:intro}

The study of permutation groups is as old as group theory itself, going back to the work of Galois on polynomials.  In a sense, group theory \emph{is} the study of permutation groups; indeed, as soon as Cayley gave the axioms of a group, he showed that every abstract group $G$ embeds in the symmetric group~$\S_G$.  In the standard embedding, each $g\in G$ is mapped to the right translation $\rho_g:G\to G$, acting by $a\rho_g = ag$ for all $a\in G$.  Thus, when $G$ is finite, it embeds in a finite symmetric group~$\S_n$.  The smallest such $n$ is called the \emph{(minimum permutation) degree} of~$G$.  This invariant is notoriously difficult to calculate, even for fairly simple groups, and can interact unintuitively with basic constructions such as (semi)direct products and quotients; see for example \cite{EP1988,EH2016,EST2010,HW2002,Johnson1971,KP2000,Saunders2010,Wright1975, Saunders2014,KP1989,BSS2017,OPU2024}.

The proof of Cayley's Theorem works without modification to show that any finite monoid~$M$ (and hence any finite semigroup) embeds in a finite full transformation monoid $\T_n$, which consists of all self-maps of $\{1,\ldots,n\}$.  The smallest such $n$ is called the \emph{(minimum transformation) degree} of $M$, and is denoted~$\deg(M)$.  
Much of the initial work on transformation degrees was undertaken by Easdown and Schein \cite{Easdown1987,Easdown1988,Easdown1992,Schein1988,Schein1992}; for some earlier work on transformation representations see \cite{Stoll1944,Tully1961,Slover1965,Hoehnke1963}. See also \cite{BT1991,EEMP2019,EEM2017,Holt2010,ACMT2023, LPRR2002,LPRR1998,LM1990,FP1997} for connections to computational algebra, and \cite{Clifford1942,Munn1957,Munn1964,Steinberg2016,AMSV2009,Steinberg2008,Steinberg2006,MS2023,Ponizovskii1956} for related studies on linear representations of semigroups and monoids.

It can again be difficult to calculate the degree of a given semigroup, even in very simple cases, and there are very few papers in the literature on this topic.  As illustrative examples, fix an integer $p\geq1$, and let~$L_p$ and~$R_p$ be  \emph{left-} and \emph{right-zero semigroups} of size~$p$, whose products obey the respective rules $xy=x$ or $xy=y$.  Easdown showed in \cite{Easdown1992} (see also \cite{GM2008}) that
\begin{align*}
\deg(L_p) &= \min\set{n\geq1}{r^{n-r}\geq p \text{ for some } 1\leq r\leq n},\\
\text{and}\qquad \deg(R_p) &= \min\set{n\geq1}{\textstyle{\prod}\lam\geq p \text{ for some integer partition $\lam\vdash n$}},
\end{align*}
%
and posed the problem of calculating the degree of an arbitrary \emph{rectangular band} $L_p\times R_q$.  
This problem was open for more than thirty years, and was eventually solved in~\cite{CEFMPQ2024}, where it was translated to a question concerning proper colourings of uniform hypergraphs, leading to:
\[
\deg(L_p\times R_q) = \min_{r\geq2}\be_r(p,q), \WHERE \be_r(p,q) = {\min}\bigset{n}{\pi_r(n-\lceil\log_rp\rceil)\geq q}.
\]
Here $\pi_r(m)$ is the maximum product $\prod\lam$ of a partition $\lam\vdash m$ with $r$ parts.  

The current article concerns the class of \emph{diagram monoids}, which consist of various kinds of set partitions that are depicted and multiplied diagrammatically.  These monoids, and related algebras and categories built from them, have origins and applications in many fields of mathematics and science.  See the seminal studies \cite{Brauer1937,BH2014,Jones1994_2,Martin1994,TL1971,MM2014}, and the surveys \cite{Martin2008,Koenig2008}; more background, and many more references, can be found in the introduction to \cite{EG2017}.

The diagram monoids we consider here are all submonoids of the partition monoid $\P_n$.  This monoid, which will be defined below, contains a copy of the full transformation monoid $\T_n$, as well as its \emph{opposite monoid} $\T_n\opp$.  Margolis and Steinberg have recently shown \cite{MS2023} that $\deg(\T_n\opp) = 2^n$, and this leads to the lower bound of $\deg(\P_n)\geq2^n$.  As far as we are aware, no upper bound for $\deg(\P_n)$ has been given in the literature, apart from $\deg(\P_n)\leq|\P_n|$, which comes directly from Cayley's Theorem.
However, since $\P_n$ is fundamental \cite{EMRT2018}, it follows from \cite[Theorem 5]{Imaoka1980} that~$\P_n$ can be faithfully represented in $\T_P\times\T_P\opp$, where $P = P(\P_n)$ is the set of \emph{projections} of~$\P_n$ (again, see below for the precise definitions).  Combining this with the result from \cite{MS2023} mentioned above, this leads to the upper bound of $\deg(\P_n) \leq |P| + 2^{|P|}$.  For example, taking $n=3$, and combining this last inequality with $\deg(\P_n)\geq2^n$ yields:
\[
8 \leq \deg(\P_3) \leq 4194\ 326.
\]
This upper bound is not an improvement on Cayley's Theorem, which gives ${\deg(\P_3) \leq |\P_3| = 203}$.  As we will see, the exact value of $\deg(\P_3)$ turns out to be $22$, which is a special case of a general formula we establish.  Specifically, we will show that
\begin{equation}\label{eq:degPn}
\deg(\P_n) = 1 + \frac{B(n+2) - B(n+1) + B(n)}2,
\end{equation}
where here $B(n)$ is the $n$th Bell number.  The `$1+$' in the above formula can be removed to obtain the minimum \emph{partial} transformation degree, which we denote by $\deg'(\P_n)$.

We also obtain analogous results for several other important families of diagram monoids---namely the Brauer monoid $\B_n$, the partial Brauer monoid $\PB_n$, the planar partition monoid~$\PP_n$, the Motzkin monoid $\M_n$, and the Temperley--Lieb monoid $\TL_n$---in terms of equally-famous number sequences, such as Catalan and Motzkin numbers.
These formulae (stated for partial transformation degrees) are summarised in Table \ref{tab:Q1}; calculated values are given in~Table~\ref{tab:Q2}.

The process of discovering the above formula for $\deg(\P_n)$ was somewhat circuitous, and involved several distinct steps.  Due to the lack of general theoretical techniques for calculating minimal degrees,\footnote{The aforementioned paper of Margolis and Steinberg \cite{MS2023} contains some very general results on minimal partial transformation representations of \emph{Rhodes semisimple} semigroups, i.e.~finite semigroups possessing faithful semisimple matrix representations.  As noted in \cite{EMRT2018}, the minimal ideals of our diagram monoids are non-trivial square bands, and for reasons explained in \cite{MS2023} this is an obstruction to Rhodes-semisimplicity.  Consequently, the results of \cite{MS2023} cannot be applied here, and it was necessary to develop tools for moving beyond the semisimple case.} we began with some computational experiments using the Semigroups package for~GAP~\cite{GAP,Semigroups}, as well as some recent advances described in \cite[Section 5.7]{ACMT2023}.  It has long been known~\cite{Tully1961} that any finite monoid $M$ acts on the quotient $M/\si$ for any right congruence $\si$.  The minimum degree of a \emph{faithful} `right congruence action' is denoted $\degrc(M)$, and this is an upper bound for~$\deg(M)$.  Moreover, the action of $M$ on $M/\si$ is faithful precisely when $\si$ contains no non-trivial two-sided congruence, and for any such faithful right congruence $\si$ we have
\[
\deg(M)\leq\degrc(M)\leq|M/\si|.
\]
The algorithms described in \cite{ACMT2023} allowed us to compute minimal faithful right congruences of~$\P_n$ for~${n\leq3}$, and obtain upper bounds for $\degrc(\P_n)$ for $4\leq n\leq7$.  It was not easy for us to comprehend the minimal congruences initially, but eventually we were able to understand a near-minimal one, and this led us to discover a general construction that will be outlined in Section \ref{subsect:si}.  In the case of $\P_n$, the congruence had two kinds of classes: one was a right ideal, and the others were all pieces of Green's $\L$-classes.  These are in turn indexed by the set $P$ of projections, which are special idempotents that are fixed by the involution of $\P_n$.  Investigating the movement of these projections, rather than that of the $\L$-classes themselves, led us to the key idea of \emph{partial actions} of so-called \emph{regular $*$-semigroups} \cite{NS1978} on their projections; see Section~\ref{subsect:Pact}.  We were able to show that the partial action of $\P_n$ on $P$ is faithful, using the classification of two-sided congruences in \cite{EMRT2018}.  The nature of the proof of this fact in turn led us to find a faithful sub-act $Q\sub P$, consisting of projections of suitably-low `rank', yielding the upper bound $\degrc(\P_n) \leq 1+|Q|$.  A combinatorial argument then shows that this reduces to the right-hand side of \eqref{eq:degPn}.
Showing that $1+|Q|$ is also a \emph{lower} bound for $\deg(\P_n)$ required an entirely different set of new techniques.  These are necessarily technical in nature, so we defer their further discussion until Section \ref{subsect:Plow}, where we will first use them.  It is worth noting, however, that the methods discussed so far apply uniformly to all of the monoids $\P_n$, $\PB_n$, $\PP_n$ and~$\M_n$.  They can be readily adapted to treat $\TL_n$, but must be significantly modified for $\B_n$.  
%
%
%
Although our results completely settle the question for all of the best-studied diagram monoids, we are hopeful that the new techniques we develop will prove useful in further studies of transformation degrees.

The paper is organised as follows.  We begin in Section \ref{sect:prelim} with preliminaries on semigroups and diagram monoids.  Section \ref{sect:T} contains general results on transformation representations, and their connections to actions and right congruences.  
Sections \ref{sect:P} and \ref{sect:TL} then apply this to the monoids $\P_n$, $\PB_n$, $\PP_n$, $\M_n$ and $\TL_n$.  The main results are Theorems~\ref{thm:M} and \ref{thm:TL}, which show that the degree of each such monoid is $1 + |Q|$, for a suitable set $Q$ of low-rank projections.  The Brauer monoid $\B_n$ is treated in Section~\ref{sect:B}, where we must use rather different methods; see Theorem~\ref{thm:B}.  
En route to proving these theorems, we construct explicit faithful actions/representations of the stated degree; see Theorems~\ref{thm:P},~\ref{thm:M1},~\ref{thm:Bnodd} and~\ref{thm:Bneven}.  Since every such action contains a global fixed point, 
the minimum \emph{partial} transformation degree is given by $\deg'(M) = \deg(M) - 1$.
Finally, in Section~\ref{sect:C} we give combinatorial formulae for~$|Q|$ for each the monoids~$\P_n$,~$\PB_n$,~$\PP_n$,~$\M_n$ and~$\TL_n$.  
As mentioned above, formulae and computed values are given in Tables~\ref{tab:Q1} and~\ref{tab:Q2}.



%

\begin{table}[ht]
\begin{center}
\begin{tabular}{lll}
\hline
Monoid $M$ & Validity & Minimum partial transformation degree $\deg'(M)$ \\
\hline
\rule[-6pt]{0pt}{20pt}
$\P_n$ & $n\geq2$ & $\tfrac{B(n+2)-B(n+1)+B(n)}2$  \\
\hline
\rule[-6pt]{0pt}{20pt}
$\PB_n$ & $n\geq2$ &  $\tfrac{I(n+2)}2$  \\
\hline
\rule[-6pt]{0pt}{18pt}
\multirow{2}{*}{$\B_n$}   & $n\geq3$ odd & $\tfrac{n+1}2\cdot n!!$  \\
\hhline{|~|-|-|}
\rule[-6pt]{0pt}{20pt}
 & $n\geq4$ even & $\tfrac{(n+4)(n+2)}8\cdot (n-1)!!$   \\
\hline
$\PP_n$ & $n\geq2$ &  $C(n+2)-2C(n+1)+C(n)$ \\
\hline
$\M_n$ & $n\geq2$ &  $M(n+2)-M(n+1)$  \\
\hline
\multirow{2}{*}{$\TL_n$} & $n=2k-1\geq3$ & $C(k+1)-C(k)$  \\
\hhline{|~|-|-|}
 & $n=2k\geq4$ & $C(k+2)-2C(k+1)+C(k)$   \\
\hline
\end{tabular}
\caption{Formulae for the minimum partial transformation degree, $\deg'(M)$, for diagram monoids~$M$, valid for the stated values of $n$.  For each such $M$ and $n$, the minimum transformation degree is equal to $\deg(M) = 1 + \deg'(M)$.  Here $B(n)$, $I(n)$, $C(n)$ and $M(n)$ are the~$n$th Bell, involution, Catalan and Motzkin numbers, and $m!!=m(m-2)(m-4)\cdots1$ for odd $m$.  See Theorems~\ref{thm:M},~\ref{thm:TL} and~\ref{thm:B}, and Propositions~\mbox{\ref{prop:QPn}--\ref{prop:QTLn}}.}
\label{tab:Q1}
\end{center}
\end{table}


\begin{table}[ht]
\begin{center}
\nc\colll{lightgray}
\scalebox{0.93}{
\begin{tabular}{lrrrrrrrrrrrl}
\hline
$n$ & 0 & 1 & 2 & 3 & 4 & 5 & 6 & 7 & 8 & 9 & 10 & OEIS\\
\hline
$\deg'(\P_n)$ & \textcolor{\colll}{1} & \textcolor{\colll}{1} & 6 & 21 & 83 & 363 & 1733 & 8942 & 49 484 & 291 871 & 1825 501  & A087649 \\
\hline
$\deg'(\PB_n)$ & \textcolor{\colll}{1} & \textcolor{\colll}{1} & 5 & 13 & 38 & 116 & 382 & 1310 & 4748 & 17 848 & 70 076  & A001475 \\
\hline
\multirow{2}{*}{$\deg'(\B_n)$}  & \textcolor{\colll}{1} &  & \textcolor{\colll}{2} &  & 18 &  & 150 &  & 1575 &  & 19 845 & $\frac13\times$A001194 \\
  &  & \textcolor{\colll}{1} &  & 6 &  & 45 &  & 420 &  & 4725 &  & A001879 \\
\hline
$\deg'(\PP_n)$ & \textcolor{\colll}{1} & \textcolor{\colll}{1} & 6 & 19 & 62 & 207 & 704 & 2431 & 8502 & 30 056 & 107 236  & A026012 \\
\hline
$\deg'(\M_n)$ & \textcolor{\colll}{1} & \textcolor{\colll}{1} & 5 & 12 & 30 & 76 & 196 & 512 & 1353 & 3610 & 9713  & A002026 \\
\hline
\multirow{2}{*}{$\deg'(\TL_n)$}  & \textcolor{\colll}{1} &  & \textcolor{\colll}{1} &  & 6 &  & 19 &  & 62 &  & 207  & A026012 \\
  &  & \textcolor{\colll}{1} &  & 3 &  & 9 &  & 28 &  & 90 &   & A000245 \\
\hline
\end{tabular}
}
\caption{Calculated values of $\deg'(M)$ for diagram monoids $M$, and their corresponding sequence numbers on the OEIS \cite{OEIS}.  Black entries are those for which the formulae in Table \ref{tab:Q1} hold.  For these entries we also have $\deg(M) = 1+\deg'(M)$.}
\label{tab:Q2}
\end{center}
\end{table}

\section{Preliminaries}\label{sect:prelim}

We begin by recalling the preliminary ideas and results we need concerning semigroups (Section~\ref{subsect:M}), regular $*$-semigroups (Section \ref{subsect:RSS}) and diagram monoids (Section \ref{subsect:D}).  For more basic background on semigroup theory, see for example \cite{Howie1995,CPbook}.

\subsection{Semigroups}\label{subsect:M}

Let $S$ be a semigroup, and let $S^1$ be the monoid completion of $S$.  So $S^1=S$ if $S$ is a monoid, or else $S=S\cup\{1\}$, where $1$ is a symbol not belonging to $S$, acting as an adjoined identity element.  Green's $\L$, $\R$ and $\J$ pre-orders and equivalences are defined, for $a,b\in S$, by
\begin{align*}
a\leqL b &\iff a\in S^1b ,&  a\mr\L b&\iff S^1a=S^1b,\\
a \leqR b &\iff a\in bS^1,& a \mr\R b&\iff aS^1= bS^1,\\
a \leqJ b &\iff a\in S^1bS^1,& a \mr\J b&\iff S^1aS^1= S^1bS^1.
\end{align*}
So, for example, $a\mr\L b$ holds when either $a=b$, or else $a=sb$ and $b=ta$ for some $s,t\in S$; similar comments apply to $\R$ and $\J$.
We also have the $\H$ and~$\D$ relations, defined by $\H=\L\cap\R$ and $\D=\L\vee\R$, where the latter is the join of $\L$ and $\R$ in the lattice of all equivalences of~$S$.  If $S$ is finite, then $\D=\J$.  The $\R$-class of an element $a\in S$ is denoted by $R_a$, and similarly for $\L$-classes, etc.  The set~$S/\R$ of all $\R$-classes is partially ordered by
\begin{equation}\label{eq:leqR}
R_a\leq R_b \iff a\leqR b \qquad\text{for $a,b\in S$.}
\end{equation}
A semigroup is \emph{stable} if
\begin{equation}\label{eq:stab}
sa \mr\J a \iff sa \mr\L a \AND as \mr\J a \iff as\mr\R a \qquad\text{for all $a,s\in S$.}
\end{equation}
Any finite semigroup is stable.

An equivalence relation $\si$ on a semigroup $S$ is a \emph{right congruence} if it is \emph{right-compatible}, meaning that
\[
(a,b)\in\si \implies (as,bs)\in\si \qquad\text{for all $a,b,s\in S$.}
\]
\emph{Left congruences} are defined symmetrically.  A \emph{(two-sided) congruence} is an equivalence that is both a left and right congruence.  For example, $\L$ is a right congruence, and $\R$ is a left congruence.  The \emph{trivial} and \emph{universal} congruences are respectively denoted
\[
\De_S = \set{(a,a)}{a\in S} \AND \nab_S = S\times S.
\]
For a set of pairs $\Si\sub S\times S$, we write $\Si^\sharp$ for the (two-sided) congruence of $S$ generated by~$\Si$, i.e.~the intersection of all congruences containing $\Si$.  When $\Si=\{(a,b)\}$ for some $a,b\in S$, we write $(a,b)^\sharp = \Si^\sharp$; such a congruence is called \emph{principal}.

A \emph{right ideal} of a semigroup $S$ is a subset $I\sub S$ such that $IS\sub I$.  \emph{Left ideals} and \emph{(two-sided) ideals} are defined analogously.  
Any left, right or two-sided ideal is a union of $\L$-, $\R$- or $\J$-classes, respectively.
The \emph{Rees right congruence} associated to a right ideal $I$ is defined by
\[
\RR_I = \nab_I \cup \De_S = \set{(x,y)\in S\times S}{x=y \text{ or } x,y\in I}.
\]
As special cases we have $\RR_\es = \De_S$ and $\RR_S = \nab_S$.

\subsection{Regular $*$-semigroups}\label{subsect:RSS}

A \emph{regular $*$-semigroup} is a semigroup $S$ with an additional unary operation $S\to S:a\mt a^*$ satisfying the identities
\[
a^{**}=a = aa^*a \AND (ab)^*=b^*a^* \qquad\text{for all $a,b\in S$.}
\]
These were introduced by Nordahl and Scheiblich in \cite{NS1978}, and the diagram monoids considered here are key examples.  The set of \emph{projections} of a regular $*$-semigroup $S$ is denoted
\[
P(S) = \set{p\in S}{p^2=p=p^*}.
\]
It is well known (see for example \cite{Petrich1985}) that $P(S) = \set{aa^*}{a\in S} = \set{a^*a}{a\in S}$, and that
\begin{equation}\label{eq:LR*}
a\mr\L b \iff a^*a=b^*b \AND a\mr\R b \iff aa^*=bb^* \qquad\text{for all $a,b\in S$.}
\end{equation}
Since every element of $S$ is $\R$-related to a unique projection (namely $a\mr\R aa^*$), with a similar statement for the $\L$ relation, it follows that any $\D$-class $D$ of $S$ contains $|D/\L| = |D/\R|$ projections, and that $|P(S)| = |S/\L| = |S/\R|$.

\subsection{Diagram monoids}\label{subsect:D}

Fix a non-negative integer $n\geq0$, and write $\bn=\{1,\ldots,n\}$ and $\bn'=\{1',\ldots,n'\}$.  The \emph{partition monoid}~$\P_n$ is the set of all set partitions of $\bn\cup\bn'$.  Such a partition is identified with any graph on vertex set $\bn\cup\bn'$ whose connected components are the blocks of the partition; the vertices from $\bn$ and~$\bn'$ are drawn on an upper and lower row, respectively, ordered $1<\cdots<n$ and~$1'<\cdots<n'$.

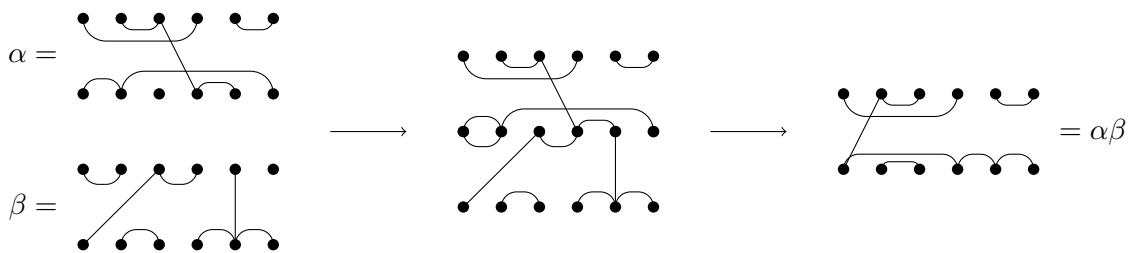
\begin{figure}[ht]
\begin{center}
\begin{tikzpicture}[scale=.5]

\begin{scope}[shift={(0,0)}]	
\uvs{1,...,6}
\lvs{1,...,6}
\uarcx14{.6}
\uarcx23{.3}
\uarcx56{.3}
\darc12
\darcx26{.6}
\darcx45{.3}
\stline34
\draw(0.6,1)node[left]{$\alpha=$};
\draw[->](7.5,-1)--(9.5,-1);
\end{scope}

\begin{scope}[shift={(0,-4)}]	
\uvs{1,...,6}
\lvs{1,...,6}
\uarc12
\uarc34
\darc45
\darc56
\darc23
\stline31
\stline55
\draw(0.6,1)node[left]{$\beta=$};
\end{scope}

\begin{scope}[shift={(10,-1)}]	
\uvs{1,...,6}
\lvs{1,...,6}
\uarcx14{.6}
\uarcx23{.3}
\uarcx56{.3}
\darc12
\darcx26{.6}
\darcx45{.3}
\stline34
\draw[->](7.5,0)--(9.5,0);
\end{scope}

\begin{scope}[shift={(10,-3)}]	
\uvs{1,...,6}
\lvs{1,...,6}
\uarc12
\uarc34
\darc45
\darc56
\stline31
\stline55
\darc23
\end{scope}

\begin{scope}[shift={(20,-2)}]	
\uvs{1,...,6}
\lvs{1,...,6}
\uarcx14{.6}
\uarcx23{.3}
\uarcx56{.3}
\darc14
\darc45
\darc56
\stline21
\darcx23{.2}
\draw(6.4,1)node[right]{$=\alpha\beta$};
\end{scope}

\end{tikzpicture}
\caption{Multiplication of partitions $\al,\be\in\P_6$, with the product graph $\Pi(\al,\be)$ in the middle.}
\label{fig:P6}
\end{center}
\end{figure}

The product of $\al,\be\in\P_n$ is defined as follows.  We first write $\bn''=\{1'',\ldots,n''\}$, and we let:
\bit
\item $\al^\vee$ be the graph on vertex set $\bn\cup\bn''$ obtained by changing each lower vertex $x'$ of $\al$ to $x''$,
\item $\be^\wedge$ be the graph on vertex set $\bn''\cup\bn'$ obtained by changing each upper vertex $x$ of $\be$ to~$x''$,
\item $\Pi(\al,\be)$ be the graph on vertex set $\bn\cup\bn''\cup\bn'$ whose edge set is the union of the edge sets of $\al^\vee$ and $\be^\wedge$.
\eit
We call $\Pi(\al,\be)$ the \emph{product graph} of $\al$ and $\be$, and we draw it with vertices $1''<\cdots<n''$ in a new middle row.  The product $\al\be\in\P_n$ is then the unique partition of $\bn\cup\bn'$ such that $x,y\in\bn\cup\bn'$ belong to the same block of $\al\be$ if and only if they belong to the same connected component of~$\Pi(\al,\be)$.  An example calculation is given in Figure \ref{fig:P6} in the case $n=6$.  The identity element of $\P_n$ is $\id_n = \custpartn{1,2,4}{1,2,4}{
\stline11
\stline22
\stline44
\udotted24
\ddotted24
}$, and the group of units is (isomorphic to) the symmetric group $\S_n$.

\begin{figure}[ht]
\begin{center}
\begin{tikzpicture}[scale=1.05]
\node[rounded corners,rectangle,draw,fill=orange!20] (P) at (3,8) {$\P_n$};
\node[rounded corners,rectangle,draw,fill=orange!20] (PB) at (3,6) {$\PB_n$};
\node[rounded corners,rectangle,draw,fill=orange!20] (B) at (3,4) {$\B_n$};
\node[rounded corners,rectangle,draw,fill=orange!20] (S) at (3,2) {$\S_n$};
\node[rounded corners,rectangle,draw,fill=orange!20] (PP) at (6,6) {$\PP_n$};
\node[rounded corners,rectangle,draw,fill=orange!20] (M) at (6,4) {$\M_n$};
\node[rounded corners,rectangle,draw,fill=orange!20] (TL) at (6,2) {$\TL_n$};
\node[rounded corners,rectangle,draw,fill=orange!20] (1) at (6,0) {$\id_n$};
\draw (P)--(PB)--(B)--(TL)--(M)--(PP)--(P) (PB)--(M) (B)--(S)--(1)--(TL);
\end{tikzpicture}
\qquad\qquad\qquad\qquad
\begin{tikzpicture}[scale=1]
\node[rounded corners,rectangle,draw,fill=orange!20] (P) at (3,8) {$\custpartn{1,2,3,4,5,6}{1,2,3,4,5,6}{\uarcx14{.6}\uarcx23{.3}\uarcx56{.3}\darc12\darcx26{.6}\darcx45{.3}\stline34}$};
\node[rounded corners,rectangle,draw,fill=orange!20] (PB) at (3,6) {$\custpartn{1,2,3,4,5,6}{1,2,3,4,5,6}{\uarcx13{.5}\uarc56\darc45\stline23}$};
\node[rounded corners,rectangle,draw,fill=orange!20] (B) at (3,4) {$\custpartn{1,2,3,4,5,6}{1,2,3,4,5,6}{\stline13\stline42\uarc23\uarc56\darcx16{.8}\darc45}$};
\node[rounded corners,rectangle,draw,fill=orange!20] (PP) at (6,6) {$\custpartn{1,2,3,4,5,6}{1,2,3,4,5,6}{\uarcx14{.6}\uarcx23{.2}\darc45\darc12\darc23\stline43\stline11\stline56}$};
\node[rounded corners,rectangle,draw,fill=orange!20] (M) at (6,4) {$\custpartn{1,2,3,4,5,6}{1,2,3,4,5,6}{\uarcx13{.5}\uarc56\darc45\stline43}$};
\node[rounded corners,rectangle,draw,fill=orange!20] (TL) at (6,2) {$\custpartn{1,2,3,4,5,6}{1,2,3,4,5,6}{\uarcx12{.4}\uarc45\darc34\darcx25{.8}\stlines{3/1,6/6}}$};
\node[rounded corners,rectangle,draw,fill=orange!20] (S) at (3,2) {$\custpartn{1,2,3,4,5,6}{1,2,3,4,5,6}{\stline13\stline21\stline32\stline44\stline56\stline65}$};
\node[rounded corners,rectangle,draw,fill=orange!20] (1) at (6,0) {$\custpartn{1,2,3,4,5,6}{1,2,3,4,5,6}{\stline11\stline22\stline33\stline44\stline55\stline66}$};
\draw (P)--(PB)--(B)--(TL)--(M)--(PP)--(P) (PB)--(M) (B)--(S)--(1)--(TL);
\end{tikzpicture}
\caption{Submonoids of $\P_n$ (left) and representative elements from each submonoid (right).}
\label{fig:submonoids}
\end{center}
\end{figure}
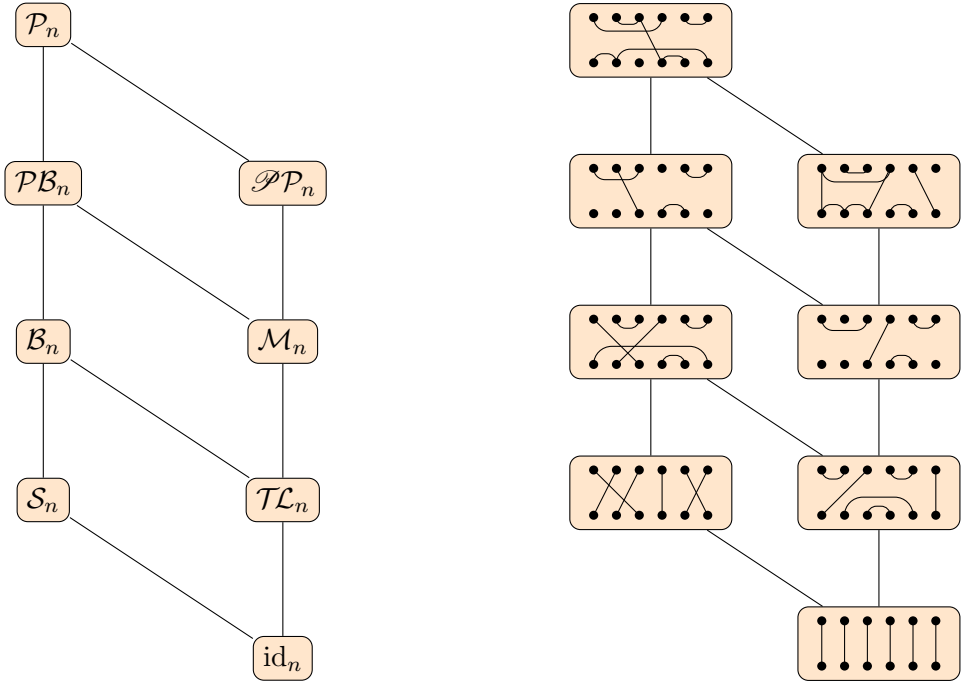

We now recall the definitions of the submonoids of $\P_n$ we will work with.  The first two are:
\bit
\item $\PB_n = \set{\al\in\P_n}{\text{each block of $\al$ has size $\leq2$}}$, the \emph{partial Brauer monoid},
\item $\B_n = \set{\al\in\P_n}{\text{each block of $\al$ has size $2$}}$, the \emph{Brauer monoid}.
\eit
A partition from $\P_n$ is \emph{planar} if it is represented by a graph whose edges are all contained within the rectangle spanned by the vertices, and have no crossings.  For example, $\be\in\P_6$ from Figure~\ref{fig:P6} is planar, but $\al$ is not.  We then have the further three submonoids of $\P_n$:
\bit
\item $\PP_n = \set{\al\in\P_n}{\al\text{ is planar}}$, the \emph{planar partition monoid},
\item $\M_n = \PP_n\cap\PB_n$, the \emph{Motzkin monoid},
\item $\TL_n = \PP_n\cap\B_n$, the \emph{Temperley--Lieb monoid} (sometimes called the \emph{Jones monoid}).
\eit
The containments among these monoids, along with sample elements, are shown in Figure \ref{fig:submonoids}.  It is well known (see for example \cite{Jones1994_2,HR2005}) that $\PP_n\cong\TL_{2n}$, via an isomorphism $\al\mt\widetilde\al$ illustrated in Figure \ref{fig:PPTL}.

\newcommand{\uarcxx}[4]{\draw[#4](#1,2)arc(180:270:#3) (#1+#3,2-#3)--(#2-#3,2-#3) (#2-#3,2-#3) arc(270:360:#3);}
\newcommand{\darcxx}[4]{\draw[#4](#1,0)arc(180:90:#3) (#1+#3,#3)--(#2-#3,#3) (#2-#3,#3) arc(90:0:#3);}
\newcommand{\stlinex}[3]{\draw[#3](#1,2)--(#2,0);}

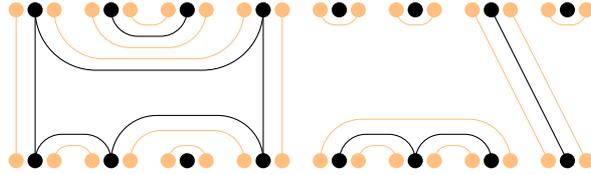
\begin{figure}[ht]
\begin{center}
\begin{tikzpicture}[scale=1]
\nc\coll{orange!50}
\uarcxx{1.25}{3.75}{.65}{\coll}
\uarcxx{1.75}{3.25}{.5}{\coll}
\uarcxx{2.25}{2.75}{.2}{\coll}
\uarcxx{4.75}{5.25}{.2}{\coll}
\uarcxx{5.75}{6.25}{.2}{\coll}
\uarcxx{7.75}{8.25}{.2}{\coll}
\uarcx14{.8}
\uarcx23{.35}
\darcxx{2.75}{3.25}{.2}{\coll}
\darcxx{1.25}{1.75}{.2}{\coll}
\darcxx{5.25}{5.75}{.2}{\coll}
\darcxx{6.25}{6.75}{.2}{\coll}
\darcxx{2.25}{3.75}{.4}{\coll}
\darcxx{4.75}{7.25}{.55}{\coll}
\stlinex{.75}{.75}{\coll}
\stlinex{4.25}{4.25}{\coll}
\stlinex{6.75}{7.75}{\coll}
\stlinex{7.25}{8.25}{\coll}
\darcx12{.35}
\darcx24{.6}
\darcx56{.35}
\darcx67{.35}
\stlines{1/1,4/4,7/8}
\foreach \x in {1,...,8} {
\fill (\x,2)circle(.1); \fill (\x,0)circle(.1);
\fill[\coll] (\x+.25,2)circle(.1);
\fill[\coll] (\x-.25,2)circle(.1);
\fill[\coll] (\x+.25,0)circle(.1);
\fill[\coll] (\x-.25,0)circle(.1);
}
\end{tikzpicture}
\caption{A planar partition $\al\in\PP_8$ (black), with its corresponding Temperley--Lieb element~${\widetilde\al\in\TL_{16}}$ (orange), illustrating the isomorphism $\PP_n\to\TL_{2n}$.}
\label{fig:PPTL}
\end{center}
\end{figure}

A non-empty subset $\es\not=X\sub\bn\cup\bn'$ is called:
\bit
\item an \emph{upper non-transversal} if $X\sub\bn$, i.e.~if $X$ contains only un-dashed vertices,
\item a \emph{lower non-transversal} if $X\sub\bn'$, i.e.~if $X$ contains only dashed vertices, or
\item a \emph{transversal} otherwise, i.e.~if $X$ contains both dashed and un-dashed vertices.
\eit
For $\al\in\P_n$ we define
\bit
\item $\dom(\al) = \set{x\in\bn}{x\text{ is contained in a transversal of }\al}$, the \emph{domain} of $\al$,
\item $\codom(\al) = \set{x\in\bn}{x'\text{ is contained in a transversal of }\al}$, the \emph{codomain} of $\al$,
\item $\ker(\al) = \set{(x,y)\in\bn\times\bn}{x\text{ and } y \text{ are contained in the same block of }\al}$, the \emph{kernel} of $\al$,
\item $\coker(\al) = \set{(x,y)\in\bn\times\bn}{x'\text{ and } y' \text{ are contained in the same block of }\al}$,  the \emph{cokernel} of $\al$.
\eit
We also define $\rank(\al)$, the \emph{rank} of $\al$, to be the number of transversals of $\al$.  For example, with $\al\in\P_6$ as in Figure \ref{fig:P6}, we have (using an obvious notation for equivalences):
\begin{align*}
\dom(\al) &= \{2,3\} , & \ker(\al) &= (1,4\mid2,3\mid5,6), \\
\codom(\al) &= \{4,5\} , & \coker(\al) &= (1,2,6\mid3\mid4,5), & \rank(\al)&=1.
\end{align*}
For $\al\in\P_n$ we use the tabular notation
\begin{equation}\label{eq:al}
\al = \begin{partn}{6} A_1&\cdots&A_r&C_1&\cdots&C_s\\ \hhline{~|~|~|-|-|-} B_1&\cdots&B_r&D_1&\cdots&D_t\end{partn}
\end{equation}
to indicate that $\al$ has transversals $A_i\cup B_i'$ ($i=1,\ldots,r$), upper non-transversals $C_i$ ($i=1,\ldots,s$) and lower non-transversals $D_i'$ ($i=1,\ldots,t$).  Note that one or two of $r,s,t$ could be $0$ in \eqref{eq:al}, but not all three (unless $n=0$).  
When we use this tabular notation, we always assume that the transversals are ordered so that ${\min(A_1)<\cdots<\min(A_r)}$.
When $\al\in\PP_n$ is planar, this in fact implies $A_1<\cdots<A_r$ and $B_1<\cdots<B_r$.  (Here for subsets $X,Y\sub\bn$ we write $X<Y$ to indicate that $x<y$ for all $x\in X$ and $y\in Y$.)

For $\al\in\P_n$ as in \eqref{eq:al}, we define
\[
\al^* = \begin{partn}{6} B_1&\cdots&B_r&D_1&\cdots&D_t\\ \hhline{~|~|~|-|-|-}A_1&\cdots&A_r&C_1&\cdots&C_s \end{partn},
\]
which is the partition obtained by interchanging dashed and un-dashed elements of $\bn\cup\bn'$.  Diagrammatically, $\al^*$ is obtained by reflecting any graph representing $\al$ in a horizontal axis.  For example, with $\al\in\P_6$ from Figure \ref{fig:P6} we have $\al^* = \custpartn{1,...,6}{1,...,6}{
\darcx14{.6}
\darcx23{.3}
\darcx56{.3}
\uarc12
\uarcx26{.6}
\uarcx45{.3}
\stline43}$.  It is easy to see that this gives~$\P_n$ the structure of a regular $*$-monoid, i.e.~that
\[
\al^{**}=\al=\al\al^*\al \AND (\al\be)^*=\be^*\al^* \qquad\text{for all $\al,\be\in\P_n$.}
\]
The projections of $\P_n$, i.e.~the elements $\ve\in\P_n$ satisfying $\ve^2=\ve=\ve^*$, have the form
\[
\ve = \begin{partn}{6} A_1&\cdots&A_r&C_1&\cdots&C_s\\ \hhline{~|~|~|-|-|-} A_1&\cdots&A_r&C_1&\cdots&C_s\end{partn}.
\]
Each of $\PB_n$, $\B_n$, $\PP_n$, $\M_n$ and $\TL_n$ is a regular $*$-submonoid of $\P_n$.  So too is $\S_n$, in which we have $\al^*=\al^{-1}$.

\section{Transformation monoids and representations}\label{sect:T}

For a set $X$, the \emph{full transformation monoid} $\T_X$ consists of all self-maps of $X$ under composition.  For $f\in\T_X$ and $x\in X$ we write $xf$ for the image of $x$ under $f$, and compose transformations left to right.  When $X=\{1,\ldots,n\}$ for a positive integer $n$, we write $\T_X=\T_n$.  

A \emph{transformation representation} of a semigroup $S$ is a homomorphism $S\to\T_X$ for some set~$X$.  When the representation is injective we say it is \emph{faithful}, and that $S$ \emph{embeds} in $\T_X$.  Cayley's Theorem states that any semigroup~$S$ embeds in $\T_{S^1}$; the proof is the observation that~$S$ acts faithfully on $S^1$ by right translations \cite[Theorem 1.1.2]{Howie1995}.  The \emph{minimum transformation degree} of a finite semigroup $S$ is defined to be
\[
\deg(S) = \min\set{n\geq1}{S\text{ embeds in }\T_n}.
\]

Our main goal in this paper is to compute this degree parameter for several families of diagram monoids.  In this section we establish some of the theoretical groundwork for doing this.  
We begin by recalling the connections between transformation representations and actions (Section \ref{subsect:act}), partial representations (Section \ref{subsect:pact}) and right congruences (Section \ref{subsect:RC}).  
In Section~\ref{subsect:si} we give a useful construction of right congruences from $\R$-classes.
Finally, we show in Section~\ref{subsect:Pact} how to build families of (partial) actions on projections of regular $*$-semigroups.

\subsection{Transformation representations and actions}\label{subsect:act}

In all that follows, it will be convenient to view transformation representations and degrees through the lens of actions.  For more detailed background, we refer to \cite{KKM2000}.

Recall that a \emph{(right) action} of a semigroup $S$ on a set $X$ is a map
\[
\mu:X\times S\to X \qquad\text{for which}\qquad \mu(\mu(x,a),b) = \mu(x,ab) \qquad\text{for all $x\in X$ and $a,b\in S$.}
\]
We write $\mu$ to the left of its arguments for readability.  It is standard to abbreviate $\mu(x,a)$ to~$xa$, in which case the defining property above says $(xa)b = x(ab)$.  However, we will only occasionally use such shorthand notation, and will always indicate when we are doing so.  One reason for this is that we will often have to deal with several actions at once.   Another reason is that a semigroup $S$ will often act on a subset $X\sub S$ in such a way that $\mu(x,a)$ is \emph{not} the product $xa$ of $x$ and $a$ in $S$.
In the literature,~$X$ is often referred to as an \emph{$S$-act} or an \emph{$S$-set}.
The \emph{degree} of the action $\mu$ is defined to be~$|X|$.  
We say~$\mu$ is a \emph{monoid action} if $S$ is a monoid with identity $1$ and $\mu(x,1)=x$ for all~$x$.
%

Given an action $\mu:X\times S\to X$, a subset $Y\sub X$ is called a \emph{sub-act} if $\mu(y,a)\in Y$ for all $y\in Y$ and $a\in S$.  In this case, $\mu$ restricts to an action $\mu\restr_Y:Y\times S\to Y$.
For an arbitrary subset~$Y$ of $X$, we define the sub-act
\[
\la Y\ra_\mu = Y\cup\set{\mu(y,a)}{y\in Y,\ a\in S},
\]
which we call the \emph{sub-act generated by $Y$}.  (Of course `$Y\cup{}$' can be deleted if $\mu$ is a monoid action.)  If $X = \la x\ra_\mu$ for some $x\in X$, we say that $\mu$ is \emph{monogenic}.  Monogenic acts are also called \emph{cyclic} or \emph{strictly cyclic} in the literature

Consider actions $\mu_1:X_1\times S\to X_1$ and $\mu_2:X_2\times S\to X_2$.  We say these are \emph{isomorphic}, and write $\mu_1\cong\mu_2$, if there is a bijection $\xi:X_1\to X_2$ such that
\[
(\mu_1(x,a))\xi = \mu_2(x\xi,a) \qquad\text{for all $x\in X_1$ and $a\in S$.}
\]
We say sub-acts $Y_1\sub X_1$ and $Y_2\sub X_2$ are isomorphic, and write $Y_1\cong Y_2$ if the restrictions are isomorphic: $\mu_1\restr_{Y_1} \cong \mu_2\restr_{Y_2}$.  

Consider again actions $\mu_1:X_1\times S\to X_1$ and $\mu_2:X_2\times S\to X_2$, where this time we assume that~$X_1$ and~$X_2$ are disjoint.  The (disjoint) union
\[
\mu_1\sqcup\mu_2:(X_1\sqcup X_2)\times S\to X_1\sqcup X_2
\]
is then an action.  More generally, suppose we have isomorphic sub-acts $Y_1\sub X_1$ and $Y_2\sub X_2$, as witnessed by the bijection $\xi:Y_1\to Y_2$.  We define $X = X_1\sqcup_\xi X_2$ to be the quotient of $X_1\sqcup X_2$ by the equivalence relation that equates $y$ and $y\xi$ for each $y\in Y_1$.  Writing $[x]\in X$ for the equivalence class of each $x\in X_1\sqcup X_2$, we then have a well-defined action
\begin{equation}\label{eq:pushout}
\mu :X\times S\to X \GIVENBY \mu([x],a) = \begin{cases}
[\mu_1(x,a)] &\text{if $x\in X_1$}\\
[\mu_2(x,a)] &\text{if $x\in X_2$.}
\end{cases}
\end{equation}
This action $\mu$ is called the \emph{push-out along $\xi$ of $\mu_1$ and $\mu_2$}, and is denoted by $\mu_1\sqcup_\xi\mu_2$.

Given an action $\mu:X\times S\to X$, one can define a transformation representation
\begin{align*}
&\phi_\mu:S\to\T_X:a\mt f_a , && \text{where} &xf_a &= \mu(x,a) &&\text{for $a\in S$ and $x\in X$.}
\intertext{Conversely, given a transformation representation $\phi:S\to\T_X:a\mt f_a$, one can define an action}
&\mu_\phi:X\times S\to X, && \text{where} &\mu_\phi(x,a) &= xf_a &&\text{for $a\in S$ and $x\in X$.}
\end{align*}
These constructions are mutually inverse.  

The \emph{kernel} of an action $\mu:X\times S\to X$, denoted $\ker(\mu)$, is defined to be the kernel of the corresponding representation $\phi_\mu:S\to\T_X$, so
\[
\ker(\mu) = \ker(\phi_\mu) = \set{(a,b)\in S\times S}{\mu(x,a) = \mu(x,b) \text{ for all }x\in X}.
\]
We say the action $\mu$ is \emph{faithful} if $\phi_\mu$ is faithful, which is equivalent to $\ker(\mu)=\De_S$.  Thus, we also have
\[
\deg(S) = \min\set{n\geq1}{S\text{ has a faithful action of degree }n}.
\]

\subsection{Partial actions and representations}\label{subsect:pact}

The \emph{partial transformation monoid} $\PT_X$ consists of all partial self-maps of $X$ under (relational) composition.  Let $-$ be a symbol not belonging to $X$, and write $X^- = X\cup\{-\}$.  Then $\PT_X$ is isomorphic to the submonoid of $\T_{X^-}$ consisting of all transformations that fix $-$.  In particular, we have embeddings $\T_n\hookrightarrow\PT_n\hookrightarrow\T_{n+1}$ for any $n$.

One can of course speak of partial transformation representations and degrees.  Denote the \emph{minimum partial transformation degree} of a finite semigroup $S$ by
\[
\deg'(S) = \min\set{n\geq1}{S\text{ embeds in }\PT_n}.
\]
It follows from the above-mentioned embeddings that $\deg'(S)\leq\deg(S)\leq\deg'(S)+1$, and that $\deg'(S) = \deg(S)-1$ holds precisely when there is a minimum-degree faithful transformation representation with a global fixed point.  This will be the case for every representation/action we construct in Sections \ref{sect:P}--\ref{sect:B}.

\subsection{Actions and right congruences}\label{subsect:RC}

One special family of actions (and hence transformation representations) come from right congruences, and these will be especially important in the current work.  To describe them, fix some right congruence~$\si$ on a semigroup $S$.  Write $[x]_\si$ for the $\si$-class of $x\in S$, and let $S/\si = \set{[x]_\si}{x\in S}$ be the set of all such classes.  We then have a well-defined action
\[
\mu_\si: (S/\si)\times S\to S/\si \GIVENBY \mu_\si([x]_\si,a) = [xa]_\si \qquad\text{for $a,x\in S$.}
\]
When $S$ is a monoid, $\mu_\si$ is generated by $[1]_\si$, so is monogenic.  Conversely, we have the following, which is a special case of \cite[Proposition 1.1]{Tully1961}.  We say an action $\mu$ is a \emph{right congruence action} if $\mu\cong\mu_\si$ for some right congruence $\si$.  

\begin{lemma}\label{lem:mono}
Any monogenic monoid action is a right congruence action.
\end{lemma}

The \emph{minimum right congruence degree} of $S$ is defined by
\[
\degrc(S) = \min\set{n\geq1}{S\text{ has a faithful right congruence action of degree }n}.
\]
Of course we have $\deg(S) \leq \degrc(S)$.

The next result is a special case of \cite[Proposition 1.2]{Tully1961}, formulated in a way that is convenient for our purposes.  

\begin{prop}\label{prop:phi}
If $\si$ is a right congruence of a monoid $S$, then $\ker(\mu_\si)$ is the largest two-sided congruence of~$S$ contained in~$\si$.  Consequently,~$\mu_\si$ is faithful if and only if $\si$ contains no non-trivial two-sided congruences.
\end{prop}

We say that an action $\mu:X\times S\to X$ \emph{separates} a pair $(a,b)\in S\times S$ if $\mu(x,a)\not=\mu(x,b)$ for some $x\in X$.  This is equivalent to the transformations $f_a,f_b\in\T_X$ under the associated transformation representation $\phi_\mu:S\to\T_X$ being different.

A two-sided congruence $\si$ of $S$ is \emph{minimal} if $\si\not=\De_S$, and the only congruence properly contained in $\si$ is $\De_S$.  Such a congruence is necessarily principal.  If $S$ is finite, then every non-trivial congruence contains a minimal congruence.

\begin{lemma}\label{lem:act1}
Let $S$ be a finite monoid, and let $\Ga\sub S\times S$ be such that $\set{(a,b)^\sharp}{(a,b)\in\Ga}$ consists of all the minimal two-sided congruences of $S$.  Then a semigroup action $\mu:X\times S\to X$ is faithful if and only if it separates each pair from $\Ga$.
\end{lemma}

\pf
If a pair $(a,b)\in\Ga$ was not separated by $\mu$, then we would have $\mu(x,a)=\mu(x,b)$ for all $x\in X$, meaning that $(a,b)\in\ker(\mu)$.  But then $\ker(\mu)\not=\De_S$, so $\mu$ is not faithful.

Conversely, suppose $\mu$ is not faithful, so that $\ker(\mu)\not=\De_S$.  We can then fix some minimal congruence $\si\sub\ker(\mu)$, and by assumption we have $\si=(a,b)^\sharp$ for some $(a,b)\in\Ga$.  Since $(a,b)\in\si\sub\ker(\mu)$ we have $\mu(x,a)=\mu(x,b)$ for all $x\in X$, so that $\mu$ does not separate $(a,b)$.
\epf

We say an equivalence relation $\si$ on a set $A$ \emph{separates} a subset $B\sub A$ if $(a,b)\not\in\si$ for distinct~${a,b\in B}$.

\begin{lemma}\label{lem:act3}
Let $S$ be a semigroup, and suppose $a,b\in S$ and $C\sub S$ are such that any right congruence on $S$ separating $\{a,b\}$ also separates $C$.  
Also let $\mu:X\times S\to X$ be an action with $\mu(x,a)\not=\mu(x,b)$ for some (fixed) $x\in X$.  Then the map $C \to X:c\mt\mu(x,c)$ is injective.
\end{lemma}

\pf
It is easy to check that the relation $\si = \set{(s,t)\in S\times S}{\mu(x,s)=\mu(x,t)}$ is a right congruence.  The assumption $\mu(x,a)\not=\mu(x,b)$ says that~$\si$ separates $\{a,b\}$, so it follows that~$\si$ separates $C$.  But this says precisely that the stated map is injective.
\epf

\subsection{Right congruences from $\R$-classes}\label{subsect:si}

For some of our later applications, we show how to build congruences on a semigroup starting from a specified $\R$-class.  These will have the form $\RR_I\vee\L^U$ for a carefully chosen right ideal~$I$ and subsemigroup $U$ of $S$.  Here $\RR_I$ is a Rees right congruence, and $\L^U$ denotes \emph{Green's relative~$\L$ relation}, introduced in \cite{Wallace1963}, and defined for $a,b\in S$ by
\[
a \mr\L^U b \iff U^1a = U^1b.
\]
This is again a right congruence, and as special cases we have $\L^S = \L$ and $\L^\es = \De_S$.

Throughout this section, let $R$ be a fixed $\R$-class of a semigroup $S$.  Define the sets
\begin{equation}\label{eq:TKI}
T = \set{a\in S}{aR\sub R} \COMMA K = \set{a\in S}{R\leq R_a} \AND I = S\sm K = \set{a\in S}{R\not\leq R_a},
\end{equation}
where $\leq$ is the ordering on $\R$-classes in \eqref{eq:leqR}.  It is easy to check that~$T$ is a (possibly empty) subsemigroup of $S$, and that $I$ is a (possibly empty) right ideal.


\begin{lemma}\label{lem:TK}
If $z$ is an arbitrary element of $R$, then
\ben
\item \label{TK1} $K=\set{a\in S}{z\leqR a}$,
\item \label{TK2} $T = \set{a\in S}{az\mr\R z}$,
\item \label{TK3} $T = \set{a\in S}{az = z}$ if $S$ is stable and $R$ is $\H$-trivial.
\een
\end{lemma}

\pf
\firstpfitem{\ref{TK1}}  Since $R=R_z$, this follows immediately from the definitions.  

\pfitem{\ref{TK2}}  Let $a\in S$.  If $aR\sub R$, then $az\in aR\sub R = R_z$, which says $az\mr\R z$.

Conversely, suppose $az\mr\R z$, and let $x\in R$ be arbitrary, so that $x\mr\R z$.  Since $\R$ is a left congruence it follows that $ax\mr\R az\mr\R z$, which says $ax\in R$, showing that $aR\sub R$.

\pfitem{\ref{TK3}}  Suppose $S$ is stable and $R$ is $\H$-trivial.  Given part \ref{TK2}, it is enough to show that ${az \mr\R z \implies az=z}$.  So suppose $az\mr\R z$.  Since $\R\sub\J$ we have $az\mr\J z$, and it follows from stability (see \eqref{eq:stab}) that $az\mr\L z$, and then from $\H$-triviality that $az=z$.
\epf

In what follows, we typically use parts \ref{TK1} and \ref{TK2} of Lemma \ref{lem:TK} without explicit reference.  If $\si$ is an equivalence on a set $X$, then for any subset $Y\sub X$ we denote by $\si\restr_Y = \si\cap(Y\times Y)$ the restriction of $\si$ to $Y$.

\begin{lemma}\label{lem:LIK}
For any subsemigroup $U\sub T$, we have $\L^U \sub \nab_I\cup\nab_K$, and consequently ${\L^U = \L^U\restr_I \cup \L^U\restr_K}$.
\end{lemma}

\pf
Let $(a,b)\in\L^U$.  By symmetry, and since $\L^U$ is an equivalence, it suffices to show that $a\in K \implies b\in K$.  So suppose $a\in K$.  Also let $z\in R$, so that $z\leqR a$, which gives $z=as$ for some $s\in S^1$.  Since $(a,b)\in\L^U$ we have $b=ua$ for some $u\in U^1$.  Since $U\sub T$ we have $uz\mr\R z$.  It follows that $z \mr\R uz = uas = bs \leqR b$, which gives $z\leqR b$, i.e.~$b\in K$.
\epf

With $R$, $T$, $K$ and $I$ as above, and for any subsemigroup $U\sub T$, we define the right congruence
\[
\si = \RR_I \vee \L^U,
\]
where here $\RR_I = \nab_I\cup\De_S = \nab_I \cup \De_K$ is the Rees right congruence, and where $\vee$ denotes the join in the lattice of equivalences of $S$.  It follows from Lemma \ref{lem:LIK} that in fact
\begin{equation}\label{eq:si}
\si = \nab_I \cup \L^U\restr_K = \bigset{(a,b)}{a,b\in I \text{ or } [a,b\in K \text{ and } U^1a=U^1b]}.
\end{equation}
In particular, we have $S/\si = \{I\} \cup (K/\L^U)$.

\subsection[Regular $*$-semigroups and (partial) actions on projections]{\boldmath Regular $*$-semigroups and (partial) actions on projections}\label{subsect:Pact}

We now show how the projections of a regular $*$-semigroup can be used to define (partial) actions, and hence transformation representations.

Let $S$ be a regular $*$-semigroup, and let $P=P(S) = \set{p\in S}{p^2=p=p^*}$ be its set of projections.  For $p\in P$ and $a\in S$ we write
\[
p^a = a^*pa = (pa)^*pa \in P.
\]
Since $(p^a)^b = p^{ab}$ for all $p\in P$ and $a,b\in S$, this defines an action
\begin{equation}\label{eq:act}
P\times S\to P:(p,a)\mt p^a.
\end{equation}
Now suppose $Q\sub P$ is a set of projections that is closed under the action \eqref{eq:act}, meaning that $p^a\in Q$ for all $p\in Q$ and $a\in S$.  Also suppose $\pre$ is a left-compatible pre-order on $S$ containing~$\leqR$, meaning that:
\bit
\item $\pre$ is reflexive and transitive,
\item $a\pre b \implies sa\pre sb$ for all $a,b,s\in S$, and
\item $as\pre a$ for all $a,s\in S$.
\eit
Let ${\approx} = {\pre} \cap {\succeq}$ be the equivalence on $S$ induced by $\pre$, so that $a\approx b$ if and only if $a\pre b$ and~${b\pre a}$.  Note that $\approx$ is a left congruence containing $\R$.  Let $-$ be a symbol not belonging to~$P$, let $Q^- = Q\cup\{-\}$, and define
\begin{equation}\label{eq:muQ}
\mu:Q^-\times S\to Q^- \BY \mu(p,a)=
\begin{cases}
p^a &\text{if $p\in Q$ and $p\approx pa$}\\
- &\text{otherwise.}
\end{cases}
\end{equation}
So in particular $\mu(-,a)=-$ for all $a\in S$.

\begin{prop}\label{prop:muQ}
If $S$ is a regular $*$-semigroup, and if $Q$, $\pre$ and $\approx$ are as above, then \eqref{eq:muQ} determines an action $\mu:Q^-\times S\to Q^-$.  If $S$ is a monoid, then $\mu$ is a monoid action.
\end{prop}

\pf
For the first assertion (the second is clear), we must show that
\[
 \mu(p,ab) = \mu(\mu(p,a),b) \qquad\text{for all $p\in Q^-$ and $a,b\in S$.}
\]
This is clear if $p=-$, so for the rest of the proof we assume that $p\in Q$.  Given that \eqref{eq:act} determines an action of $S$ on $Q$, it is in fact enough to show that
\[
\mu(p,ab) \in Q \IFF \mu(\mu(p,a),b) \in Q.
\]
Following the definitions, this amounts to showing that
\[
p\approx pab \IFF [p\approx pa \text{ and } p^a\approx p^ab]. 
\]
For the forward implication, suppose $p\approx pab$.  Since ${\leqR}\sub{\pre}$, we then have
\[
p \approx pab \pre pa \pre p,
\]
so that $p\approx pa\approx pab$.  From $pa\approx pab$ we have $a^*pa\approx a^*pab$ (as $\approx$ is left-compatible), i.e.~$p^a \approx p^ab$.  

Conversely, suppose $p\approx pa$ and $p^a\approx p^ab$.  From the latter (and left-compatibility) we have $pa\cdot p^a \approx pa\cdot p^ab$.  But
\[
pa\cdot p^a = pa \cdot a^*pa = pa\cdot (pa)^* \cdot pa = pa,
\]
so the previous conclusion says $pa\approx pab$.  Combined with $p\approx pa$, it follows that $p\approx pab$.
\epf

\begin{rem}\label{rem:leqR}
\ben
\item \label{leqR1} One obvious choice for the pre-order $\pre$ is $\leqR$ itself, in which case $\approx$ is $\R$.  Also taking the obvious $Q=P$, this leads to the action 
\begin{equation}\label{eq:RLM}
P^-\times S\to P^- \GIVENBY (p,a)\mt \begin{cases}
p^a &\text{if $p\mr\R pa$}\\
- &\text{otherwise.}
\end{cases}
\end{equation}
It turns out that this is equivalent to a well-known construction in the \emph{semilocal theory} of finite semigroups; see \cite[Chapters 7 and 8]{Arbib1968}, \cite[Section 4.6]{RS2009} and \cite[Appendix A]{MS2023}.  First consider the action of a \emph{finite} regular $*$-semigroup $S$ on $P=P(S)$ given in~\eqref{eq:act}.  We have already observed that the set~${S/{\L}}$ of $\L$-classes of $S$ is in one-one correspondence with $P$, via $p\mt L_p$.  Since $pa \mr\L a^*pa = p^a$ for all $p\in P$ and $a\in S$ (as $(pa)^*pa = p^a$; cf.~\eqref{eq:LR*}), it follows that this action is equivalent to that of $S$ on~${S/{\L}}$ given by $(L_p,a)\mt L_{pa}$.  For any $\D$-class $D$ of $S$, this restricts to a partial action of $S$ on~${D/{\L}}$ given by
\[
(L_p,a)\mt \begin{cases}
L_{pa} = L_{p^a} &\text{if $pa\in D$}\\
- &\text{otherwise.}
\end{cases}
\]
This is the \emph{Right Letter Mapping (RLM)} action associated to $D$.
Since $S$ is stable (as it is finite), the condition $pa\in D(=D_p)$ is equivalent to $pa\mr\R p$, and hence this RLM action of $S$ on~${D/{\L}}$ is equivalent to that of $S$ on $P(D)^-$ given by
\[
(p,a)\mt \begin{cases}
p^a &\text{if $p\mr\R pa$}\\
- &\text{otherwise.}
\end{cases}
\]
The union of these actions over all $\D$-classes $D$ of $S$ (amalgamating the common fixed point, $-$, of each) is precisely that given above in \eqref{eq:RLM}.

Finally, it is also worth noting that in the special case of an \emph{inverse} semigroup $S$, the partial transformation representation $S\to\PT_P$ associated to the action in \eqref{eq:RLM} is precisely Munn's representation of~$S$, used in his construction of \emph{fundamental} inverse semigroups \cite{Munn1970}.

\item \label{leqR2} Another obvious choice for $\pre$ is $\nab_S$, in which case $\approx$ is also equal to $\nab_S$.  Again taking $Q=P$, and in this case keeping in mind that $p\approx pa$ for all $p\in P$ and $a\in S$, the action~\eqref{eq:muQ} essentially reduces to \eqref{eq:act}, with the symbol $-$ acting as an adjoined fixed point.
\een
\end{rem}

\begin{rem}
Keeping Remark \ref{rem:leqR}\ref{leqR1} in mind, one could think of the results of this section as an extension of certain aspects of the semilocal theory, in the special case of regular $*$-semigroups.  It would be interesting to consider this in more generality.  As noted in the aforementioned remark, one could replace the conjugation action on projections by the standard action on $\L$-classes, $(L_x,a)\mt L_{xa}$ (which exists in any semigroup), and one would need a left-compatible pre-order $\pre$, as above.  In the case of partition monoids, a useful such pre-order comes from a so-called \emph{Ehresmann structure} \cite{Lawson1991,EG2021}, as we note in Section \ref{subsect:Pup}, and this could be a natural starting point for more general investigations.
\end{rem}

\begin{rem}
(Partial) actions on projections will be used in the next two sections to construct minimum-degree faithful representations/actions for many of our diagram monoids.  However, it is worth noting that the degree of a general regular $*$-semigroup need not have anything to do with the number of its projections.  For example, a billion-by-billion rectangular band has a billion projections, but degree $64$ \cite{CEFMPQ2024}.  
\end{rem}

\section{The partition, partial Brauer, planar partition and Motzkin monoids}\label{sect:P}

With the theoretical results of Section \ref{sect:T} now established, we now come to our main task: the calculation of the minimum transformation degrees of our diagram monoids.

Throughout this section, we let $M$ denote any of the partition monoid $\P_n$, the partial Brauer monoid $\PB_n$, the planar partition monoid $\PP_n$ or the Motzkin monoid $\M_n$.  The Temperley--Lieb monoid $\TL_n$ will be treated in Section \ref{sect:TL}, and the Brauer monoid $\B_n$ in Section \ref{sect:B}.  Our goal here is to prove the following:

\begin{thm}\label{thm:M}
If $n\geq2$, and if $M$ is any of $\P_n$, $\PB_n$, $\PP_n$ or $\M_n$, then
\[
\deg'(M) = |Q| \AND \deg(M) = \degrc(M) = 1+|Q|,
\]
where $Q = \set{\ve\in P(M)}{\rank(\ve) \leq 2}$.  Formulae for $|Q|$ can be found in Propositions \ref{prop:QPn}--\ref{prop:QMn}.
\end{thm}

To prove the theorem, we show that $1+|Q|$ is an upper bound for $\degrc(M)$ in Section \ref{subsect:Pup} (see Theorem~\ref{thm:P}), and a lower bound for $\deg(M)$ in Section \ref{subsect:Plow} (see Proposition \ref{prop:P}).  In fact, Theorem~\ref{thm:P} gives an explicit right congruence action of degree $1+|Q|$.  Since this action has a global fixed point it will follow (as explained in Section \ref{subsect:pact}) that $\deg'(M) = \deg(M) - 1$.

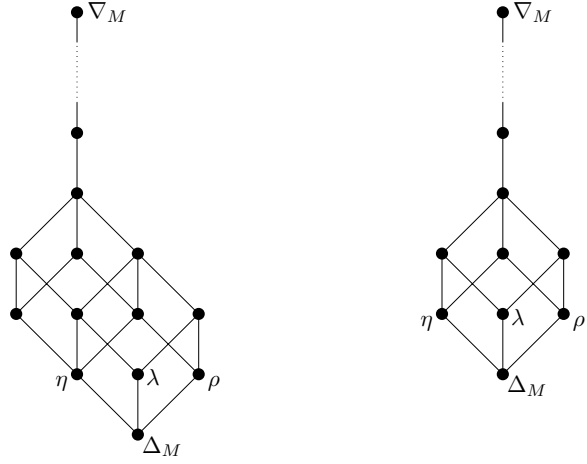
\begin{figure}[ht]
\begin{center}
\scalebox{0.8}{
\begin{tikzpicture}[scale=1]
\begin{scope}[shift={(0,0)}]	
\foreach \x/\y in {
0/0,
-1/1,0/1,1/1,
-1/2,0/2,1/2,
0/3,0/4,0/6
} {\fill (\x,\y)circle(.1);}
\draw (0,0)--(-1,1)--(-1,2)--(0,3)--(1,2)--(1,1)--(0,0)--(0,1)--(-1,2)
(-1,1)--(0,2)--(1,1)
(1,2)--(0,1)
(0,2)--(0,4.5)
(0,5.5)--(0,6)
;
\draw[dotted] (0,4.5)--(0,5.5);
\node () at (0.4,-0.2) {$\De_M$};
\node () at (-1.25,0.85) {$\eta$};
\node () at (0.25,0.95) {$\lam$};
\node () at (1.25,0.85) {$\rho$};
\node () at (0.5,6) {$\nab_M$};
\end{scope}
\begin{scope}[shift={(-7,0)}]	
\foreach \x/\y in {
0/0,
-1/1,0/1,1/1,
-1/2,0/2,1/2,
0/3,0/4,0/6,
1/-1,1/0,2/0,2/1
} {\fill (\x,\y)circle(.1);}
\draw (0,0)--(-1,1)--(-1,2)--(0,3)--(1,2)--(1,1)--(0,0)--(0,1)--(-1,2)
(-1,1)--(0,2)--(1,1)
(1,2)--(0,1)
(0,2)--(0,4.5)
(0,5.5)--(0,6)
(0,0)--(1,-1)--(1,0)--(0,1)
(1,1)--(2,0)--(2,1)--(1,2)
(1,-1)--(2,0)
(1,0)--(2,1)
;
\draw[dotted] (0,4.5)--(0,5.5);
\node () at (1+0.4,-0.2-1) {$\De_M$};
\node () at (1-1.25,0.85-1) {$\eta$};
\node () at (1+0.25,0.95-1) {$\lam$};
\node () at (1+1.25,0.85-1) {$\rho$};
\node () at (0.5,6) {$\nab_M$};
\end{scope}
\end{tikzpicture}
}
\caption{The congruence lattice $\Cong(M)$, for $M=\P_n$ or $\PB_n$ (left), and for $M=\PP_n$ or $\M_n$ (right).}
\label{fig:CongM}
\end{center}
\end{figure}

The proof of the theorem requires some understanding of the (two-sided) congruences of $M$ (being $\P_n$, $\PB_n$, $\PP_n$ or $\M_n$).  These were classified in \cite[Theorems~5.4,~6.1 and~7.3]{EMRT2018}.  We do not need the full classification here, but we do need to know that $M$ has precisely three \emph{minimal} (non-trivial) congruences:
\begin{align}
\nonumber \lam &= \De_M \cup \set{(\al,\be)\in M\times M}{\rank(\al)=\rank(\be)=0,\ \al\mr\L\be},\\
\nonumber \rho &= \De_M \cup \set{(\al,\be)\in M\times M}{\rank(\al)=\rank(\be)=0,\ \al\mr\R\be},\\
\label{eq:lamrhomu} \eta &= \De_M \cup \set{(\al,\be)\in M\times M}{\rank(\al),\rank(\be)\leq1,\ \wh\al=\wh\be}.
\end{align}
(These were respectively denoted $\lam_0$, $\rho_0$ and $\mu_1$ in \cite{EMRT2018}.)  The congruence $\eta$ involves the mapping
\[
\P_n\to\P_n:\al = \begin{partn}{6} A_1&\cdots&A_r&C_1&\cdots&C_s\\ \hhline{~|~|~|-|-|-} B_1&\cdots&B_r&D_1&\cdots&D_t\end{partn} \mt \wh\al = \begin{partn}{6} A_1&\cdots&A_r&C_1&\cdots&C_s\\ \hhline{-|-|-|-|-|-} B_1&\cdots&B_r&D_1&\cdots&D_t\end{partn}.
\]
Since the above congruences are minimal, they are generated by any non-trivial pair they contain.  Thus, for example, we have
\begin{equation}\label{eq:lrmgen}
\lam = (\ze,\al)^\sharp \COMMA
\rho = (\ze,\be)^\sharp \AND
\eta = (\ze,\ga)^\sharp ,
\end{equation}
in terms of the partitions
\begin{equation}\label{eq:zabc}
\ze = \custpartn{1,2,3,6}{1,2,3,6}{\udotted36 \ddotted36} \COMMA
\al = \custpartn{1,2,3,6}{1,2,3,6}{\uarc12\udotted36 \ddotted36} \COMMA
\be = \custpartn{1,2,3,6}{1,2,3,6}{\darc12\udotted36 \ddotted36} \AND
\ga = \custpartn{1,2,3,6}{1,2,3,6}{\stline11\udotted36 \ddotted36} .  
\end{equation}
Note here that $\be=\al^*$ (and $\al=\be^*$), while $\ze=\ze^*$ and $\ga=\ga^*$.
Figure \ref{fig:CongM} shows a Hasse diagram of the congruence lattice $\Cong(M)$, which is the set of all congruences of $M$, ordered by inclusion.  For more detailed diagrams see \cite[Figures~5 and~6]{EMRT2018}.

\subsection{Upper bound}\label{subsect:Pup}

Recall that $M$ denotes any of $\P_n$, $\PB_n$, $\PP_n$ or $\M_n$.  For the rest of this section we write
\[
P=P(M) = \set{\ve\in M}{\ve^2=\ve=\ve^*} \AND P_r = \set{\ve\in P}{\rank(\ve)=r} \qquad\text{for $0\leq r\leq n$.}
\]
We also define
\begin{equation}\label{eq:Q}
Q = Q(M) = P_0\cup P_1\cup P_2 = \set{\ve\in P}{\rank(\ve)\leq2}.
\end{equation}
Since $\rank(\ve^\al)=\rank(\al^*\ve\al) \leq\rank(\ve)$ for all $\ve\in P$ and $\al\in M$, it is clear that $Q$ is closed under the action~\eqref{eq:act}.

We now define the relation $\pre$ on $M$ by
\begin{equation}\label{eq:pre}
\al\pre\be \iff \ker(\al) \supseteq \ker(\be)  \qquad\text{for $\al,\be\in M$.}
\end{equation}
This relation played a crucial role in \cite{EG2021}, in connection with so-called \emph{Ehresmann structures} on~$\P_n$.  Of particular relevance to the current situation, it was shown in the proof of \cite[Lemma~4.6]{EG2021} that $\pre$ is left-compatible, and it is clearly transitive.  We obtain ${\leqR} \sub {\pre}$ from the identity $\ker(\al\be)\supseteq\ker(\al)$.  It then follows from Proposition \ref{prop:muQ} that we have an action
\begin{equation}\label{eq:muQM}
\mu:Q^-\times M\to Q^- \GIVENBY \mu(\ve,\al)=
\begin{cases}
\ve^\al &\text{if $\ve\in Q$ and $\ker(\ve) = \ker(\ve\al)$}\\
- &\text{otherwise.}
\end{cases}
\end{equation}

During this section, it will be convenient in some circumstances to omit singleton blocks in the tabular notation for partitions.  For example, we write $\al = \begin{partn}{6} A_1&\cdots&A_r&\multicolumn{3}{c}{} \\ \hhline{~|~|~|-|-|-} B_1&\cdots&B_r&D_1&\cdots&D_s\end{partn}$ to indicate that the upper non-transversals of $\al$ are all singletons.  
As concrete examples, the partitions~$\al$,~$\be$ and~$\ga$ in \eqref{eq:zabc} can be denoted as $\al = \begin{partn}{1} 1,2\ \\ \hhline{-} \ \end{partn}$, $\be = \begin{partn}{1} \ \\ \hhline{-} 1,2\ \end{partn}$ and $\ga = \begin{partn}{1} 1 \\ \hhline{~} 1 \end{partn}$.

In the remainder of this section, an important role will be played by the map
\begin{equation}\label{eq:olve}
P\to M: 
\ve = \begin{partn}{6} A_1&\cdots&A_r&C_1&\cdots&C_s\\ \hhline{~|~|~|-|-|-} A_1&\cdots&A_r&C_1&\cdots&C_s\end{partn}
\mt 
\ol\ve = \begin{partn}{6} 1&\cdots&r&\multicolumn{3}{c}{}\\ \hhline{~|~|~|-|-|-} A_1&\cdots&A_r&C_1&\cdots&C_s\end{partn}.
\end{equation}
Here as usual we assume that $\min(A_1)<\cdots<\min(A_r)$, which ensures that $\ol\ve$ is planar whenever~$\ve$ is.  Keeping in mind the convention regarding singleton blocks, we note that $\ker(\ol\ve)=\De_{\bn}$ for all $\ve\in P$.
Because of the identity $\ve = \ol\ve^*\ol\ve$, the map $\ve\mt\ol\ve$ is injective.
It is worth noting that this map need not map projections to projections.  In fact, one can check that $\ol\ve$ is a projection if and only if $\ol\ve=\ve$ is a partial identity, i.e.~has the form $\begin{partn}{3} a_1&\cdots&a_r\\ a_1&\cdots&a_r\end{partn}$ for some $a_1,\ldots,a_r\in\bn$.

\begin{thm}\label{thm:P}
Let $n\geq2$, let $M$ be any of $\P_n$, $\PB_n$, $\PP_n$ or $\M_n$, and let $\mu:Q^-\times M\to Q^-$ be the action in \eqref{eq:muQM}.  Then $\mu$ is faithful and monogenic, and consequently
\[
\deg(M) \leq \degrc(M) \leq 1+|Q|.
\]
\end{thm}

\pf
As explained above, the minimal congruences of $M$ are $(\ze,\al)^\sharp$, $(\ze,\be)^\sharp$ and $(\ze,\ga)^\sharp$, where $\ze,\al,\be,\ga\in M$ are as in \eqref{eq:zabc}.  Thus, by Lemma \ref{lem:act1} we can establish faithfulness of $\mu$ by showing that it separates each of $(\ze,\al)$, $(\ze,\be)$ and $(\ze,\ga)$.  For this, we define $\pi = \custpartn{1,2,3,6}{1,2,3,6}{\stline11\stline22\udotted36 \ddotted36}\in Q$, and one can check that
\[
\mu(\pi, \ze) = \ze \COMMA
\mu(\pi, \al) = {-} \COMMA
\mu(\pi, \be) =  \custpartn{1,2,3,6}{1,2,3,6}{\uarc12\darc12\udotted36 \ddotted36} \AND
\mu(\pi, \ga) = \ga,
\]
which are all distinct.

Given Lemma \ref{lem:mono}, it remains to check that $\mu$ is monogenic, and for this we claim that $Q^- = \la\pi\ra_\mu$, where $\pi\in Q$ is as above.  We have already observed that $-=\mu(\pi,\al)$, and for $\ve\in Q$ it is easy to check that $\ve = \mu(\pi,\ol\ve)$, where $\ol\ve$ is as in \eqref{eq:olve}.
\epf

\begin{rem}\label{rem:MPTP}
If instead we took ${\pre} = {\leqR}$ and $Q=P$, we would obtain the alternative action ${P^-\times M\to P^-}$, as in Remark \ref{rem:leqR}\ref{leqR1}, which we will here denote by $\mu'$.  It turns out that $\mu'$ is not faithful for any of the diagram monoids in Theorem~\ref{thm:P}.  Indeed, for any $\al\in M$ with $\rank(\al)=0$, and for  $\ve\in Q^-$, one can check that
\[
\mu'(\ve,\al) = \begin{cases}
\al^*\al &\text{if $\ve\in P_0$}\\
- &\text{otherwise.}
\end{cases}
\]
Combining this with \eqref{eq:LR*}, it follows that $(\al,\be)\in\ker(\mu')$ if $\al\mr\L\be$.  This is to say that $\lam \sub \ker(\mu')$, where $\lam$ is the congruence of $M$ in \eqref{eq:lamrhomu}.  
On the other hand, with $\be,\ga,\ze\in M$ as in~\eqref{eq:zabc}, we have
\[
\mu'(\ze,\ze) = \ze \COMMA \mu'(\ze,\be) = \custpartn{1,2,3,6}{1,2,3,6}{\uarc12\darc12\udotted36 \ddotted36} \COMMA
\mu'(\ga,\ze) = - \AND \mu'(\ga,\ga) = \ga,
\]
which shows that neither $(\ze,\be)$ nor $(\ze,\ga)$ belongs to $\ker(\mu')$.  Since $(\ze,\be)\in\rho$ and $(\ze,\ga)\in\eta$, it follows that $\rho\not\sub\ker(\mu')$ and $\eta\not\sub\ker(\mu')$.  Thus, keeping the structure of the congruence lattice $\Cong(M)$ in mind (see Figure \ref{fig:CongM}), it follows that in fact $\ker(\mu')=\lam$.
\end{rem}

\subsection{Lower bound}\label{subsect:Plow}

We continue to fix the monoid $M$, being one of $\P_n$, $\PB_n$, $\PP_n$ or $\M_n$, with $n\geq2$.  Our goal now is to show that $1+|Q|$ is a lower bound for $\deg(M)$, thus completing the proof of Theorem~\ref{thm:M}.

In what follows, we fix the elements $\ze,\al,\be,\ga\in M$ from \eqref{eq:zabc}, as well as the map ${P=P(M)\to M:\ve\mt\ol\ve}$ from \eqref{eq:olve}.  For $X\sub P$ we write $\ol X = \set{\ol\ve}{\ve\in X}$, noting that $|\ol X|=|X|$, as $\ve\mt\ol\ve$ is injective.  Also recall that $P_r = \set{\ve\in P}{\rank(\ve)=r}$ denotes the set of projections of rank $0\leq r\leq n$.

\begin{lemma}\label{lem:Msi}
Let $\si$ be a right congruence of $M$. 
\ben
\item \label{Msi1} If $\si$ separates $\{\ze,\al\}$, then it separates $\ol P_2$.
\item \label{Msi2} If $\si$ separates $\{\ze,\be\}$, then it separates $\ol P_0$.
\item \label{Msi3} If $\si$ separates $\{\ze,\ga\}$, then it separates $\ol P_1$.
\een
\end{lemma}

\pf
\firstpfitem{\ref{Msi1}}  
Aiming to prove the contrapositive, suppose $(\ol\ve_1,\ol\ve_2)\in\si$ for distinct $\ve_1,\ve_2\in P_2$, and write
\begin{equation}\label{eq:ve12a}
\ol\ve_1 = \begin{partn}{5} 1&2&\multicolumn{3}{c}{}\\ \hhline{~|~|-|-|-} A&B&C_1&\cdots&C_s\end{partn}
\ANd
\ol\ve_2 = \begin{partn}{5} 1&2&\multicolumn{3}{c}{}\\ \hhline{~|~|-|-|-} D&E&F_1&\cdots&F_t\end{partn}.
\end{equation}
We must show that $(\ze,\al)\in\si$.  

\pfcase1  Suppose first that $A\cup B\not= D\cup E$.  Without loss of generality, we can fix some $a\in A$ and $b\in B$ such that at least one of $a,b$ does not belong to $D\cup E$.  Then with $\th = \begin{partn}{1} a,b\ \\ \hhline{-} \ \end{partn}\in M$, we have $(\al,\ze) = (\ol\ve_1\th,\ol\ve_2\th)\in\si$.

\pfcase2  Next suppose $A\cup B=D\cup E$, but $A\not=D$ (and $B\not=E$).  (Since $\min(A)<\min(B)$ and $\min(D)<\min(E)$, it is impossible to have $A=E$ and $B=D$.)  Without loss of generality, we can fix some $a\in A$ and $b\in B$ such that $a,b\in D$ or $a,b\in E$.  Then with $\th = \begin{partn}{1} a,b\ \\ \hhline{-} \ \end{partn}\in M$, we again have $(\al,\ze) = (\ol\ve_1\th,\ol\ve_2\th)\in\si$.

\pfcase3  Finally, suppose $A=D$ and $B=E$.  Since $\ve_1\not=\ve_2$, we can assume without loss of generality that there exist $x,y\in C_1$ such that $x\in F_1$ and $y\in F_2$.  Also let $a=\min(A)$ and $b=\min(B)$, noting that $a<b$.  Then with $\th = \begin{partn}{2} a,x & b,y \\ \hhline{-|-} \multicolumn{2}{c}{} \end{partn}$ or $\begin{partn}{2} a,y & b,x \\ \hhline{-|-} \multicolumn{2}{c}{} \end{partn}$, whichever is planar (and hence belongs to $M$), we again have $(\al,\ze) = (\ol\ve_1\th,\ol\ve_2\th)\in\si$.  

\pfitem{\ref{Msi2}}  Suppose $(\ol\ve_1,\ol\ve_2)\in\si$ for distinct $\ve_1,\ve_2\in P_0$; this time we must show that $(\ze,\be)\in\si$.  Without loss of generality, we can fix some $(x,y)\in\coker(\ve_1)\sm\coker(\ve_2)$, say with $x<y$.  Then with $\th = \begin{partn}{2} x & y \\ \hhline{~|~} 1&2 \end{partn}\in M$ we have $(\be,\ze) = (\ol\ve_1\th,\ol\ve_2\th) \in \si$.

\pfitem{\ref{Msi3}}  
Suppose $(\ol\ve_1,\ol\ve_2)\in\si$ for distinct $\ve_1,\ve_2\in P_1$, and write
\begin{equation}\label{eq:ve12b}
\ol\ve_1 = \begin{partn}{4} 1&\multicolumn{3}{c}{}\\ \hhline{~|-|-|-} A&B_1&\cdots&B_s\end{partn}
\ANd
\ol\ve_2 = \begin{partn}{4} 1&\multicolumn{3}{c}{}\\ \hhline{~|-|-|-} C&D_1&\cdots&D_t\end{partn}.
\end{equation}

\pfcase1  Suppose first that $A\not=C$.  Without loss of generality, we can fix some $a\in A\sm C$.  Then with $\th = \binom a1\in M$ we have $(\ga,\ze) = (\ol\ve_1\th,\ol\ve_2\th) \in \si$.

\pfcase2  Now suppose $A=C$, and fix some $a\in A$.  Without loss of generality, we can also fix some ${(x,y)\in\coker(\ve_1)\sm\coker(\ve_2)}$.  Then with $\th = \begin{partn}{2} x&a,y \\ \hhline{~|-} 1&\ \end{partn}$ or $\begin{partn}{2} y&a,x \\ \hhline{~|-} 1&\ \end{partn}$, whichever is planar, we have $(\ga,\ze) = (\ol\ve_1\th,\ol\ve_2\th) \in \si$.
\epf

Here then is the last part of the proof of Theorem \ref{thm:M}:

\begin{prop}\label{prop:P}
If $n\geq2$, and if $M$ is any of $\P_n$, $\PB_n$, $\PP_n$ or $\M_n$, then
\[
\deg(M) \geq 1+|Q|.
\]
\end{prop}

\pf
Let $\mu:X\times M\to X$ be a faithful action.  We prove the result by showing that ${|X|\geq1+|Q|}$.  Throughout the proof it will be convenient to write $x\de = \mu(x,\de)$ for $x\in X$ and $\de\in M$.

Since $\mu$ is faithful, it follows from Lemma \ref{lem:act1} and \eqref{eq:lrmgen} that it separates the pairs $(\ze,\al)$, $(\ze,\be)$ and~$(\ze,\ga)$.  Thus, we can fix elements ${x_0,x_1,x_2\in X}$ for which
\begin{equation}\label{eq:x0x1x2}
x_0\ze\not=x_0\be \COMMA x_1\ze\not=x_1\ga \AND x_2\ze\not=x_2\al.
\end{equation}
We note that $x_0,x_1,x_2$ need not be distinct.  We also define the sets
\[
Y_i = x_i\ol P_i = \set{x_i\ol\ve}{\ve\in P_i} \sub X \qquad\text{for $i=0,1,2$.}
\]
By Lemmas \ref{lem:act3} and \ref{lem:Msi}, each map $P_i\to Y_i:\ve\mt x_i\ol\ve$ is a bijection, and hence
\begin{equation}\label{eq:YiPi}
|Y_i|=|P_i| \qquad\text{for $i=0,1,2$.}
\end{equation}
Our strategy for obtaining $|X|\geq1+|Q|$ involves showing that
\[
|Y_0\cup Y_1\cup Y_2| = |Q| = |P_0|+|P_1|+|P_2|,
\]
and that $X\sm(Y_0\cup Y_1\cup Y_2)$ is non-empty.

We begin by claiming that
\begin{equation}\label{eq:Y1Y2ze}
(Y_1\cup Y_2) \cap \set{x\ze}{x\in X} = \es.
\end{equation}
To prove this, suppose to the contrary that $x_i\ol\ve = x\ze$ for some $\ve\in P_i$ and $x\in X$, where $i=1$~or~$2$.  Since $\ze = \ze\ol\ve^*$, it then follows that 
\[
x\ze = x\ze\ol\ve^* = x_i\ol\ve\ol\ve^* = \begin{cases}
x_1\ga &\text{if $i=1$}\\
x_2\pi &\text{if $i=2$,}
\end{cases}
\]
where again $\pi = \custpartn{1,2,3,6}{1,2,3,6}{\stline11\stline22\udotted36 \ddotted36}$.  The $i=1$ and $i=2$ cases lead respectively to
\[
x_1\ze = x_1\ga\ze = x\ze\ze = x\ze = x_1\ga
\OR
x_2\ze = x_2\pi\ze = x\ze\ze = x\ze = x\ze\al = x_2\pi\al = x_2\al,
\]
both contradicting \eqref{eq:x0x1x2}.

Now that we have proved \eqref{eq:Y1Y2ze}, our next claim is that
\begin{equation}\label{eq:Y1Y2}
Y_1\cap Y_2 = \es.
\end{equation}
To prove this, suppose to the contrary that $x_1\ol\ve_1 = x_2\ol\ve_2$ for some $\ve_1\in P_1$ and $\ve_2\in P_2$.  By considering the form of $\ol\ve_1$ and $\ol\ve_2$, one can check that $\ol\ve_1 \ol\ve_2^*$ is equal to one of
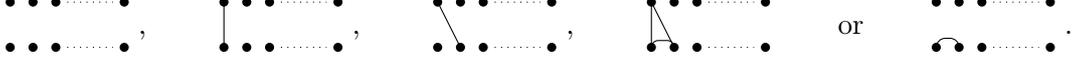
\[
\custpartn{1,2,3,6}{1,2,3,6}{\udotted36 \ddotted36} \COMMA 
\custpartn{1,2,3,6}{1,2,3,6}{\stline11\udotted36 \ddotted36} \COMMA 
\custpartn{1,2,3,6}{1,2,3,6}{\stline12\udotted36 \ddotted36} \COMMA 
\custpartn{1,2,3,6}{1,2,3,6}{\stline11\stline12\darcx12{0.3}\udotted36 \ddotted36} \OR
\custpartn{1,2,3,6}{1,2,3,6}{\darc12\udotted36 \ddotted36}.
\]
In any case, it follows that $\ol\ve_1 \ol\ve_2^*\al = \ze$.  On the other hand, we have $\ol\ve_2 \ol\ve_2^*\al = \pi\al = \al$.  Combining these, we obtain
\[
x_1\ze = (x_1\ol\ve_1) \ol\ve_2^*\al = (x_2\ol\ve_2) \ol\ve_2^*\al = x_2\al.
\]
From $\ol\ve_1\ze=\ol\ve_2\ze=\ze$, we also have
\[
x_1\ze = (x_1\ol\ve_1)\ze = (x_2\ol\ve_2)\ze = x_2\ze.
\]
Combining the last two conclusions, it follows that $x_2\al = x_2\ze$, which contradicts \eqref{eq:x0x1x2}.  

Now that we have proved \eqref{eq:Y1Y2}, our next claim is that
\begin{equation}\label{eq:Y1Y2Y0}
(Y_1\cup Y_2) \cap Y_0 = \es.
\end{equation}
To prove this, suppose to the contrary that $x_i\ol\ve_i = x_0\ol\ve_0$ for some $\ve_0\in P_0$ and $\ve_i\in P_i$, where $i=1$~or~$2$.  We then have $x_i\ol\ve_i\ol\ve_i^* = x_0\ol\ve_0\ol\ve_i^*$.  Now, 
\[
\ol\ve_i\ol\ve_i^* = \ga \text{ or } \pi \text{ (for $i=1$ or $2$, respectively)}
\AND
\ol\ve_0\ol\ve_i^* = \ze \text{ or } \be,
\]
but we note that $\ol\ve_0\ol\ve_i^* = \be$ is only possible when $i=2$.  We see then that one of the following holds:
\[
x_1\ga = x_0\ze \COMMA 
x_2\pi = x_0\ze \OR 
x_2\pi = x_0\be.
\]
Keeping in mind $\ga=\ol\ga$ and $\pi=\ol\pi$, the first two options contradict \eqref{eq:Y1Y2ze}, so suppose instead that $x_2\pi = x_0\be$.  As noted above, this case arises when $i=2$ and $\ol\ve_0\ol\ve_2^* = \be$, and so our original assumption was that $x_2\ol\ve_2 = x_0\ol\ve_0$.  Putting this all together, we have
\[
x_2\ze = x_2\pi\ze = x_0\be\ze = 
x_0\ze = (x_0\ol\ve_0)\ol\ve_0^* = (x_2\ol\ve_2)\ol\ve_0^* = x_2(\ol\ve_0\ol\ve_2^*)^* = x_2\be^* = x_2\al,
\]
which contradicts \eqref{eq:x0x1x2}.  This completes the proof of \eqref{eq:Y1Y2Y0}.

Combining \eqref{eq:Y1Y2} and \eqref{eq:Y1Y2Y0}, we see that $Y_0$, $Y_1$ and $Y_2$ are pairwise disjoint.  Given \eqref{eq:YiPi}, it follows that the subset $Y_0\cup Y_1\cup Y_2$ of $X$ has size $|P_0|+|P_1|+|P_2| = |Q|$.  Thus, it remains only to show that
\begin{equation}\label{eq:XY0Y1Y2}
X \sm (Y_0\cup Y_1\cup Y_2) \not= \es.
\end{equation}
This is certainly true if either of $x_1\ze$ or $x_2\ze$ does not belong to $Y_0\cup Y_1\cup Y_2$.  Given \eqref{eq:Y1Y2ze}, the only alternative is that $x_1\ze$ and $x_2\ze$ both belong to $Y_0$, so we now assume that this is the case.  Thus, for $i=1,2$ we have $x_i\ze = x_0\ol\ve_i$ for some $\ve_i\in P_0$, and then $x_i\ze = x_i\ze\ze = x_0\ol\ve_i\ze = x_0\ze$, so in fact
\begin{equation}\label{eq:xize}
x_0\ze = x_1\ze = x_2\ze.
\end{equation}
We can complete the proof of \eqref{eq:XY0Y1Y2}, and hence of the proposition, by showing that
\[
x_2\al \not\in Y_0\cup Y_1\cup Y_2.
\]
To do so, suppose to the contrary that $x_2\al\in Y_i$ for some $i=0,1,2$, so that $x_2\al = x_i\ol\ve$ for some $\ve\in P_i$.  Keeping \eqref{eq:xize} in mind, we then have
\[
x_2\al = x_2\al\ze = x_i\ol\ve\ze = x_i\ze = x_2\ze,
\]
contradicting \eqref{eq:x0x1x2}.  As noted above, this completes the proof.
\epf

\section{The Temperley--Lieb monoid}\label{sect:TL}

In this section we consider the Temperley--Lieb monoid $\TL_n$, our main result being the following:

\begin{thm}\label{thm:TL}
For $n\geq3$ we have
\[
\deg'(\TL_n) = |Q| \AND \deg(\TL_n) = \degrc(\TL_n) = 1 + |Q|, 
\]
where
\[
Q = \begin{cases}
P_0\cup P_2\cup P_4 &\text{if $n$ is even}\\
P_1\cup P_3 &\text{if $n$ is odd.}
\end{cases}
\]
A formula for $|Q|$ can be found in Proposition \ref{prop:QTLn}.
\end{thm}

Here as usual $P_r = \set{\ve\in P(\TL_n)}{\rank(\ve)=r}$ for $0\leq r\leq n$, but we note that this set is non-empty precisely when $r\equiv n \pmod 2$.
The even case of Theorem \ref{thm:TL} in fact follows from Theorem~\ref{thm:M}, and the previously-mentioned isomorphism $\PP_n \cong \TL_{2n}$, which maps ${P_r(\PP_n) \to P_{2r}(\TL_{2n})}$; see Figure~\ref{fig:PPTL}.

We are therefore left to deal with the odd case of Theorem \ref{thm:TL}.  For this it is again convenient to work with an isomorphic copy of $\TL_{2n-1}$.  Specifically, it was explained in \cite[Section 1]{HR2005} that~$\TL_{2n-1}$ is isomorphic to the monoid
\[
M = \set{\al\in\PP_n}{1 \text{ and } 1' \text{ belong to the same block of }\al}.
\]
Indeed, the isomorphism $\PP_n\to\TL_{2n}$ maps an element $\al\in M$ to a Temperley--Lieb diagram $\widetilde\al\in\TL_{2n}$ with transversal $\{1,1'\}$, the set of which is clearly isomorphic to $\TL_{2n-1}$; again see Figure~\ref{fig:PPTL}.

We fix the above monoid $M(\cong\TL_{2n-1})$ for the remainder of the section.  We also write
\[
P = P(M) \ANd Q = P_1\cup P_2, \WHERe P_r = \set{\ve\in P}{\rank(\ve)=r} \qquad\text{for $1\leq r\leq n$.}
\]
(Because of its defining property, $M$ contains no partitions of rank $0$, but it contains partitions of every rank in $\{1,\ldots,n\}$, regardless of parity.)  We therefore wish to prove the following:

\begin{thm}\label{thm:TL2n-1}
For $n\geq2$ we have
\[
\deg'(M) = |Q| \AND \deg(M) = \degrc(M) = 1+|Q|, \WHERE Q = P_1\cup P_2.
\]
\end{thm}

Our strategy for proving this is similar to that for Theorem \ref{thm:M}.  We will therefore be somewhat briefer.

As before, we need to understand the minimal congruences of $M$.  Since $M\cong\TL_{2n-1}$, it follows from \cite[Theorem 9.1]{EMRT2018} that there are two of these, which we again denote by
\begin{align*}
\lam &= \De_M \cup \set{(\al,\be)\in M\times M}{\rank(\al)=\rank(\be)=1,\ \al\mr\L\be},\\
 \rho &= \De_M \cup \set{(\al,\be)\in M\times M}{\rank(\al)=\rank(\be)=1,\ \al\mr\R\be}.
\end{align*}
We also have $\lam = (\ze,\al)^\sharp$ and $\rho = (\ze,\be)^\sharp$, where 
\begin{equation}\label{eq:zab}
\ze = \custpartn{1,2,3,6}{1,2,3,6}{\stline11\udotted36 \ddotted36} \COMMA
\al = \custpartn{1,2,3,6}{1,2,3,6}{\stline11\stline21\uarcx12{.3}\udotted36 \ddotted36} \AND
\be = \custpartn{1,2,3,6}{1,2,3,6}{\stline11\stline12\darcx12{.3}\udotted36 \ddotted36}  .  
\end{equation}
We still have the pre-order $\pre$ on $M$, as in \eqref{eq:pre}, and its associated equivalence ${\approx}={\pre}\cap{\succeq}$, leading to the representation $\mu:Q^-\times M\to Q^-$, as in \eqref{eq:muQM}.  In the proofs that follow we make use of the fact that the map $\ve\mt\ol\ve$ from \eqref{eq:olve} still maps $P$ into $M$.  

\begin{thm}[cf.~Theorem \ref{thm:P}]\label{thm:M1}
For $n\geq2$, the action $\mu:Q^-\times M\to Q^-$ is faithful and monogenic, and consequently
\[
\deg(M) \leq \degrc(M) \leq 1+|Q|.
\]
\end{thm}

\pf
Taking $\pi = \custpartn{1,2,3,6}{1,2,3,6}{\stline11\stline22\udotted36 \ddotted36} \in Q$, faithfulness of $\mu$ follows from the fact that
\[
\mu(\pi,\ze) = \ze \COMMA \mu(\pi,\al) = - \AND \mu(\pi,\be) =  \custpartn{1,2,3,6}{1,2,3,6}{\stline11\stline22\uarc12\darc12\udotted36 \ddotted36}
\]
are distinct.  Monogenicity follows from the fact that $\ve = \mu(\pi,\ol\ve)$ for all $\ve\in Q$.
\epf

\begin{lemma}[cf.~Lemma \ref{lem:Msi}]\label{lem:MMsi}
Let $\si$ be a right congruence of $M$.
\ben
\item \label{MMsi1} If $\si$ separates $\{\ze,\al\}$, then it separates $\ol P_2$.
\item \label{MMsi2} If $\si$ separates $\{\ze,\be\}$, then it separates $\ol P_1$.
\een
\end{lemma}

\pf
\firstpfitem{\ref{MMsi1}}  Suppose $(\ol\ve_1,\ol\ve_2)\in\si$ for distinct $\ve_1,\ve_2\in P_2$, and write $\ol\ve_1$ and $\ol\ve_2$ as in \eqref{eq:ve12a}, noting that $1\in A<B$ and $1\in D<E$.  We must show that $(\ze,\al)\in\si$.

\pfcase1  Suppose first that $A\not=D$.  Without loss of generality, fix some $a\in A\sm D$, and let $b\in B$ be arbitrary, noting that $1<a<b$.  It follows that $\th = \begin{partn}{2} 1&a,b \\ \hhline{~|-} 1&\ \end{partn} \in M$.  We then have $(\al,\ze) = (\ol\ve_1\th,\ol\ve_2\th) \in \si$.  (Note that by the form of~$\ol\ve_2$ and~$\th$, the product $\ol\ve_2\th$ could only be equal to $\al$ or $\ze$.  It could only be equal to $\al$ if the edge~$\{a,b\}$ of $\th$ connected $D$ and $E$.  Since $a\not\in D$, this could only be the case if $a\in E$ and $b\in D$, but this is impossible since $a<b$ and $D<E$.)

\pfcase2  Next suppose $A=D$ but $B\not=E$, and without loss of generality, fix $b\in B\sm E$.  Then with $\th = \binom{1,b}{1}\in M$, we have $(\al,\ze) = (\ol\ve_1\th,\ol\ve_2\th) \in \si$.

\pfcase3  Finally, suppose $A=D$ and $B=E$, and without loss of generality fix some ${(x,y) \in \coker(\ve_1)\sm\coker(\ve_2)}$.  We assume that $x<y$, and we also let $b\in B$ be arbitrary.  By planarity of $\ve_1$, we either have $1<x<y<b$ or $1<b<x<y$.  We then define $\th = \begin{partn}{2} 1,x&b,y \\ \hhline{~|-} 1&\ \end{partn}$ or $\begin{partn}{2} 1,y&b,x \\ \hhline{~|-} 1&\ \end{partn}$, respectively, and we have $\th\in M$, and $(\al,\ze) = (\ol\ve_1\th,\ol\ve_2\th) \in \si$.

\pfitem{\ref{MMsi2}}  Suppose $(\ol\ve_1,\ol\ve_2)\in\si$ for distinct $\ve_1,\ve_2\in P_1$, and write $\ol\ve_1$ and $\ol\ve_2$ as in \eqref{eq:ve12b}.

\pfcase1  Suppose first that $A\not=C$, and without loss of generality fix some $a\in A\sm C$, noting that $a\not=1$.  Then with $\th = \begin{partn}{2} 1 & a \\ \hhline{~|~} 1&2 \end{partn}\in M$ we have $(\be,\ze) = (\ol\ve_1\th,\ol\ve_2\th) \in \si$.

\pfcase2  Now suppose $A=C$, and without loss of generality fix some $(x,y) \in \coker(\ve_1)\sm\coker(\ve_2)$ with $x<y$.  Then with $\th = \begin{partn}{2} 1,x & y \\ \hhline{~|~} 1&2 \end{partn}\in M$ we have $(\be,\ze) = (\ol\ve_1\th,\ol\ve_2\th) \in \si$.
\epf

The next result completes the proof of Theorem \ref{thm:TL2n-1}, and hence of Theorem \ref{thm:TL}.

\begin{prop}[cf.~Proposition \ref{prop:P}]
If $n\geq2$, then $\deg(M)\geq1+|Q|$.  
\end{prop}

\pf
Let $\mu:X\times M\to X$ be a faithful action, denoted $\mu(x,\de)=x\de$, as in the proof of Proposition~\ref{prop:P}.  We must show that $|X|\geq1+|Q|$.  Since $\mu$ is faithful, we can fix elements $x_1,x_2\in X$ such that 
\begin{equation}\label{eq:x1x2}
x_1\ze \not= x_1\be \AND x_2\ze \not= x_2\al.
\end{equation}
We define the sets $Y_i = x_i\ol P_i = \set{x_i\ol\ve}{\ve\in P_i}$ for $i=1,2$, noting that $|Y_i|=|P_i|$ by Lemmas~\ref{lem:act3} and~\ref{lem:MMsi}.  

We first claim that
\begin{equation}\label{eq:Y2}
Y_2\cap\set{x\ze}{x\in X} = \es.
\end{equation}
Indeed, suppose to the contrary that $x_2\ol\ve = x\ze$ for some $\ve\in P_2$ and $x\in X$.  Again writing ${\pi = \custpartn{1,2,3,6}{1,2,3,6}{\stline11\stline22\udotted36 \ddotted36}}$, we have $x_2\pi = x_2\ol\ve\ol\ve^* = x\ze\ol\ve^* = x\ze$.  It follows that
\[
x_2\ze = x_2\pi\ze = x\ze\ze = x\ze = x\ze\al = x_2\pi\al = x_2\al,
\]
contradicting \eqref{eq:x1x2}.

Next we claim that
\begin{equation}\label{eq:Y1Y2b}
Y_1\cap Y_2 = \es.
\end{equation}
Indeed, suppose to the contrary that $x_1\ol\ve_1 = x_2\ol\ve_2$ for some $\ve_1\in P_1$ and $\ve_2\in P_2$.  Then
\[
x_2\pi = x_2\ol\ve_2\ol\ve_2^* = x_1\ol\ve_1\ol\ve_2^*.
\]
Noting that $\ol\ve_1\ol\ve_2^* = \ze$ or $\be$, it follows that $x_2\pi = x_1\ze$ or $x_2\pi = x_1\be$.  The first contradicts \eqref{eq:Y2}, so suppose instead that $x_2\pi = x_1\be$.  This occurs when $\ol\ve_1\ol\ve_2^* = \be$, from which it follows that $\ol\ve_2\ol\ve_1^* = \be^* = \al$.  We then calculate
\[
x_1\ze = x_1\ol\ve_1\ol\ve_1^* = x_2\ol\ve_2\ol\ve_1^* = x_2\al \AND x_1\ze = x_1\be\ze = x_2\pi\ze = x_2\ze,
\]
so that $x_2\al = x_2\ze$, contradicting \eqref{eq:x1x2}.

Now that we have proved \eqref{eq:Y1Y2b}, it follows that $|Y_1\cup Y_2| = |Y_1|+|Y_2| = |P_1|+|P_2| = |Q|$, so it remains to show that
\begin{equation}\label{eq:XY1Y2}
X\sm(Y_1\cup Y_2) \not= \es.
\end{equation}
This is certainly true if $x_2\ze\not\in Y_1\cup Y_2$, so suppose this is not the case.  It follows from \eqref{eq:Y2} that $x_2\ze\in Y_1$, so we have $x_2\ze = x_1\ol\ve$ for some $\ve\in P_1$.  We then have $x_2\ze = x_2\ze\ze = x_1\ol\ve\ze = x_1\ze$.  We will complete the proof of \eqref{eq:XY1Y2}, and hence of the proposition, by showing that 
\[
x_2\al\not\in Y_1\cup Y_2.
\]
To do so, suppose to the contrary that $x_2\al = x_i\ol\ve_i$ for some $\ve_i\in P_i$, where $i=1$ or $2$.  Then
\[
x_2\al = x_2\al\ze = x_i\ol\ve_i\ze = x_i\ze = x_2\ze,
\]
using $x_1\ze=x_2\ze$ in the last step.  This again contradicts \eqref{eq:x1x2}.
\epf

\section{The Brauer monoid}\label{sect:B}

We now come to our last diagram monoid, the Brauer monoid $\B_n$.  Our main result here is the following, stated in terms of the projection sets $P_r = \set{\ve\in P(\B_n)}{\rank(\ve)=r}$:

\begin{thm}\label{thm:B}
For $n\geq3$ we have
\[
\deg'(\B_n) = \begin{cases}
|P_1| + 3|P_3| &\text{if $n$ is odd}\\
|P_0| + 2|P_2| + 3|P_4| &\text{if $n$ is even,}
\end{cases}
\AND \deg(\B_n) = 1 + \deg'(\B_n).
\]
If $n$ is odd, then $\deg(\B_n) = \degrc(\B_n)$.
\end{thm}

It was shown in \cite[Theorem 8.4]{EG2017} that $|P_r| = \tbinom nr (n-r-1)!!$.  Here as usual for a positive odd integer $m$ we define the double factorial ${m!! = m(m-2)(m-4)\cdots1}$, and by convention $(-1)!!=1$.  One can then easily check that Theorem \ref{thm:B} leads to the explicit formulae
\begin{equation}\label{eq:degBn}
\deg'(\B_n) = \begin{cases}
\frac{n+1}2\cdot n!! &\text{if $n\geq3$ is odd}\\
\frac{(n+4)(n+2)}8\cdot (n-1)!! &\text{if $n\geq4$ is even.}
\end{cases}
\end{equation}

\begin{rem}\label{rem:B}
As we will see, the methods we use here are necessarily rather different from those of previous sections.  For one thing, as indicated by the theorem itself, we generally have the strict inequality $\deg(\B_n) < \degrc(\B_n)$ for even $n\geq4$.
For example, it follows from the theorem that $\deg(\B_4) = 19$, and from GAP computations~\cite{GAP,Semigroups} that $\degrc(\B_4) = 22$.    
Another dissimilarity is that $\deg'(\B_n)$ is not simply a sum of $|P_r|$ parameters, but rather a (non-trivial) linear combination.  The reason for this will become clear as we progress.
\end{rem}

Our strategy for proving Theorem \ref{thm:B} is as follows.  After gathering the required preliminaries in Section \ref{subsect:Bprelim}, we build a right congruence $\si$ in Section \ref{subsect:Bsi}, using the construction from Section~\ref{subsect:si}.  We then use this to show that the claimed value for $\deg(\B_n)$ is an upper bound in Sections~\ref{subsect:Bup_odd} and~\ref{subsect:Bup_even} for odd and even $n$, respectively; see Theorems \ref{thm:Bnodd} and \ref{thm:Bneven}.  The even case is more involved, and the faithful action of the desired degree is constructed as a push-out of two smaller actions (only one of which involves $\si$).  Finally, we show that the claimed value of $\deg(\B_n)$ is a lower bound in Sections \ref{subsect:Blow_odd} and \ref{subsect:Blow_even}, again split by parity; see Propositions \ref{prop:Blowodd} and \ref{prop:Bloweven}.  As usual, the minimum-degree faithful actions we construct have global fixed points, so it follows that~${\deg'(\B_n) = \deg(\B_n) - 1}$.

\subsection{Preliminaries}\label{subsect:Bprelim}

In what follows, we will need the characterisation of Green's relations on $\B_n$ from \cite[Proposition~3.1]{DE2018}.  For $\al,\be\in\B_n$ we have
\begin{align*}
\al\leqL \be &\iff \coker(\al) \supseteq\coker(\be) ,&  \al\mr\L \be&\iff \coker(\al) =\coker(\be),\\
\al \leqR \be &\iff \ker(\al)\supseteq\ker(\be),& \al \mr\R \be&\iff \ker(\al)=\ker(\be),\\
\al \leqJ \be &\iff \rank(\al)\leq\rank(\be),& \al \mr\J \be&\iff \rank(\al)=\rank(\be),
\end{align*}
and of course $\D=\J$ as $\B_n$ is finite.  Keeping in mind that the rank of a Brauer partition has the same parity as $n$, it follows that the $\D$-classes and non-empty ideals of $\B_n$ are the sets
\[
D_r = \set{\al\in\B_n}{\rank(\al)=r} \AND I_r = \set{\al\in\B_n}{\rank(\al)\leq r} 
\]
for $0\leq r\leq n$ with $r\equiv n$ (mod $2$).  Note that $I_r = D_r\cup D_{r-2}\cup D_{r-4}\cup\cdots$.

\begin{rem}\label{rem:preBn}
It is also worth noting that $\leqR$ is precisely the pre-order $\pre$ used in Section \ref{sect:P}; see~\eqref{eq:pre}.  As in Remark~\ref{rem:MPTP}, it follows that the resulting action of $\B_n$ on $P^-=P\cup\{-\}$ is not faithful, hence the need for the different approach we take here.
\end{rem}

Throughout this section we write 
\[
P = P(\B_n) = \set{\ve\in\B_n}{\ve^2=\ve=\ve^*} \AND P_r = P_r(\B_n) = P\cap D_r \qquad\text{for $0\leq r\leq n$,}
\]
noting that $P_r$ is non-empty precisely when $r\equiv n$ (mod $2$).  We also write
\[
p_r = |P_r| = \tbinom nr (n-r-1)!!.
\]

As before, we will need to understand the minimal congruences of $\B_n$, as described in \cite[Section~8]{EMRT2018}.  For odd~$n$, there are two minimal congruences, which we will denote by
\[
\lam = \De_{\B_n} \cup \L\restr_{D_1} \AND \rho = \De_{\B_n} \cup \R\restr_{D_1},
\]
and we have $\lam = (\ze,\al)^\sharp$ and $\rho = (\ze,\be)^\sharp$, where
\begin{equation}\label{eq:Bzab}
\ze = \custpartn{1,2,3,4,5,8,9}{1,2,3,4,5,8,9}{\uarc23\uarc45\uarc89\darc23\darc45\darc89\udotted58\ddotted58\stline11} \COMMA
\al = \custpartn{1,2,3,4,5,8,9}{1,2,3,4,5,8,9}{\uarc12\uarc45\uarc89\darc23\darc45\darc89\udotted58\ddotted58\stline31} \AND
\be = \custpartn{1,2,3,4,5,8,9}{1,2,3,4,5,8,9}{\uarc23\uarc45\uarc89\darc12\darc45\darc89\udotted58\ddotted58\stline13}.
\end{equation}
For even $n$ there are three minimal congruences, which we will denote by
\[
\lam = \De_{\B_n} \cup \L\restr_{D_0} \COMMA \rho = \De_{\B_n} \cup \R\restr_{D_0} \AND \eta = \De_{\B_n} \cup \H\restr_{D_2},
\]
and we have $\lam = (\ze,\al)^\sharp$, $\rho = (\ze,\be)^\sharp$ and $\eta = (\ga,\de)^\sharp$, where
\begin{align}
\nonumber \ze &= \custpartn{1,2,3,4,5,6,9,10}{1,2,3,4,5,6,9,10}{\uarc12\uarc34\uarc56\uarc9{10}\darc12\darc34\darc56\darc9{10}\udotted69\ddotted69}, 
& \al &= \custpartn{1,2,3,4,5,6,9,10}{1,2,3,4,5,6,9,10}{\uarcx14{.6}\uarcx23{.3}\uarc56\uarc9{10}\darc12\darc34\darc56\darc9{10}\udotted69\ddotted69}, 
& \ga &= \custpartn{1,2,3,4,5,6,9,10}{1,2,3,4,5,6,9,10}{\uarc34\uarc56\uarc9{10}\darc34\darc56\darc9{10}\stline11\stline22\udotted69\ddotted69}, \\[5mm]
\label{eq:Bzabcd} && \be &= \custpartn{1,2,3,4,5,6,9,10}{1,2,3,4,5,6,9,10}{\darcx14{.6}\darcx23{.3}\darc56\darc9{10}\uarc12\uarc34\uarc56\uarc9{10}\udotted69\ddotted69}, 
& \de &= \custpartn{1,2,3,4,5,6,9,10}{1,2,3,4,5,6,9,10}{\uarc34\uarc56\uarc9{10}\darc34\darc56\darc9{10}\stline12\stline21\udotted69\ddotted69}.
\end{align}

\subsection{A right congruence}\label{subsect:Bsi}

We will soon have to split our consideration into cases according to the parity of $n$.  However, a certain right congruence $\si$ that is parity-independent will play an important role, and we introduce this here.  This will have the form $\si = \RR_I\vee\L^U = \nab_I\cup\L^U\restr_K$, as in Section \ref{subsect:si}; see~\eqref{eq:si}.

Throughout this section we assume that $n\geq3$, and we fix $k = \lfloor n/2\rfloor$, so that~${n=2k}$~or~${2k+1}$.  We also fix the partition
\[
\ze = \custpartn{1,2,5,6}{1,2,5,6}{\uarc12\uarc56\darc12\darc56\udotted25\ddotted25}
\OR
\ze = \custpartn{1,2,3,6,7}{1,2,3,6,7}{\uarc23\uarc67\darc23\darc67\udotted36\ddotted36\stline11}.
\]
We write $\ka = \ker(\ze) = \coker(\ze)$, and we denote the non-trivial $\ka$-classes by
\[
Z_1<\cdots<Z_k, \qquad\text{so}\qquad Z_i = \begin{cases}
\{2i,2i+1\} &\text{if $n$ is odd}\\
\{2i-1,2i\} &\text{if $n$ is even.}
\end{cases}
\]
When $n$ is even, these are \emph{all} the $\ka$-classes, but when $n$ is odd, the last $\ka$-class is $\{1\}$.
Define the $\R$-class $R = R_\ze$, and let $T$, $K$ and $I = \B_n\sm K$ be as in \eqref{eq:TKI}.  We then have
\[
R = \set{\al\in\B_n}{\ker(\al) = \ka} \COMMa K = \set{\al\in\B_n}{\ker(\al) \sub \ka} \ANd I = \set{\al\in\B_n}{\ker(\al) \not\sub \ka}.
\]
Note that the condition $\ker(\al)\sub\ka$ says that the non-trivial $\ker(\al)$-classes are among the $Z_i$.
Since~$\B_n$ is stable (as it is finite), and since $D_\ze$ (being $D_0$ or $D_1$) is $\H$-trivial, it follows from Lemma~\ref{lem:TK}\ref{TK3} that
\[
T = \set{\al\in\B_n}{\al\ze = \ze}.
\]
This set is not easy to describe concisely, but we will instead work with its sub(semi)group
\[
U = T\cap\S_n = \set{\al\in\S_n}{\al\ze = \ze} = \set{\al\in\S_n}{(\forall i\in\bk)(\exists j\in\bk) \ Z_i\al = Z_j}.
\]
This subgroup of $\S_n$ is the stabiliser of the partition $\bn/\ka$, and is (isomorphic to) the wreath product $C_2\wr\S_k$.  In later arguments (in which we will need to vary the underlying parameter $n$), we will also denote this group by $U = \G_n$, and we note that $|\G_n| = 2^kk!$.  

In what follows, it will be convenient in some circumstances to omit some upper and/or lower non-transversals when using the tabular notation for Brauer partitions, although we always list \emph{all} transversals.  Specifically, if we write $\al = \begin{partn}{6} a_1&\cdots&a_r&A_1&\cdots&A_s\\ \hhline{~|~|~|-|-|-} b_1&\cdots&b_r&B_1&\cdots&B_t\end{partn}$, and if the unlisted upper vertices are $c_1<\cdots<c_u$ (there must be an even number of these), then we assume the remaining upper blocks of $\al$ are $\{c_1,c_2\},\{c_3,c_4\},\ldots,\{c_{u-1},c_u\}$.  A similar statement holds for unlisted lower vertices/blocks.
Thus, for example, the elements $\al$ and $\be$ in \eqref{eq:Bzab} can be denoted as~$\al = \binom31$ and~$\be=\binom13$, while~$\al$ in~\eqref{eq:Bzabcd} can be written as either $\al = \begin{partn}{1} 1,4 \\ \hhline{-} \ \end{partn}$ or $\begin{partn}{1} 2,3 \\ \hhline{-} \ \end{partn}$.

The coming arguments will also use the mapping
\begin{equation}\label{eq:olal}
\B_n\to K: \al = \begin{partn}{6} a_1&\cdots&a_r&A_1&\cdots&A_s\\ \hhline{~|~|~|-|-|-} b_1&\cdots&b_r&B_1&\cdots&B_s\end{partn} 
\mt
\ol\al = \begin{partn}{6} 1&\cdots&r&\multicolumn{3}{c}{}\\ \hhline{~|~|~|-|-|-} b_1&\cdots&b_r&B_1&\cdots&B_s\end{partn}.
\end{equation}
We assume that \emph{all} the blocks of $\al$ have been listed here, and as usual we have $a_1<\cdots<a_r$.  The upper non-transversals of $\ol\al$ have been omitted, however, as per the above convention, and these are in fact $Z_{k-s+1},\ldots,Z_k$.  
Note also that $\ol\al = \th\al$, where $\th\in\S_n$ is any permutation for which:
\begin{equation}\label{eq:th}
\text{$i\th = a_i$ for each $i\in\br$,} \AND \text{$Z_{k-s+i}\th = A_i$ for each $i\in\bs$.}
\end{equation}

\begin{lemma}\label{lem:olal}
For any $\al\in K$ we have $\al \mr\L^U \ol\al$.
\end{lemma}

\pf
Let $\al\in K$.  Since $\ker(\al)\sub\ka$, we can write $\al = \begin{partn}{6} a_1&\cdots&a_r&Z_{i_1}&\cdots&Z_{i_s}\\ \hhline{~|~|~|-|-|-} b_1&\cdots&b_r&B_1&\cdots&B_s\end{partn}$.  Also writing ${\bk\sm\{i_1,\ldots,i_s\}} = \{j_1<\cdots<j_t\}$, we note that:
\bit
\item $Z_{j_1} = \{a_1,a_2\}$, $Z_{j_2} = \{a_3,a_4\}$, and so on, when $n$ is even, or
\item $Z_{j_1} = \{a_2,a_3\}$, $Z_{j_2} = \{a_4,a_5\}$, and so on, when $n$ is odd, in which case we also have $a_1 = 1$.
\eit
It follows that any permutation $\th$ as in \eqref{eq:th} belongs to~$U$, and $\ol\al=\th\al$ gives the claim.
\epf

We now let
\begin{equation}\label{eq:siBn}
\si = \RR_I\vee\L^U = \nab_I \cup \L^U\restr_K  = \bigset{(\al,\be)}{\al,\be\in I \text{ or } [\al,\be\in K \text{ and } U\al = U\be]\ }
\end{equation}
be the right congruence of $\B_n$ from \eqref{eq:si}.  Of course $I$ is a $\si$-class.  By definition, and since ${\L^U\sub\L(=\L^{\B_n})}$, every other $\si$-class is contained in $L\cap K$ for some $\L$-class $L$ of $\B_n$.  It follows from Lemma~\ref{lem:LIK} that every such $\si$-class is also an $\L^U$-class.  The next result tells us how many such classes are contained in $L\cap K$.  For the proof, we define a map
\[
\phi:\B_n\to \S_n\cup\S_{n-2}\cup\S_{n-4}\cup\cdots,
\]
as follows.  Let $\al\in\B_n$, and write $\dom(\al)=\{a_1<\cdots<a_r\}$ and $\codom(\al)=\{b_1<\cdots<b_r\}$.  Then we define $\phi(\al)$ to be the permutation $f\in\S_r$ for which $\al$ contains the transversal $\{a_i,b_{if}'\}$ for each $i$.

\begin{lemma}\label{lem:Gr}
Let $L$ be an $\L$-class of $\B_n$, and let the common rank of its elements be $r$.  Then
\[
|(L\cap K)/\si| = |(L\cap K)/\L^U| = r!/|\G_r|.
\]
Consequently,
\[
|(D_r\cap K)/\si| = |(D_r\cap K)/\L^U| = p_rr!/|\G_r|.
\]
\end{lemma}

\pf
Since $D_r$ contains $p_r$ $\L$-classes, it suffices to prove the first claim.  We do this by showing that 
\[
\al\mr\si\be \IFF \phi(\al)\phi(\be)^{-1} \in \G_r \qquad\text{for all $\al,\be\in L\cap K$.}
\]
So let $\al,\be\in L\cap K$, write $f=\phi(\al)$ and $g=\phi(\be)$, let ${\codom(\al)=\codom(\be) = \{b_1<\cdots<b_r\}}$, and let (all) the non-trivial $\coker(\al)=\coker(\be)$-classes be $B_1,\ldots,B_s$.  We then have
\begin{equation}\label{eq:olalolbe}
\ol\al = \begin{partn}{6} 1&\cdots&r&\multicolumn{3}{c}{}\\ \hhline{~|~|~|-|-|-} b_{1\!f}&\cdots&b_{r\!f}&B_1&\cdots&B_s\end{partn}
\AND
\ol\be = \begin{partn}{6} 1&\cdots&r&\multicolumn{3}{c}{}\\ \hhline{~|~|~|-|-|-} b_{1g}&\cdots&b_{rg}&B_1&\cdots&B_s\end{partn}.
\end{equation}
Since $\al,\be\in K$, it follows from Lemma \ref{lem:olal} (and the definition of $\si$) that
\[
\al\mr\si\be \Iff \al\mr\L^U\be \Iff \ol\al\mr\L^U\ol\be \Iff \ol\al \in U\ol\be.
\]
So it remains to show that
\begin{equation}\label{eq:olaltholbe}
\ol\al \in U\ol\be \IFF f g^{-1}\in\G_r.
\end{equation}
For the forward implication, suppose $\ol\al = \th\ol\be$ for some $\th\in U$.  Examining \eqref{eq:olalolbe}, we see that~$\th$ fixes $\br=\{1,\ldots,r\}$ set-wise, and so $\vt=\th\restr_\br\in\S_r$.  Moreover, it follows from $\th\in U = \G_n$ that~$\vt\in\G_r$.  Now, for every $i\in\br$, the product $\th\ol\be$ contains the transversal $\{i,b_{i\vt g}'\}$.  Keeping in mind that $\th\ol\be = \ol\al$, and comparing with~\eqref{eq:olalolbe}, it follows that $f = \vt g$, i.e.~$f g^{-1} = \vt\in \G_r$.

Conversely, if $f g^{-1}\in\G_r$, then $\ol\al = \th\ol\be$, where $\th\in U$ is such that:
\[
i\th = if g^{-1}\text{ for $1\leq i\leq r$} \AND j\th = j \text{ for $r+1\leq j\leq n$.}
\]
This completes the proof of \eqref{eq:olaltholbe}, and hence of the lemma.
\epf

\begin{rem}
Using $|\G_r| = 2^ss!$, where $s=\lfloor r/2\rfloor$, one can show that $r!/|\G_r| = r!!$ or $(r-1)!!$ for $r$ odd or even, respectively.
\end{rem}

\subsection{Upper bound -- odd case}\label{subsect:Bup_odd}

We are now ready to construct faithful actions of the degree specified in Theorem \ref{thm:B}.  Here we consider the case that $n\geq3$ is odd, and the even case will be treated in Section \ref{subsect:Bup_even}.  We continue to fix the partition $\ze= \custpartn{0,1,2,5,6}{0,1,2,5,6}{\uarc12\uarc56\darc12\darc56\udotted25\ddotted25\stline00}$, and the right congruence $\si$ from \eqref{eq:siBn}.  
Let
\[
\Om = \bigset{[\al]_\si}{\al\in I_3} , \WHERE I_3 = D_1\cup D_3.
\]
Since $I_3$ is an ideal of $\B_n$, the action of $\B_n$ on $\B_n/\si$ restricts to an action
\begin{equation}\label{eq:muBn}
\mu : \Om\times\B_n \to \Om \GIVENBY \mu([\al]_\si,\be) = [\al]_\si\be = [\al\be]_\si.
\end{equation}

\begin{thm}\label{thm:Bnodd}
For odd $n\geq3$, the action $\mu : \Om\times\B_n \to \Om$ is faithful and monogenic, and consequently
\[
\deg(\B_n) \leq \degrc(\B_n) \leq 1 + 3p_3 + p_1.
\]
\end{thm}

\pf
As explained in Section \ref{subsect:Bprelim}, the minimal congruences of $\B_n$ are $(\ze,\al)^\sharp$ and $(\ze,\be)^\sharp$, where~$\al$ and~$\be$ are as in \eqref{eq:Bzab}.  Thus, by Lemma \ref{lem:act1}, we can establish faithfulness by showing that $\mu$ separates $(\ze,\al)$ and $(\ze,\be)$.  For this we let ${\pi = \custpartn{1,2,3,4,5,8,9}{1,2,3,4,5,8,9}{\uarc45\uarc89\darc45\darc89\udotted58\ddotted58\stline11\stline22\stline33}}$, and we note that 
\[
[\pi]_\si\ze = [\ze]_\si \COMMA
[\pi]_\si\al = [\al]_\si = I \AND
[\pi]_\si\be = [\be]_\si .
\]
These are distinct because $\ze,\be\in K$, and $[\ze]_\si=U\ze=\{\ze\}$.

For monogenicity, we claim that $\Om$ is generated by $[\pi]_\si$.  For this, we have seen that $I = [\pi]_\si\al$.  Any other $\si$-class has the form $[\de]_\si$ for some $\de\in I_3\cap K$, and by Lemma \ref{lem:olal} we have
\[
[\de]_\si = [\ol\de]_\si = [\pi\ol\de]_\si = [\pi]_\si\ol\de.
\]

Finally, using Lemma \ref{lem:Gr}, we have
\[
\degrc(\B_n) \leq |\Om| = 1 + |(D_3\cap K)/\si| + |(D_1\cap K)/\si| = 1 + 3p_3 + p_1.  \qedhere
\]
\epf

\subsection{Upper bound -- even case}\label{subsect:Bup_even}

Now we consider the case that $n=2k\geq4$ is even.  We fix the partition $\ze= \custpartn{1,2,5,6}{1,2,5,6}{\uarc12\uarc56\darc12\darc56\udotted25\ddotted25}$, and the right congruence $\si$ from \eqref{eq:siBn}.  
Unlike the odd case, minimum-degree faithful actions of~$\B_n$ are not always monogenic; cf.~Remark \ref{rem:B}.  Consequently, we will define such an action by first defining two separate actions $\mu_1$ and $\mu_2$, and then taking an appropriate push-out $\mu_1\sqcup_\xi\mu_2$.

In all that follows, an important role will be played by the set
\[
J = I_2\cap K = I_2 \sm I = \set{\al\in\B_n}{\rank(\al)\leq2,\ \ker(\al)\sub\ka}.
\]

\begin{lemma}\label{lem:J}
$J$ is a right ideal of $\B_n$.
\end{lemma}

\pf
Let $\al\in J$ and $\be\in\B_n$; we must show that $\al\be\in J$, which amounts to showing that $\al\be\in K$, as we certainly have $\al\be\in I_2$.  Now $\al\be$ contains all the upper blocks of $\al$, and these are either
\[
Z_1,\ldots,Z_k \OR Z_1,\ldots,Z_{i-1},Z_{i+1},\ldots,Z_k \text{\ \ for some $i\in\bk$.}
\]
In the first case, the upper blocks of $\al\be$ are of course $Z_1,\ldots,Z_k$.  In the second case, the upper blocks of $\al\be$ are the same as for $\al$ if $\rank(\al\be)=2$, or else $Z_1,\ldots,Z_k$ if $\rank(\al\be)=0$.  In all cases, it follows that $\al\be\in K$.
\epf

Now, in addition to $\ze = \custpartn{1,2,3,4,7,8}{1,2,3,4,7,8}{\uarc12\uarc34\uarc78\darc12\darc34\darc78\udotted47\ddotted47}$, we also define $\ga = \custpartn{1,2,3,4,7,8}{1,2,3,4,7,8}{\uarc34\uarc78\darc34\darc78\stline11\stline22\udotted47\ddotted47}$.  It follows from the proof of Lemma \ref{lem:J} that 
\begin{equation}\label{eq:Om1}
\Om_1 = R_\ga\cup R_\ze = \set{\al\in\B_n}{\al\mr\R\ga \text{ or } \al\mr\R\ze}
\end{equation}
is a right ideal of $\B_n$.  Consequently, we have an action
\[
\mu_1:\Om_1\times\B_n\to\Om_1 \GIVENBY \mu_1(\al,\be) = \al\be.
\]

To define the second action, we need to work a little harder.  We first define a map
\[
\B_n \to I_0=D_0 : \al = \begin{partn}{6} a_1&\cdots&a_r&A_1&\cdots&A_s\\ \hhline{~|~|~|-|-|-} b_1&\cdots&b_r&B_1&\cdots&B_s\end{partn}  
\mt
\wh\al = \begin{partn}{6} a_1,a_2&\cdots&a_{r-1},a_r&A_1&\cdots&A_s\\ \hhline{-|-|-|-|-|-} b_1,b_2&\cdots&b_{r-1,}b_r&B_1&\cdots&B_s\end{partn} .
\]
In particular, note that $\wh\al=\al$ for $\al\in D_0$.  By \cite[Theorem 8.4]{EMRT2018}, the relation
\[
\chi = \De_{\B_n} \cup \set{(\al,\be)\in I_2\times I_2}{\wh\al = \wh\be}
\]
is a (two-sided) congruence on $\B_n$.  We now define the relation
\[
\tau = \nab_I \sqcup \chi\restr_J \sqcup \L^U\restr_{K\sm J},
\]
and we note that $K\sm J = K\cap(D_4\cup D_6\cup\cdots\cup D_n)$.

\begin{lemma}
$\tau$ is a right congruence on $\B_n$.
\end{lemma}

\pf
Since $\B_n = I\sqcup K$ and $K=J\sqcup(K\sm J)$, $\tau$ is the union of equivalences on disjoint sets, and is therefore an equivalence itself.  For right-compatibility, let $(\al,\be)\in\tau$ and $\th\in\B_n$; we must show that $(\al\th,\be\th)\in\tau$.  This is clear if $(\al,\be)\in\nab_I$.

Next suppose $(\al,\be)\in\chi\restr_J$.  Since $\chi$ is a congruence we have $(\al\th,\be\th)\in\chi$, and since $J$ is a right ideal we have $\al\th,\be\th\in J$.  Together these give $(\al\th,\be\th)\in\chi\restr_J\sub\tau$.

Finally suppose $(\al,\be)\in\L^U\restr_{K\sm J}$.  Since $\L^U\restr_{K\sm J}\sub\L^U\restr_K\sub\si$, it follows that $(\al\th,\be\th)$ belongs to $\si$, and hence to either $\nab_I$ or $\L^U\restr_K$.  In the first case we are done, so suppose $(\al\th,\be\th)\in\L^U\restr_K$.  Since $\L^U\sub\D$ we have $\rank(\al\th)=\rank(\be\th) = r$, say.  
\bit
\item If $r\geq4$ then $(\al\th,\be\th)\in\L^U\restr_{K\sm J} \sub \tau$.
\item If $r=0$ then $\al\th$ and $\be\th$ belong to $D_0\cap K = R$, and are $\L^U$- and hence $\L$-related.  Since~$D_0$ is $\H$-trivial, it follows that $\al\th = \be\th$ in this case.
\item Finally, if $r=2$, then from $\al\th,\be\th\in D_2\cap K(\sub J)$ and $(\al\th,\be\th)\in\L^U\sub\L$, we can write
\[
\al\th = \begin{partn}{5} 2i-1&2i&\multicolumn{3}{c}{}\\ \hhline{~|~|-|-|-} a&b&C_1&\cdots&C_{k-1}\end{partn}
\AND
\be\th = \begin{partn}{5} 2j-1&2j&\multicolumn{3}{c}{}\\ \hhline{~|~|-|-|-} a&b&C_1&\cdots&C_{k-1}\end{partn}
\qquad\text{for some $1\leq i,j\leq k$.}
\]
We then have $\widehat{\al\th} = \begin{partn}{4} \multicolumn{4}{c}{} \\ \hhline{-|-|-|-} a,b&C_1&\cdots&C_{k-1}\end{partn} = \widehat{\be\th}$, so that $(\al\th,\be\th)\in\chi\restr_J\sub\tau$.  \qedhere
\eit
\epf

Now let
\[
\Om_2 = \bigset{[\al]_\tau}{\al\in I_4}, \WHERE I_4 = D_0\cup D_2\cup D_4.
\]
Since $I_4$ is an ideal, the action of $\B_n$ on $\B_n/\tau$ restricts to an action
\[
\mu_2 : \Om_2\times\B_n \to \Om_2 \GIVENBY \mu_2([\al]_\tau,\be) = [\al\be]_\tau.
\]

It follows from the definition of $\tau$ that
\begin{equation}\label{eq:Om2}
\Om_2 = \{I\} \cup (D_4\cap K)/\L^U \cup \Om_3 , \WHERE \Om_3 = \bigset{[\al]_\tau}{\al\in J} = \bigset{[\al]_\tau}{\al\in R_\ze},
\end{equation}
and that each $\tau$-class from $\Om_3$ is in fact a $\chi$-class.  Since $\wh\al=\al$ for all $\al\in R_\ze$, the map
\[
\xi:R_\ze \to \Om_3: \al\mt[\al]_\tau
\]
is a bijection.  Since $R_\ze$ is a right ideal of $\B_n$, it follows that $R_\ze$ and $\Om_3$ are sub-acts of $\Om_1$ and~$\Om_2$, respectively, and it is easy to see that $\xi$ is in fact an isomorphism of these sub-acts.  As in \eqref{eq:pushout}, we can therefore form the push-out
\[
\mu = \mu_1\sqcup_\xi\mu_2 : \Om\times\B_n\to\Om, \WHERE \Om = \Om_1\sqcup_\xi\Om_2.
\]

\begin{thm}\label{thm:Bneven}
For even $n\geq4$, the action $\mu: \Om\times\B_n\to\Om$ is faithful, and consequently
\[
\deg(\B_n) \leq 1 + 3p_4 + 2p_2 + p_0 .
\]
\end{thm}

\pf
As explained in Section \ref{subsect:Bprelim}, the minimal congruences of $\B_n$ are $(\ze,\al)^\sharp$, $(\ze,\be)^\sharp$ and~$(\ga,\de)^\sharp$, where~these partitions are as in \eqref{eq:Bzabcd}.  Thus, by Lemma \ref{lem:act1}, we can establish faithfulness by showing that $\mu$ separates $(\ze,\al)$, $(\ze,\be)$ and $(\ga,\de)$.  Since $\mu = \mu_1\sqcup_\xi\mu_2$, we can do this by showing that $\mu_1$ separates $(\ze,\be)$ and $(\ga,\de)$, while $\mu_2$ separates $(\ze,\al)$.  For the former we have
\[
\mu_1(\ze,\ze) = \ze \COMMA
\mu_1(\ze,\be) = \be \COMMA
\mu_1(\ga,\ga) = \ga \AND
\mu_1(\ga,\de) = \de ,
\]
which are all distinct.
For the latter, and with $\pi = \custpartn{1,2,3,4,5,6,9,10}{1,2,3,4,5,6,9,10}{\uarc56\uarc9{10}\darc56\darc9{10}\stline11\stline22\stline33\stline44\udotted69\ddotted69}$, we have
\[
\mu_2([\pi]_\tau,\ze) = [\ze]_\tau \AND \mu_2([\pi]_\tau,\al) = [\al]_\tau = I.
\]

It now follows that
\[
\deg(\B_n) \leq |\Om| = |\Om_1\sqcup_\xi\Om_2| = |\Om_1| + |\Om_2| - |\Om_3|.
\]
Combining \eqref{eq:Om1} and~\eqref{eq:Om2} with Lemma \ref{lem:Gr} we have
\[
|\Om_1| = |R_\ga|+|R_\ze| = 2p_2 + p_0 \AND |\Om_2| = 1 + 3p_4 + |\Om_3|.
\]
It follows that $|\Om| = 1 + 3p_4 + 2p_2 + p_0$.  
\epf

\subsection{Lower bound -- odd case}\label{subsect:Blow_odd}

We now turn to the task of showing that the claimed value for $\deg(\B_n)$ in Theorem \ref{thm:B} is a lower bound, again treating the odd and even cases in separate sections.

Here we assume $n=2k+1\geq3$ is odd, and we fix $\ze,\al,\be\in\B_n$, as in \eqref{eq:Bzab}.  For a subset ${X\sub P = P(\B_n)}$ we write $\ol X = \set{\ol\ve}{\ve\in X}$, where here $\ol\ve$ is as in \eqref{eq:olal}.  We will also use cycle notation for permutations from $\S_n\sub\B_n$, so that for example $(1,2,3) = \custpartn{1,2,3,4,7}{1,2,3,4,7}{\stline12\stline23\stline31\stline44\stline77\udotted47\ddotted47}$.  We fix the cyclic subgroup
\[
\C_3 = \la (1,2,3)\ra = \{\id_n,(1,2,3),(1,3,2)\} \leq\S_n.
\]

\begin{lemma}\label{lem:Bsiodd}
Let $\si$ be a right congruence of $\B_n$, where $n\geq3$ is odd.
\ben
\item \label{Bsiodd1} If $\si$ separates $\{\ze,\al\}$, then it separates $\C_3\ol P_3 = \set{\xi\ol\ve}{\xi\in\C_3,\ \ve\in P_3}$.
\item \label{Bsiodd2} If $\si$ separates $\{\ze,\be\}$, then it separates $\ol P_1$.
\een
\end{lemma}

\pf
\firstpfitem{\ref{Bsiodd1}}  As usual, we prove the contrapositive.  So suppose $(\xi_1\ol\ve_1,\xi_2\ol\ve_2)\in\si$ for some ${\xi_1,\xi_2\in\C_3}$ and $\ve_1,\ve_2\in P_3$, with $\xi_1\ol\ve_1 \not= \xi_2\ol\ve_2$; we must show that $(\ze,\al)\in\si$.  Write
\[
\xi_1\ol\ve_1 = \begin{partn}{6} 1&2&3&\multicolumn{3}{c}{}\\ \hhline{~|~|~|-|-|-} a_1&a_2&a_3&A_1&\cdots&A_{k-1}\end{partn}
\AND
\xi_2\ol\ve_2 = \begin{partn}{6} 1&2&3&\multicolumn{3}{c}{}\\ \hhline{~|~|~|-|-|-} b_1&b_2&b_3&B_1&\cdots&B_{k-1}\end{partn},
\]
and let $A=\{a_1,a_2,a_3\}$ and $B=\{b_1,b_2,b_3\}$.  We now split into cases; in each we define an element $\th\in\B_n$ for which $\{\xi_1\ol\ve_1\th,\xi_2\ol\ve_2\th\} = \{\ze,\al\}$.  

\pfcase1  Suppose first that $A\not=B$.  If $a_1\not\in\{b_1,b_2\}$ or $a_3\not\in\{b_2,b_3\}$, then we take $\th = \begin{partn}{2} a_1&b_1,b_2 \\ \hhline{~|-} 1&\ \end{partn}$ or $\th = \begin{partn}{2} a_3&b_2,b_3 \\ \hhline{~|-} 1&\ \end{partn}$, respectively.  The cases in which $b_1\not\in\{a_1,a_2\}$ or $b_3\not\in\{a_2,a_3\}$ are analogous.  This leave us to consider the case in which
\[
a_1\in\{b_1,b_2\} \COMMA a_3\in\{b_2,b_3\} \COMMA b_1\in\{a_1,a_2\} \AND b_3\in\{a_2,a_3\}.
\]
Since $A\not=B$, this implies that $a_2\not\in B$ and $b_2\not\in A$, so in fact $a_1=b_1$ and $a_3=b_3$.  We can then assume without loss of generality that $B_1=\{a_2,x\}$ for some $x\in\bn$, and we take $\th = \begin{partn}{2} x&a_1,a_2 \\ \hhline{~|-} 1&\ \end{partn}$.

\pfcase2  Next suppose $a_i=b_i$ for $i=1,2,3$.  Since $\xi_1\ol\ve_1 \not= \xi_2\ol\ve_2$, we can assume without loss of generality that $A_1=\{u,v\}$, $B_1=\{u,x\}$ and $B_2=\{v,y\}$ for distinct $u,v,x,y\in\bn$.  We then take $\th = \begin{partn}{3} x&a_1,u&a_2,v \\ \hhline{~|-|-} 1&\multicolumn{2}{c}{} \end{partn}$.

\pfcase3  Up to symmetry, the final case to consider is where $b_1=a_2$, $b_2=a_3$ and $b_3=a_1$, and we then take $\th = \begin{partn}{2} a_1&a_2,a_3 \\ \hhline{~|-} 1&\ \end{partn}$.

\pfitem{\ref{Bsiodd2}}  Suppose $(\ol\ve_1,\ol\ve_2)\in\si$ for distinct $\ve_1,\ve_2\in P_1$.  Since then $\coker(\ve_1)\not=\coker(\ve_2)$, there exist distinct $u,v,w\in\bn$ such that $(u,v)\in\coker(\ve_1)$ and $(v,w)\in\coker(\ve_2)$.  Then with $\th = \begin{partn}{3} u&v&w \\ \hhline{~|~|~} 1&2&3\end{partn}$ we have $(\be,\ze) = (\ol\ve_1\th,\ol\ve_2\th) \in \si$.
\epf

We can now complete the proof of Theorem \ref{thm:B} in the odd case:

\begin{prop}\label{prop:Blowodd}
If $n\geq3$ is odd, then $\deg(\B_n) \geq 3p_3+p_1+1$.
\end{prop}

\pf
Let $\mu:X\times \B_n\to X$ be a faithful action, denoted $\mu(x,\de)=x\de$.  We must show that ${|X|\geq3p_3+p_1+1}$.  Since $\mu$ is faithful, and since the minimal congruences of $\B_n$ are generated by the pairs $(\ze,\al)$ and $(\ze,\be)$, we can fix elements $x_1,x_3\in X$ such that 
\begin{equation}\label{eq:x1x3}
x_1\ze \not= x_1\be \AND x_3\ze \not= x_3\al.
\end{equation}
We define the sets 
\[
Y_1 = x_1\ol P_1 = \set{x_1\ol\ve}{\ve\in P_1} \AND Y_3 = x_3\C_3\ol P_3 = \set{x_3\xi\ol\ve}{\xi\in\C_3,\ \ve\in P_3},
\]
noting that $|Y_1|=p_1$ and $|Y_3|=|\C_3\ol P_3|=3p_3$ by Lemmas \ref{lem:act3} and \ref{lem:Bsiodd}.  Thus, we can show that ${|X|\geq3p_3+p_1+1}$ by showing that 
\begin{equation}\label{eq:Y1Y3}
Y_1\cap Y_3=\es \AND X\sm(Y_1\cup Y_3)\not=\es.
\end{equation}
To prove the first assertion, suppose to the contrary that $x_1\ol\ve_1=x_3\xi\ol\ve_3$ for some $\ve_i\in P_i$ and $\xi\in\C_3$.  Then with $\th_1 = \ol\ve_3^*\xi^{-1}\al$ and $\th_2 = \ol\ve_3^*\xi^{-1}\ze$, we have
\[
\xi\ol\ve_3\th_1 = \al \AND \ol\ve_1\th_1 = \ol\ve_1\th_2 = \xi\ol\ve_3\th_2 = \ze.
\]
It follows that
\[
x_3\al = (x_3\xi\ol\ve_3)\th_1 = (x_1\ol\ve_1)\th_1 = x_1\ze = (x_1\ol\ve_1)\th_2 = (x_3\xi\ol\ve_3)\th_2 = x_3\ze,
\]
contradicting \eqref{eq:x1x3}.

This leaves us to prove the second assertion in \eqref{eq:Y1Y3}.  This is certainly true if $x_3\ze\not\in Y_1\cup Y_3$, so suppose otherwise.  

If $x_3\ze\in Y_3$, say with $x_3\ze = x_3\xi\ol\ve$ for $\xi\in\C_3$ and $\ve\in P_3$, then with $\th = \ol\ve^*\xi^{-1}\al$ we have $\ze\th=\ze$ and $\xi\ol\ve\th = \al$, and this gives $x_3\ze = (x_3\ze)\th = (x_3\xi\ol\ve)\th = x_3\al$, again contradicting \eqref{eq:x1x3}.

So we must instead have $x_3\ze\in Y_1$, say with $x_3\ze = x_1\ol\ve$ for $\ve\in P_1$, and it follows that
\[
x_3\ze = x_3\ze\ze = x_1\ol\ve\ze = x_1\ze.
\]
We complete the proof by showing that $x_3\al\not\in Y_1\cup Y_3$.  Indeed, if $x_3\al\in Y_1$, say with $x_3\al = x_1\ol\ve_1$ for $\ve_1\in P_1$, then 
\[
x_3\al = x_3\al\ze = x_1\ol\ve_1\ze = x_1\ze = x_3\ze,
\]
contradicting \eqref{eq:x1x3}.  But if $x_3\al\in Y_3$, say with $x_3\al = x_3\xi\ol\ve_3$ for $\xi\in\C_3$ and $\ve_3\in P_3$, then
\[
x_3\al = x_3\al \cdot \ol\ve_3^*\xi^{-1}\ze = x_3\xi\ol\ve_3 \cdot \ol\ve_3^*\xi^{-1}\ze = x_3\ze,
\]
contradicting \eqref{eq:x1x3}.
\epf

\subsection{Lower bound -- even case}\label{subsect:Blow_even}

We now assume $n=2k\geq4$ is even, and we fix $\ze,\al,\be,\ga,\de\in\B_n$, as in \eqref{eq:Bzabcd}.  In addition to the subgroup $\C_3 = \la (1,2,3)\ra = \{\id_n,(1,2,3),(1,3,2)\}$ of $\S_n$, we also fix $\C_2 = \la(1,2)\ra = \{\id_n,(1,2)\}$.

\begin{lemma}\label{lem:Bsieven}
Let $\si$ be a right congruence of $\B_n$, where $n\geq4$ is even.
\ben
\item \label{Bsieven1} If $\si$ separates $\{\ze,\al\}$, then it separates $\C_3\ol P_4$.
\item \label{Bsieven2} If $\si$ separates $\{\ze,\be\}$, then it separates $\ol P_0$.
\item \label{Bsieven3} If $\si$ separates $\{\ga,\de\}$, then it separates $\C_2\ol P_2$.
\een
\end{lemma}

\pf
\firstpfitem{\ref{Bsieven1}}  Suppose $(\xi_1\ol\ve_1,\xi_2\ol\ve_2)\in\si$ for some ${\xi_1,\xi_2\in\C_3}$ and $\ve_1,\ve_2\in P_4$, with $\xi_1\ol\ve_1 \not= \xi_2\ol\ve_2$, and write
\[
\xi_1\ol\ve_1 = \begin{partn}{7} 1&2&3&4&\multicolumn{3}{c}{}\\ \hhline{~|~|~|~|-|-|-} a_1&a_2&a_3&a_4&A_1&\cdots&A_{k-2}\end{partn}
\AND
\xi_2\ol\ve_2 = \begin{partn}{7} 1&2&3&4&\multicolumn{3}{c}{}\\ \hhline{~|~|~|~|-|-|-} b_1&b_2&b_3&b_4&B_1&\cdots&B_{k-2}\end{partn}.
\]
We must show that $(\ze,\al)\in\si$.  We now consider various cases; in each we define an element $\th\in\B_n$ for which $\{\xi_1\ol\ve_1\th,\xi_2\ol\ve_2\th\} = \{\ze,\al\}$.

\pfcase1  First, if
\bena\bmc2
\item \label{ab1} $\{a_1,a_2\}\cap\{b_1,b_4\} = \es$,
\item \label{ab2} $\{a_1,a_2\}\cap\{b_2,b_3\} = \es$,
\item \label{ab3} $\{a_3,a_4\}\cap\{b_1,b_4\} = \es$, or
\item \label{ab4} $\{a_3,a_4\}\cap\{b_2,b_3\} = \es$,
\emc\een
then we take 
\ref{ab1} $\th = \begin{partn}{2} a_1,a_2 & b_1,b_4 \\ \hhline{-|-} \multicolumn{2}{c}{} \end{partn}$, 
\ref{ab2} $\th = \begin{partn}{2} a_1,a_2 & b_2,b_3 \\ \hhline{-|-} \multicolumn{2}{c}{} \end{partn}$, 
\ref{ab3} $\th = \begin{partn}{2} a_3,a_4 & b_1,b_4 \\ \hhline{-|-} \multicolumn{2}{c}{} \end{partn}$, or 
\ref{ab4} $\th = \begin{partn}{2} a_3,a_4 & b_2,b_3 \\ \hhline{-|-} \multicolumn{2}{c}{} \end{partn}$, 
respectively.  

\pfcase2  Now suppose all four of 
\[
\{a_1,a_2\}\cap\{b_1,b_4\} \COMMA
\{a_1,a_2\}\cap\{b_2,b_3\} \COMMA
\{a_3,a_4\}\cap\{b_1,b_4\} \AND
\{a_3,a_4\}\cap\{b_2,b_3\} 
\]
are non-empty, noting then that this forces $\{a_1,a_2,a_3,a_4\} = \{b_1,b_2,b_3,b_4\}$.  In particular, this says that $\codom(\ve_1) = \codom(\ve_2)$, and it follows that $\ol\ve_1$ and $\ol\ve_2$ have the same transversals.  Since $\xi_1,\xi_2\in\C_3$, it follows that the ordered tuple $(b_1,b_2,b_3,b_4)$ is equal to either
\bena\bmc3 
\item \label{a1} $(a_1,a_2,a_3,a_4)$,
\item \label{a2} $ (a_2,a_3,a_1,a_4)$, or
\item \label{a3} $(a_3,a_1,a_2,a_4)$.
\emc\een
In cases \ref{a2} and \ref{a3} we take $\th = \begin{partn}{2} a_1,a_4 & a_2,a_3 \\ \hhline{-|-} \multicolumn{2}{c}{} \end{partn}$ or $\begin{partn}{2} a_1,a_2 & a_3,a_4 \\ \hhline{-|-} \multicolumn{2}{c}{} \end{partn}$, respectively.  So, finally consider case \ref{a1}, in which we have $a_i=b_i$ for $i=1,2,3,4$.  Since $\xi_1\ol\ve_1 \not= \xi_2\ol\ve_2$, we can assume without loss of generality that $A_1=\{u,v\}$, $B_1=\{u,x\}$ and $B_2=\{v,y\}$ for distinct $u,v,x,y\in\bn$, and we take $\th = \begin{partn}{4} a_1,u & a_2,v & a_3,y & a_4,x \\ \hhline{-|-|-|-} \multicolumn{4}{c}{} \end{partn}$.

\pfitem{\ref{Bsieven2}}  Suppose $(\ol\ve_1,\ol\ve_2)\in\si$ for distinct $\ve_1,\ve_2\in P_0$, and fix
\begin{equation}\label{eq:uvxy}
(u,v)\in\coker(\ve_1) \AND (u,x),(v,y)\in\coker(\ve_2), \WHERE u,v,x,y\in\bn \text{ are distinct.}
\end{equation}
Then with $\th = \begin{partn}{4} x&u&v&y \\ \hhline{~|~|~|~} 1&2&3&4\end{partn}$ we have $(\be,\ze) = (\ol\ve_1\th,\ol\ve_2\th) \in \si$.

\pfitem{\ref{Bsieven3}}  Suppose $(\xi_1\ol\ve_1,\xi_2\ol\ve_2)\in\si$ for some ${\xi_1,\xi_2\in\C_2}$ and $\ve_1,\ve_2\in P_2$, with $\xi_1\ol\ve_1 \not= \xi_2\ol\ve_2$, and write 
\[
\xi_1\ol\ve_1 = \begin{partn}{5} 1&2&\multicolumn{3}{c}{}\\ \hhline{~|~|-|-|-} a_1&a_2&A_1&\cdots&A_{k-1}\end{partn}
\AND
\xi_2\ol\ve_2 = \begin{partn}{5} 1&2&\multicolumn{3}{c}{}\\ \hhline{~|~|-|-|-} b_1&b_2&B_1&\cdots&B_{k-1}\end{partn}.
\]
If $a_1\not=b_1$ or $a_2\not=b_2$, then with $\th = \begin{partn}{2} a_1&b_1 \\ \hhline{~|~} 1&2\end{partn}$ or $\begin{partn}{2} a_2&b_2 \\ \hhline{~|~} 2&1\end{partn}$, respectively, we have ${(\ga,\de) = (\xi_1\ol\ve_1\th,\xi_2\ol\ve_2\th)\in\si}$.  Now suppose $a_1=b_1$ and $a_2=b_2$.  Since $\xi_1\ol\ve_1 \not= \xi_2\ol\ve_2$, we can fix $u,v,x,y\in\bn$ as in \eqref{eq:uvxy}, and with $\th = \begin{partn}{4} v&x&a_1,u&a_2,y \\ \hhline{~|~|-|-} 1&2&\multicolumn{2}{c}{} \end{partn}$ we have ${(\ga,\de) = (\xi_1\ol\ve_1\th,\xi_2\ol\ve_2\th)\in\si}$.
\epf

Here then is the final piece of the proof of Theorem \ref{thm:B}.

\begin{prop}\label{prop:Bloweven}
If $n\geq4$ is even, then $\deg(\B_n) \geq 3p_4 + 2p_2 + p_0 + 1$.
\end{prop}

\pf
Let $\mu:X\times \B_n\to X$ be a faithful action.  We must show that $|X|\geq3p_4 + 2p_2 + p_0 + 1$.  Since $\mu$ is faithful, and since the minimal congruences of $\B_n$ are generated by the pairs $(\ze,\al)$, $(\ze,\be)$ and $(\ga,\de)$, we can fix elements $x_0,x_2,x_4\in X$ such that 
\begin{equation}\label{eq:x0x2x4}
x_0\ze \not= x_0\be \COMMA x_2\ga \not= x_2\de \AND x_4\ze \not= x_4\al,
\end{equation}
where we again use shorthand notation for the action.
We define the sets 
\[
Y_0 = x_0\ol P_0 \COMMA  Y_2 = x_2\C_2\ol P_2 \AND Y_4 = x_4\C_3\ol P_4,
\]
noting that $|Y_0|=p_0$, $|Y_2|=2p_2$ and $|Y_4|=3p_4$, by Lemmas \ref{lem:act3} and \ref{lem:Bsieven}.  Thus, we can show that ${|X|\geq3p_4 + 2p_2 + p_0 + 1}$ by showing that 
\begin{equation}\label{eq:Y0Y2Y4X}
\text{$Y_0$, $Y_2$ and $Y_4$ are pairwise disjoint} \AND X\sm(Y_0\cup Y_2\cup Y_4)\not=\es.
\end{equation}
Before we do this, we first claim that
\begin{equation}\label{eq:Y2Y4ze}
(Y_2\cup Y_4) \cap \set{x\ze}{x\in X} = \es.
\end{equation}
To prove this, suppose to the contrary that $x\ze\in Y_2\cup Y_4$ for some $x\in X$.  

First consider the case that $x\ze\in Y_2$, and write $x\ze = x_2\xi\ol\ve$ where $\xi\in\C_2$ and $\ve\in P_2$.  Then
\[
x\ze = x\ze\ol\ve^*\xi^{-1} = x_2\xi\ol\ve\ol\ve^*\xi^{-1} = x_2\ga,
\]
and we then have $x_2\de = x_2\ga\de = x\ze\de = x\ze = x_2\ga$, contradicting \eqref{eq:x0x2x4}.

Now suppose instead that $x\ze\in Y_4$, and write $x\ze = x_4\xi\ol\ve$ where $\xi\in\C_3$ and $\ve\in P_4$.  This time we first note that $\ze\ol\ve^*\xi^{-1}$ is equal to one of 
\[
\custpartn{1,2,3,4,5,6,9,10}{1,2,3,4,5,6,9,10}{\uarc12\uarc34\uarc56\uarc9{10}\darc12\darc34\darc56\darc9{10}\udotted69\ddotted69}
\COMMA
\custpartn{1,2,3,4,5,6,9,10}{1,2,3,4,5,6,9,10}{\uarc12\uarc34\uarc56\uarc9{10}\darcx14{.6}\darcx23{.3}\darc56\darc9{10}\udotted69\ddotted69}
\OR
\custpartn{1,2,3,4,5,6,9,10}{1,2,3,4,5,6,9,10}{\uarc12\uarc34\uarc56\uarc9{10}\darcx13{.6}\darcx24{.3}\darc56\darc9{10}\udotted69\ddotted69},
\]
each of which has the form $\ze\xi_0$ for some $\xi_0\in\C_3$.  Then with $\pi = \custpartn{1,2,3,4,5,6,9,10}{1,2,3,4,5,6,9,10}{\uarc56\uarc9{10}\darc56\darc9{10}\stline11\stline22\stline33\stline44\udotted69\ddotted69}$, it follows that
\[
x\ze\xi_0 = x\ze\ol\ve^*\xi^{-1} = x_4\xi\ol\ve\ol\ve^*\xi^{-1} = x_4\pi.
\]
But then $x_4\ze = x_4\pi\ze = x\ze\xi_0\ze = x\ze = x\ze\xi_0\al = x_4\pi\al = x_4\al$, contradicting \eqref{eq:x0x2x4}.

Now that we have proved \eqref{eq:Y2Y4ze}, our next claim is that
\begin{equation}\label{eq:Y2Y4}
Y_2\cap Y_4 = \es.
\end{equation}
To prove this, suppose to the contrary that $x_2\xi_2\ol\ve_2 = x_4\xi_4\ol\ve_4$ for some $\xi_2\in\C_2$, $\xi_4\in\C_3$, $\ve_2\in P_2$ and $\ve_4\in P_4$.  Then
\[
x_2\ze = x_2\xi_2\ol\ve_2\ol\ve_4^*\xi_4^{-1}\al = x_4\xi_4\ol\ve_4\ol\ve_4^*\xi_4^{-1}\al = x_4\pi\al = x_4\al.
\]
The previous calculation is valid with $\al$ replaced by $\ze$, leading to $x_2\ze = x_4\ze$.  Combining these, it follows that $x_4\al = x_2\ze = x_4\ze$, contradicting \eqref{eq:x0x2x4}.

Now that we have proved \eqref{eq:Y2Y4}, our next claim is that
\begin{equation}\label{eq:Y0Y2Y4}
Y_0\cap(Y_2\cup Y_4) = \es.
\end{equation}
To prove this, suppose to the contrary that $x_0\ol\ve_0 = x_i\xi_i\ol\ve_i$, where $i=2$ or $4$, and where $\ve_i\in P_i$ and $\xi_i\in\C_2$ or $\C_3$ as appropriate.  We then have $x_0\ol\ve_0\ol\ve_i^*\xi_i^{-1}= x_i\xi_i\ol\ve_i\ol\ve_i^*\xi_i^{-1}$, and we note that
\[
\xi_i\ol\ve_i\ol\ve_i^*\xi_i^{-1} = \begin{cases}
\ga &\text{if $i = 2$}\\
\pi &\text{if $i = 4$,}
\end{cases}
\AND \ol\ve_0\ve_i^*\xi_i^{-1} = \ze,\ \be \text{ or } \be\de,
\]
but that $\ol\ve_0\ve_i^*\xi_i^{-1} = \be$ or $\be\de$ is only possible when $i=4$.  It follows that either
\[
x_2\ga = x_0\ze \COMMA 
x_4\pi = x_0\ze \COMMA
x_4\pi = x_0\be \OR
x_4\pi = x_0\be\de.
\]
Keeping in mind that $\ga=\ol\ga$ and $\pi=\ol\pi$, the first two options contradict \eqref{eq:Y2Y4ze}.  We now investigate the remaining two cases.
First, if $x_4\pi = x_0\be\de$, then
\[
x_4\al = x_4\pi\al = x_0\be\de\al = x_0\be\de\ze = x_4\pi\ze = x_4\ze,
\]
contradicting \eqref{eq:x0x2x4}.
So now suppose $x_4\pi = x_0\be$.  As noted above, this arises when $i=4$ (so ${x_0\ol\ve_0 = x_4\xi_4\ol\ve_4}$) and $\ol\ve_0\ol\ve_4^*\xi_4^{-1} = \be$.  Here we calculate
\[
x_4\ze = x_4\pi\ze = x_0\be\ze = x_0\ze = x_0\ol\ve_0\ol\ve_0^* = x_4\xi_4\ol\ve_4\ol\ve_0^* = x_4 (\ol\ve_0\ol\ve_4^*\xi_4^{-1})^* = x_4\be^* = x_4\al,
\]
again contradicting \eqref{eq:x0x2x4}.

Now that we have proved \eqref{eq:Y2Y4} and \eqref{eq:Y0Y2Y4}, it remains to prove the second assertion in \eqref{eq:Y0Y2Y4X}.  This is clearly the case if either $x_2\ze$ or $x_4\ze$ does not belong to $Y_0\cup Y_2\cup Y_4$.  Given \eqref{eq:Y2Y4ze}, the only other possibility is that $x_2\ze,x_4\ze\in Y_0$, and we now assume this is the case.  So for $i=2,4$ we have $x_i\ze = x_0\ol\ve_i$ for some $\ve_i\in P_0$.  It then follows that $x_i\ze = x_i\ze\ze = x_0\ol\ve_i\ze = x_0\ze$, so in fact
\[
x_4\ze = x_2\ze = x_0\ze.
\]
We complete the proof by showing that $x_4\al\not\in Y_0\cup Y_2\cup Y_4$.  To do so, suppose to the contrary that $x_4\al = x_i\xi\ol\ve$ for some $i\in\{0,2,4\}$, and some $\ve\in P_i$ and $\xi\in\C_2\cup\C_3$ as appropriate.  We then note that $\xi\ol\ve\ol\ve^*\xi^{-1} = \ze$, $\ga$ or $\pi$; in any case it follows that $\xi\ol\ve\ol\ve^*\xi^{-1}\ze = \ze$.  But then
\[
x_4\al = x_4\al\ol\ve^*\xi^{-1}\ze = x_i\xi\ol\ve\ol\ve^*\xi^{-1}\ze = x_i\ze = x_4\ze,
\]
contradicting \eqref{eq:x0x2x4}.
\epf

\section{Combinatorics: the explicit formulae}\label{sect:C}

In Section \ref{sect:B} we obtained an explicit formula for the transformation degree of the Brauer monoid~$\B_n$; see \eqref{eq:degBn}.
For the other diagram monoids $M=\P_n$, $\PB_n$, $\PP_n$, $\M_n$ and $\TL_n$, we showed in Sections~\ref{sect:P} and~\ref{sect:TL} that
\[
\deg(M)=\degrc(M) = 1+|Q| \AND \deg'(M) = |Q|
\]
for $n\geq2$ or $3$ (as appropriate), where 
\[
Q = Q(M) = \begin{cases}
P_0\cup P_1\cup P_2 &\text{if $M=\P_n$, $\PB_n$, $\PP_n$ or $\M_n$}\\
P_0\cup P_2\cup P_4 &\text{if $M = \TL_n$ for even $n$}\\
P_1\cup P_3 &\text{if $M = \TL_n$ for odd $n$.}
\end{cases}
\]
In this final section we obtain explicit formulae for $|Q|$.  These are valid for arbitrary $n\geq0$, but we note that some of the subsets $P_r$ involved in the above unions are empty for very small $n$.

Formulae for the sizes of $P_r = P_r(M) = \set{\ve\in P}{\rank(\ve)=r}$ are known (see for example \cite[Proposition~4.6]{DDE2021}), and of course summing these over appropriate $r$ yields formulae for $|Q|$.  The resulting expressions are somewhat unwieldy, however.  For example, in the case of the partition monoid $\P_n$, it follows from \cite[Proposition~4.6]{DDE2021} that
$|Q| = \sum_{r=0}^2\sum_{k=r}^n\binom kr S(n,k)$,
where $S(n,k)$ is a Stirling number of the second kind.
Fortunately, since the relevant values of $r$ are very small, an alternative combinatorial analysis of these cases is possible, and allows us to derive much simpler expressions for~$|Q|$ in terms of very fundamental number sequences, such as Bell, Catalan and Motzkin numbers.

Our strategy relies on the existence, for any $0\leq r\leq n$, of an injective mapping
\begin{equation}\label{eq:ve+}
P_r(\P_n) \to P_0(\P_{n+r}) : 
\ve = \begin{partn}{6} A_1&\cdots&A_r&C_1&\cdots&C_s\\ \hhline{~|~|~|-|-|-} A_1&\cdots&A_r&C_1&\cdots&C_s\end{partn} 
\mt \ve^+ = \begin{partn}{6} A_1\cup\{n+r\}&\cdots&A_r\cup\{n+1\}&C_1&\cdots&C_s\\ \hhline{-|-|-|-|-|-} A_1\cup\{n+r\}&\cdots&A_r\cup\{n+1\}&C_1&\cdots&C_s\end{partn}.
\end{equation}
Recall that we always assume $\min(A_1)<\cdots<\min(A_r)$ when using this tabular notation.  It follows that the map $\ve\mt\ve^+$ preserves planarity; it also maps (partial) Brauer diagrams to (partial) Brauer diagrams.  As an example, consider a typical projection $\ve\in P_3(\TL_n)$:
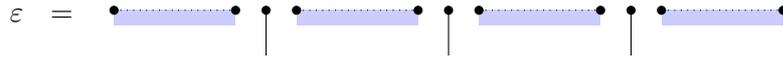
\[
\begin{tikzpicture}[scale=0.4]
\foreach \x in {1,2,3,4} {\ur{6*\x-5}{6*\x-1}{.5}};
\hstlines{6,12,18}
\foreach \x in {1,5,6,7,11,12,13,17,18,19,23} \uv\x;
\foreach \x in {24,25,26} \uvc\x{white};
\node[left] () at (0,2) {$\ve\phantom{^+} = $};
\end{tikzpicture}
\]
Here we have only drawn the top half of $\ve$, as the bottom is just the mirror image.  We then have
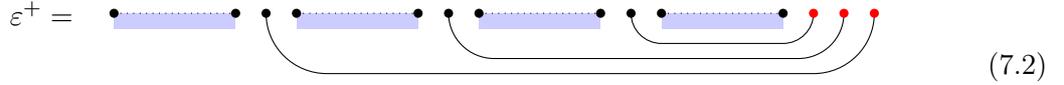
\begin{equation}\label{eq:ve+pic}
\begin{tikzpicture}[scale=0.4]
\foreach \x in {1,2,3,4} {\ur{6*\x-5}{6*\x-1}{.5}};
\foreach \x in {1,5,6,7,11,12,13,17,18,19,23} \uv\x;
\uarcx{18}{24}{1}
\uarcx{12}{25}{1.5}
\uarcx{6}{26}{2}
\foreach \x in {24,25,26} \uvc\x{red};
\node[left] () at (0,2) {$\ve^+ = $};
\end{tikzpicture}
\end{equation}
To simplify expressions in what follows, we will write 
\[
p_r(M) = |P_r(M)| \qquad\text{for any diagram monoid $M$, and any integer $r\geq0$.}
\]
We therefore have
\[
|Q(M)| = p_0(M) + p_1(M) + p_2(M) \qquad\text{if $M$ is one of $\P_n$, $\PB_n$, $\PP_n$ or $\M_n$.}
\]
Analogous statements hold for $M=\TL_n$.  Our ultimate goal, however, is to express each $|Q(M)|$ in terms of the relevant $p_0$ parameters, which themselves have the following simple forms:

\begin{prop}[{see \cite[Proposition 4.6]{DDE2021}}]\label{prop:proj}
For $n\geq0$, we have
\ben
\item \label{proj1} $p_0(\P_n) = B(n)$, the $n$th Bell number,
\item \label{proj2} $p_0(\PB_n) = I(n)$, the $n$th involution number,
\item \label{proj3} $p_0(\PP_n) = C(n)$, the $n$th Catalan number, 
\item \label{proj4} $p_0(\M_n) = M(n)$, the $n$th Motzkin number,
\item \label{proj5} $p_0(\TL_n) = C(n/2)$ for even $n$.  \qed
\een
\end{prop}

The ubiquity of the numbers occurring in Proposition \ref{prop:proj} is evidenced by their extremely low sequence numbers as  A000085, A000108, A000110 and A001006 on the OEIS \cite{OEIS}.  Note that~$I(n)$ is the number of involutions (i.e.~self-inverse permutations) of $\bn$.  These are given by the recurrence
\[
I(0)=I(1)=1 \AND I(n) = I(n-1) + (n-1)I(n-2) \qquad\text{for $n\geq2$.}
\]
We now proceed to give formulae for $|Q(M)|$, starting with $M=\P_n$.

\begin{prop}\label{prop:QPn}
For $n\geq0$ we have $\ds |Q(\P_n)| = \frac{B(n+2)-B(n+1)+B(n)}2$.
\end{prop}

\pf
Since $|Q(\P_n)| = p_0(\P_n) + p_1(\P_n) + p_2(\P_n)$, we can prove the result by showing that:
\ben
\item \label{p0Pn} $p_0(\P_n) = B(n)$, 
\item \label{p1Pn} $p_1(\P_n) = B(n+1)-B(n)$,
\item \label{p2Pn} $p_2(\P_n) = \frac{B(n+2)-3B(n+1)+B(n)}2$.
\een

\pfitem{\ref{p0Pn}}  This is part \ref{proj1} of Proposition \ref{prop:proj}.

\pfitem{\ref{p1Pn}}  The image of the map $P_1(\P_n) \to P_0(\P_{n+1}):\ve\mt\ve^+$ consists of all projections from $P_0(\P_{n+1})$ not containing the block $\{n+1\}$.  Since there are $p_0(\P_n)$ such `offending' projections, we have
\[
p_1(\P_n) = p_0(\P_{n+1}) - p_0(\P_n) = B(n+1) - B(n).
\]

\pfitem{\ref{p2Pn}}  This time we consider the map $P_2(\P_n) \to P_0(\P_{n+2}):\ve\mt\ve^+$, and calculate the size of its image, via the following steps.
\bena
\item \label{p2a} First we note that $P_0(\P_{n+2})$ contains $B(n+2)$ projections.
\item \label{p2b} From $B(n+2)$ we subtract $B(n+1)$, corresponding to the projections for which $n+1$ and $n+2$ belong to the same block. (These projections are not in the image of the $\ve\mt\ve^+$ map.)
\item \label{p2c} We then subtract a further $B(n+1)$, for the projections containing the block $\{n+1\}$.
\item \label{p2d} We then subtract $B(n+1)$ again, for those containing the block $\{n+2\}$.
\item \label{p2e} We must now add $B(n)$, for the projections containing \emph{both} blocks $\{n+1\}$ and $\{n+2\}$.
\een
At this point we are left with a set $\Si$ of $B(n+2) - 3B(n+1) + B(n)$ projections, but we must halve this total, as we have double-counted the image of the $\ve\mt\ve^+$ map.  Specifically, for 
${\ve = \begin{partn}{5} A_1&A_2&C_1&\cdots&C_s\\ \hhline{~|~|-|-|-} A_1&A_2&C_1&\cdots&C_s\end{partn} \in P_2(\P_n)}$, the set $\Si$ contains both
\[
\ve^+ = \begin{partn}{5} A_1\cup\{n+2\}&A_2\cup\{n+1\}&C_1&\cdots&C_s\\ \hhline{-|-|-|-|-} A_1\cup\{n+2\}&A_2\cup\{n+1\}&C_1&\cdots&C_s\end{partn}
\AND
\begin{partn}{5} A_1\cup\{n+1\}&A_2\cup\{n+2\}&C_1&\cdots&C_s\\ \hhline{-|-|-|-|-} A_1\cup\{n+1\}&A_2\cup\{n+2\}&C_1&\cdots&C_s\end{partn}.  \qedhere
\]
\epf

\begin{prop}\label{prop:QPBn}
For $n\geq0$ we have $\ds |Q(\PB_n)| = \frac{I(n+2)}2$.
\end{prop}

\pf
This time we show that:
\[
p_0(\PB_n) = I(n) \COMMA
p_1(\PB_n) = I(n+1)-I(n) \AND
p_2(\PB_n) = \tfrac{I(n+2)-2I(n+1)}2.
\]
The first two items are treated identically to the proof of Proposition \ref{prop:QPn}.  The third is \emph{almost} identical.  We follow steps \ref{p2a}--\ref{p2e}, with $I(k)$ in place of $B(k)$, with the only significant difference being in step~\ref{p2b}.  Here we subtract $I(n)$, rather than $I(n+1)$, as there are $I(n)$ projections in~$P_0(\PB_{n+2})$ for which $n+1$ and $n+2$ belong to the same block (as blocks have size $\leq2$).
\epf

\begin{prop}\label{prop:QPPn}
For $n\geq0$ we have $\ds |Q(\PP_n)| = C(n+2) -2C(n+1) + C(n)$.
\end{prop}

\pf
This time the result follows from:
\[
p_0(\PP_n) = C(n) \COMMa
p_1(\PP_n) = C(n+1)-C(n) \ANd
p_2(\PP_n) = C(n+2)-3C(n+1)+C(n).
\]
The proof is essentially the same as that of Proposition \ref{prop:QPn}, except that we do not need to divide by~$2$ in the third calculation, as planarity ensures there is no double-counting.
\epf

\begin{prop}\label{prop:QMn}
For $n\geq0$ we have $\ds |Q(\M_n)| = M(n+2) - M(n+1)$.
\end{prop}

\pf
This time we have
\[
p_0(\M_n) = M(n) \COMMa
p_1(\M_n) = M(n+1)-M(n) \ANd
p_2(\M_n) = M(n+2)-2M(n+1),
\]
as in Proposition \ref{prop:QPBn}, but with no double-counting.
\epf

\begin{prop}\label{prop:QTLn}
For $n\geq0$ we have
\[
\ds |Q(\TL_n)| = \begin{cases}
C(k+1) - C(k) &\text{if $n=2k-1$ is odd}\\
C(k+2) -2C(k+1) + C(k) &\text{if $n=2k$ is even.}
\end{cases}
\]
\end{prop}

\pf
Given Proposition \ref{prop:QPPn} and the isomorphism $\TL_{2k}\cong\PP_k$, we need only consider the case that $n=2k-1$ is odd.  Here we have $|Q(\TL_n)| = p_1(\TL_n) + p_3(\TL_n)$, and we claim that
\[
p_1(\TL_n) = C(k) \AND p_3(\TL_n) = C(k+1) - 2C(k).
\]
Indeed, the first holds because the map $P_1(\TL_n) \to P_0(\TL_{n+1}):\ve\mt\ve^+$ is a bijection.  (No Temperley--Lieb partition contains the block $\{n+1\}$.)

For the second, we note that the image of the map $P_3(\TL_n) \to P_0(\TL_{n+3}):\ve\mt\ve^+$ consists of all projections from $P_0(\TL_{n+3})$ not containing the blocks $\{n+1,n+2\}$ or $\{n+2,n+3\}$.  (See~\eqref{eq:ve+pic}, and note that no Temperley--Lieb partition contains the block $\{n+1,n+3\}$.)  We have $p_0(\TL_{n+3}) = C(k+1)$, and there are $p_0(\TL_{n+1}) = C(k)$ of each of the above two types of `offending' projections, with no overlap between them.
\epf

The values of $|Q(M)|$ in Propositions \ref{prop:QPn}--\ref{prop:QTLn} appear as sequences A000245, A001475, A002026, A026012 and A087649 on the OEIS \cite{OEIS}.

\begin{rem}
Denoting the \emph{forward difference} of a sequence $s(n)$ by ${\partial s(n) = s(n+1) - s(n)}$, we see that
\[
|Q(\PP_n)| = \partial^2 C(n) \COMMa
|Q(\M_n)| = \partial M(n+1) \ANd
|Q(\TL_n)| = \begin{cases}
\partial C(k) &\text{for $n=2k-1$}\\
\partial^2 C(k) &\text{for $n=2k$,}
\end{cases}
\]
where $C(n)$ and $M(n)$ denote Catalan and Motzkin numbers.
\end{rem}

\footnotesize
\def\bibspacing{-1.1pt}
\bibliography{biblio}

\begin{thebibliography}{10}

\bibitem{OEIS}
The on-line encyclopedia of integer sequences.
\newblock Published electronically at {\tt http://oeis.org/}.

\bibitem{AMSV2009}
J.~Almeida, S.~Margolis, B.~Steinberg, and M.~Volkov.
\newblock Representation theory of finite semigroups, semigroup radicals and
  formal language theory.
\newblock {\em Trans. Amer. Math. Soc.}, 361(3):1429--1461, 2009.

\bibitem{ACMT2023}
M.~Anagnostopoulou-Merkouri, R.~Cirpons, J.~D. Mitchell, and M.~Tsalakou.
\newblock Computing finite index congruences of finitely presented semigroups
  and monoids.
\newblock {\em Math. Comp.}, to appear, {\tt arXiv:2302.06295}.

\bibitem{Arbib1968}
M.~A. Arbib, editor.
\newblock {\em Algebraic theory of machines, languages, and semigroups}.
\newblock Academic Press, New York-London, 1968.
\newblock With a major contribution by Kenneth Krohn and John L. Rhodes.

\bibitem{BH2014}
G.~Benkart and T.~Halverson.
\newblock Motzkin algebras.
\newblock {\em European J. Combin.}, 36:473--502, 2014.

\bibitem{BT1991}
G.~Bijev and K.~Todorov.
\newblock On the representation of abstract semigroups by transformation
  semigroups: computer investigations.
\newblock {\em Semigroup Forum}, 43(2):253--256, 1991.

\bibitem{Brauer1937}
R.~Brauer.
\newblock On algebras which are connected with the semisimple continuous
  groups.
\newblock {\em Ann. of Math. (2)}, 38(4):857--872, 1937.

\bibitem{BSS2017}
J.~R. Britnell, N.~Saunders, and T.~Skyner.
\newblock On exceptional groups of order {$p^5$}.
\newblock {\em J. Pure Appl. Algebra}, 221(11):2647--2665, 2017.

\bibitem{CEFMPQ2024}
P.~J. Cameron, J.~East, D.~FitzGerald, J.~D. Mitchell, L.~Pebody, and
  T.~Quinn-Gregson.
\newblock Minimum degrees of finite rectangular bands, null semigroups, and
  variants of full transformation semigroups.
\newblock {\em Comb. Theory}, 3(3):Paper No. 16, 48 pp., 2023.

\bibitem{Clifford1942}
A.~H. Clifford.
\newblock Matrix representations of completely simple semigroups.
\newblock {\em Amer. J. Math.}, 64:327--342, 1942.

\bibitem{CPbook}
A.~H. Clifford and G.~B. Preston.
\newblock {\em The algebraic theory of semigroups. {V}ol. {I}}.
\newblock Mathematical Surveys, No. 7. American Mathematical Society,
  Providence, R.I., 1961.

\bibitem{DE2018}
I.~Dolinka and J.~East.
\newblock Twisted {B}rauer monoids.
\newblock {\em Proc. Roy. Soc. Edinburgh Sect. A}, 148(4):731--750, 2018.

\bibitem{DDE2021}
I.~Dolinka, I.~{\DJ}ur{\dj}ev, and J.~East.
\newblock Sandwich semigroups in diagram categories.
\newblock {\em Internat. J. Algebra Comput.}, 31(7):1339--1404, 2021.

\bibitem{Easdown1987}
D.~Easdown.
\newblock The minimal faithful degree of a fundamental inverse semigroup.
\newblock {\em Bull. Austral. Math. Soc.}, 35(3):373--378, 1987.

\bibitem{Easdown1988}
D.~Easdown.
\newblock The minimal faithful degree of a semilattice of groups.
\newblock {\em J. Austral. Math. Soc. Ser. A}, 45(3):341--350, 1988.

\bibitem{Easdown1992}
D.~Easdown.
\newblock Minimal faithful permutation and transformation representations of
  groups and semigroups.
\newblock In {\em Proceedings of the {I}nternational {C}onference on {A}lgebra,
  {P}art 3 ({N}ovosibirsk, 1989)}, volume 131 of {\em Contemp. Math.}, pages
  75--84. Amer. Math. Soc., Providence, RI, 1992.

\bibitem{EH2016}
D.~Easdown and M.~Hendriksen.
\newblock Minimal permutation representations of semidirect products of groups.
\newblock {\em J. Group Theory}, 19(6):1017--1048, 2016.

\bibitem{EP1988}
D.~Easdown and C.~E. Praeger.
\newblock On minimal faithful permutation representations of finite groups.
\newblock {\em Bull. Austral. Math. Soc.}, 38(2):207--220, 1988.

\bibitem{EEM2017}
J.~East, A.~Egri-Nagy, and J.~D. Mitchell.
\newblock Enumerating transformation semigroups.
\newblock {\em Semigroup Forum}, 95(1):109--125, 2017.

\bibitem{EEMP2019}
J.~East, A.~Egri-Nagy, J.~D. Mitchell, and Y.~P\'{e}resse.
\newblock Computing finite semigroups.
\newblock {\em J. Symbolic Comput.}, 92:110--155, 2019.

\bibitem{EG2017}
J.~East and R.~D. Gray.
\newblock Diagram monoids and {G}raham--{H}oughton graphs: {I}dempotents and
  generating sets of ideals.
\newblock {\em J. Combin. Theory Ser. A}, 146:63--128, 2017.

\bibitem{EG2021}
J.~East and R.~D. Gray.
\newblock Ehresmann theory and partition monoids.
\newblock {\em J. Algebra}, 579:318--352, 2021.

\bibitem{EMRT2018}
J.~East, J.~D. Mitchell, N.~Ru\v{s}kuc, and M.~Torpey.
\newblock Congruence lattices of finite diagram monoids.
\newblock {\em Adv. Math.}, 333:931--1003, 2018.

\bibitem{EST2010}
B.~Elias, L.~Silberman, and R.~Takloo-Bighash.
\newblock Minimal permutation representations of nilpotent groups.
\newblock {\em Experiment. Math.}, 19(1):121--128, 2010.

\bibitem{FP1997}
V.~Froidure and J.-E. Pin.
\newblock Algorithms for computing finite semigroups.
\newblock In {\em Foundations of computational mathematics ({R}io de {J}aneiro,
  1997)}, pages 112--126. Springer, Berlin, 1997.

\bibitem{GAP}
The GAP~Group.
\newblock {\em {GAP -- Groups, Algorithms, and Programming, Version 4.14.0}},
  2024, \url{https://www.gap-system.org}.

\bibitem{GM2008}
R.~Gray and J.~D. Mitchell.
\newblock Largest subsemigroups of the full transformation monoid.
\newblock {\em Discrete Math.}, 308(20):4801--4810, 2008.

\bibitem{HR2005}
T.~Halverson and A.~Ram.
\newblock Partition algebras.
\newblock {\em European J. Combin.}, 26(6):869--921, 2005.

\bibitem{Hoehnke1963}
H.-J. Hoehnke.
\newblock Zur {S}trukturtheorie der {H}albgruppen.
\newblock {\em Math. Nachr.}, 26:1--13, 1963.

\bibitem{Holt2010}
D.~F. Holt.
\newblock Enumerating subgroups of the symmetric group.
\newblock In {\em Computational group theory and the theory of groups, {II}},
  volume 511 of {\em Contemp. Math.}, pages 33--37. Amer. Math. Soc.,
  Providence, RI, 2010.

\bibitem{HW2002}
D.~F. Holt and J.~Walton.
\newblock Representing the quotient groups of a finite permutation group.
\newblock {\em J. Algebra}, 248(1):307--333, 2002.

\bibitem{Howie1995}
J.~M. Howie.
\newblock {\em Fundamentals of semigroup theory}, volume~12 of {\em London
  Mathematical Society Monographs. New Series}.
\newblock The Clarendon Press, Oxford University Press, New York, 1995.
\newblock Oxford Science Publications.

\bibitem{Imaoka1980}
T.~Imaoka.
\newblock On fundamental regular {$\ast $}-semigroups.
\newblock {\em Mem. Fac. Sci. Shimane Univ.}, 14:19--23, 1980.

\bibitem{Johnson1971}
D.~L. Johnson.
\newblock Minimal permutation representations of finite groups.
\newblock {\em Amer. J. Math.}, 93:857--866, 1971.

\bibitem{Jones1994_2}
V.~F.~R. Jones.
\newblock The {P}otts model and the symmetric group.
\newblock In {\em Subfactors ({K}yuzeso, 1993)}, pages 259--267. World Sci.
  Publ., River Edge, NJ, 1994.

\bibitem{KKM2000}
M.~Kilp, U.~Knauer, and A.~V. Mikhalev.
\newblock {\em Monoids, acts and categories}, volume~29 of {\em De Gruyter
  Expositions in Mathematics}.
\newblock Walter de Gruyter \& Co., Berlin, 2000.

\bibitem{Koenig2008}
S.~Koenig.
\newblock A panorama of diagram algebras.
\newblock In {\em Trends in representation theory of algebras and related
  topics}, EMS Ser. Congr. Rep., pages 491--540. Eur. Math. Soc., Z\"urich,
  2008.

\bibitem{KP1989}
L.~G. Kov\'acs and C.~E. Praeger.
\newblock Finite permutation groups with large abelian quotients.
\newblock {\em Pacific J. Math.}, 136(2):283--292, 1989.

\bibitem{KP2000}
L.~G. Kov\'{a}cs and C.~E. Praeger.
\newblock On minimal faithful permutation representations of finite groups.
\newblock {\em Bull. Austral. Math. Soc.}, 62(2):311--317, 2000.

\bibitem{LM1990}
G.~Lallement and R.~McFadden.
\newblock On the determination of {G}reen's relations in finite transformation
  semigroups.
\newblock {\em J. Symbolic Comput.}, 10(5):481--498, 1990.

\bibitem{Lawson1991}
M.~V. Lawson.
\newblock Semigroups and ordered categories. {I}. {T}he reduced case.
\newblock {\em J. Algebra}, 141(2):422--462, 1991.

\bibitem{LPRR1998}
S.~A. Linton, G.~Pfeiffer, E.~F. Robertson, and N.~Ru\v{s}kuc.
\newblock Groups and actions in transformation semigroups.
\newblock {\em Math. Z.}, 228(3):435--450, 1998.

\bibitem{LPRR2002}
S.~A. Linton, G.~Pfeiffer, E.~F. Robertson, and N.~Ru\v{s}kuc.
\newblock Computing transformation semigroups.
\newblock {\em J. Symbolic Comput.}, 33(2):145--162, 2002.

\bibitem{MS2023}
S.~Margolis and B.~Steinberg.
\newblock On the minimal faithful degree of {R}hodes semisimple semigroups.
\newblock {\em J. Algebra}, 633:788--813, 2023.

\bibitem{Martin1994}
P.~Martin.
\newblock Temperley-{L}ieb algebras for nonplanar statistical mechanics---the
  partition algebra construction.
\newblock {\em J. Knot Theory Ramifications}, 3(1):51--82, 1994.

\bibitem{Martin2008}
P.~Martin.
\newblock On diagram categories, representation theory and statistical
  mechanics.
\newblock In {\em Noncommutative rings, group rings, diagram algebras and their
  applications}, volume 456 of {\em Contemp. Math.}, pages 99--136. Amer. Math.
  Soc., Providence, RI, 2008.

\bibitem{MM2014}
P.~Martin and V.~Mazorchuk.
\newblock On the representation theory of partial {B}rauer algebras.
\newblock {\em Q. J. Math.}, 65(1):225--247, 2014.

\bibitem{Semigroups}
J.~D. Mitchell et~al.
\newblock {\em Semigroups - GAP package, Version 5.5.0}, Feb 2025,
  \url{http://dx.doi.org/10.5281/zenodo.592893}.

\bibitem{Munn1957}
W.~D. Munn.
\newblock Matrix representations of semigroups.
\newblock {\em Proc. Cambridge Philos. Soc.}, 53:5--12, 1957.

\bibitem{Munn1964}
W.~D. Munn.
\newblock Matrix representations of inverse semigroups.
\newblock {\em Proc. London Math. Soc. (3)}, 14:165--181, 1964.

\bibitem{Munn1970}
W.~D. Munn.
\newblock Fundamental inverse semigroups.
\newblock {\em Quart. J. Math. Oxford Ser. (2)}, 21:157--170, 1970.

\bibitem{NS1978}
T.~E. Nordahl and H.~E. Scheiblich.
\newblock Regular {$\ast $}-semigroups.
\newblock {\em Semigroup Forum}, 16(3):369--377, 1978.

\bibitem{OPU2024}
E.~A. O'Brien, S.~K. Prajapati, and A.~Udeep.
\newblock Minimal degrees for faithful permutation representations of groups of
  order {$p^6$} where {$p$} is an odd prime.
\newblock {\em J. Algebraic Combin.}, 60(2):319--388, 2024.

\bibitem{Petrich1985}
M.~Petrich.
\newblock Certain varieties of completely regular {$^\ast$}-semigroups.
\newblock {\em Boll. Un. Mat. Ital. B (6)}, 4(2):343--370, 1985.

\bibitem{Ponizovskii1956}
I.~S. Ponizovski\u{\i}.
\newblock On matrix representations of associative systems.
\newblock {\em Mat. Sb. (N.S.)}, 38(80):241--260, 1956.

\bibitem{RS2009}
J.~Rhodes and B.~Steinberg.
\newblock {\em The {$q$}-theory of finite semigroups}.
\newblock Springer Monographs in Mathematics. Springer, New York, 2009.

\bibitem{Saunders2010}
N.~Saunders.
\newblock The minimal degree for a class of finite complex reflection groups.
\newblock {\em J. Algebra}, 323(3):561--573, 2010.

\bibitem{Saunders2014}
N.~Saunders.
\newblock Minimal faithful permutation degrees for irreducible {C}oxeter groups
  and binary polyhedral groups.
\newblock {\em J. Group Theory}, 17(5):805--832, 2014.

\bibitem{Schein1988}
B.~M. Schein.
\newblock The minimal degree of noble inverse semigroups.
\newblock In {\em Contributions to general algebra, 6}, pages 247--252.
  H\"{o}lder-Pichler-Tempsky, Vienna, 1988.

\bibitem{Schein1992}
B.~M. Schein.
\newblock The minimal degree of a finite inverse semigroup.
\newblock {\em Trans. Amer. Math. Soc.}, 333(2):877--888, 1992.

\bibitem{Slover1965}
R.~E. Slover.
\newblock Representations of a semigroup.
\newblock {\em Trans. Amer. Math. Soc.}, 120:417--427, 1965.

\bibitem{Steinberg2006}
B.~Steinberg.
\newblock M\"obius functions and semigroup representation theory.
\newblock {\em J. Combin. Theory Ser. A}, 113(5):866--881, 2006.

\bibitem{Steinberg2008}
B.~Steinberg.
\newblock M\"obius functions and semigroup representation theory. {II}.
  {C}haracter formulas and multiplicities.
\newblock {\em Adv. Math.}, 217(4):1521--1557, 2008.

\bibitem{Steinberg2016}
B.~Steinberg.
\newblock {\em Representation theory of finite monoids}.
\newblock Universitext. Springer, Cham, 2016.

\bibitem{Stoll1944}
R.~R. Stoll.
\newblock Representations of finite simple semigroups.
\newblock {\em Duke Math. J.}, 11:251--265, 1944.

\bibitem{TL1971}
H.~N.~V. Temperley and E.~H. Lieb.
\newblock Relations between the ``percolation'' and ``colouring'' problem and
  other graph-theoretical problems associated with regular planar lattices:
  some exact results for the ``percolation'' problem.
\newblock {\em Proc. Roy. Soc. London Ser. A}, 322(1549):251--280, 1971.

\bibitem{Tully1961}
E.~J. Tully, Jr.
\newblock Representation of a semigroup by transformations acting transitively
  on a set.
\newblock {\em Amer. J. Math.}, 83:533--541, 1961.

\bibitem{Wallace1963}
A.~D. Wallace.
\newblock Relative ideals in semigroups. {II}. {T}he relations of {G}reen.
\newblock {\em Acta Math. Acad. Sci. Hungar.}, 14:137--148, 1963.

\bibitem{Wright1975}
D.~Wright.
\newblock Degrees of minimal embeddings for some direct products.
\newblock {\em Amer. J. Math.}, 97(4):897--903, 1975.

\end{thebibliography}
\bibliographystyle{abbrv}
\end{document}